# W. B. Vasantha Kandasamy

# *Groupoids and Smarandache Groupoids*

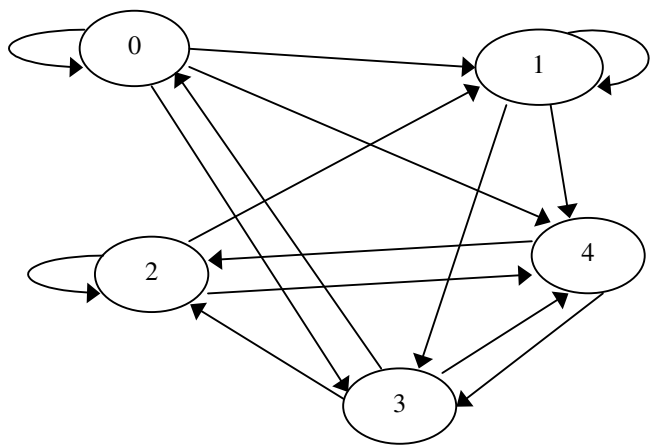

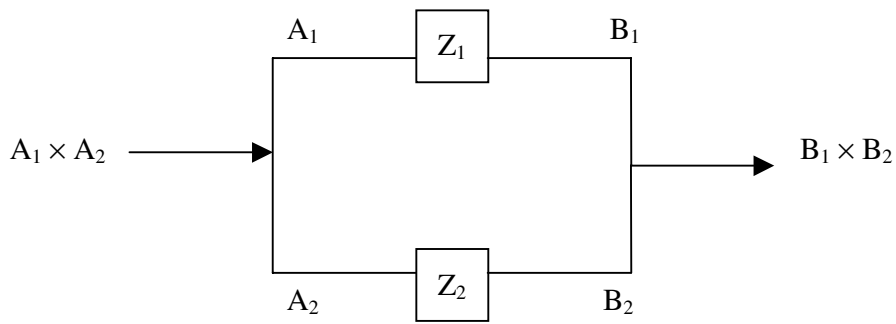

2 0 0 2

*W. B. Vasantha Kandasamy*
*Department of Mathematics*
*Indian Institute of Technology*
*Madras, Chennai – 600036, India*

# Groupoids and Smarandache Groupoids

| * | e | 0 | 1 | 2 | 3 | 4 | 5 |
|---|---|---|---|---|---|---|---|
| e | e | 0 | 1 | 2 | 3 | 4 | 5 |
| 0 | 0 | e | 3 | 0 | 3 | 0 | 3 |
| 1 | 1 | 5 | e | 5 | 2 | 5 | 2 |
| 2 | 2 | 4 | 1 | e | 1 | 4 | 1 |
| 3 | 3 | 3 | 0 | 3 | e | 3 | 0 |
| 4 | 4 | 2 | 5 | 2 | 5 | e | 5 |
| 5 | 5 | 1 | 4 | 1 | 4 | 1 | e |

**2002**

**Definition:**
**Generally, in any human field, a** *Smarandache Structure* **on a set A means a weak structure W on A such that there exists a proper subset B ⊂ A which is embedded with a stronger structure S.**

**These types of structures occur in our everyday's life, that's why we study them in this book.**

**Thus, as a particular case:**
**A** *Smarandache Groupoid* **is a groupoid G which has a proper subset S ⊂ G such that S under the operation of G is a semigroup.**

# CONTENTS

# PREFACE

The study of Smarandache Algebraic Structure was initiated in the year 1998 by Raul Padilla following a paper written by Florentin Smarandache called "Special Algebraic Structures". In his research, Padilla treated the Smarandache algebraic structures mainly with associative binary operation. Since then the subject has been pursued by a growing number of researchers and now it would be better if one gets a coherent account of the basic and main results in these algebraic structures. This book aims to give a systematic development of the basic non-associative algebraic structures viz. Smarandache groupoids. Smarandache groupoids exhibits simultaneously the properties of a semigroup and a groupoid. Such a combined study of an associative and a non associative structure has not been so far carried out. Except for the introduction of smarandacheian notions by Prof. Florentin Smarandache such types of studies would have been completely absent in the mathematical world.

Thus, Smarandache groupoids, which are groupoids with a proper subset, which is a semigroup, has several interesting properties, which are defined and studied in this book in a sequential way. This book assumes that the reader should have a good background of algebraic structures like semigroup, group etc. and a good foundation in number theory.

In Chapter 1 we just recall the basic notations and some important definitions used in this book. In Chapter 2 almost all concepts, most of them new have been introduced to groupoids in general. Since the study of groupoids and books on groupoids is meager, we in Chapter 3 introduce four new classes of groupoids using the set of modulo integers $Z_n$, $n \geq 3$ and $n < \infty$. This chapter is mainly introduced to lessen the non abstractness of this structure. In this chapter, several number theoretic techniques are used.

Chapter 4 starts with the definition of Smarandache groupoids. All properties introduced in groupoids are studied in the case of Smarandache groupoids. Several problems and examples are given in each section to make the concept easy. In Chapter 5 conditions for the new classes of groupoids built using $Z_n$ to contain Smarandache groupoids are obtained. Chapter 6 gives the application of Smarandache groupoids to semi automaton and automaton, that is to finite machines. The final chapter on research problems is the major attraction of the book as we give several open problems about groupoids. Any researcher on algebra will find them interesting and absorbing.

We have attempted to make this book a self contained one provided a reader has a sufficient background knowledge in algebra. Thus, this book will be the first one in Smarandache algebraic structures to deal with non associative operations.

I deeply acknowledge Dr. Minh Perez, because of whose support and constant encouragement this book was possible.

# CHAPTER ONE
# PRELIMINARY NOTIONS

In this chapter, we give some basic notion and preliminary concepts used in this book so as to make this book self contained. The study of groupoids is very rare and meager; the only reason the author is able to attribute to this is that it may be due to the fact that there is no natural way by which groupoids can be constructed. This book aim is two fold, firstly to construct new classes of groupoids using finite integers and define in these new classes many properties which have not been studied yet. Secondly, to define Smarandache groupoids and introduce the newly defined properties in groupoids to Smarandache groupoids. In this chapter, we recall some basic properties of integers, groupoids, Smarandache groupoids and loops.

## 1.1 Integers

We start this chapter with a brief discussion on the set of both finite and infinite integers. We mainly enumerate the properties, which will be used in this book. As concerned with notations, the familiar symbols a > b, a ≥ b, |a|, a / b, b ∤ a occur with their usual meaning.

**DEFINITION:** *The positive integer c is said to be the greatest divisor of a and b if*

1. *c is a divisor of a and of b*
2. *Any divisor of a and b is a divisor of c.*

**DEFINITION:** *The integers a and b are relatively prime if (a, b)= 1 and there exists integers m and n such that ma + nb = 1.*

**DEFINITION:** *The integer p > 1 is a prime if its only divisor are ± 1 and ± p.*

**DEFINITION:** *The least common multiple of two positive integers a and b is defined to be the smallest positive integer that is divisible by a and b and it is denoted by l.c.m (a, b) or [a, b].*

**Notation:**

1. $Z^+$ is the set of positive integers.

2. $Z^+ \cup \{0\}$ is the set of positive integers with zero.

3. $Z = Z^+ \cup Z^- \cup \{0\}$ is the set of integers where $Z^-$ is the set of negative integers.

4. $Q^+$ is the set of positive rationals.

5. $Q^+ \cup \{0\}$ is the set of positive rationals with zero.

6. $Q = Q^+ \cup Q^- \cup \{0\}$, is the set of rationals where $Q^-$ is the set of negative rationals.

*Similarly $R^+$ is the set of positive reals, $R^+ \cup (0)$ is the set of positive reals with zero and the set of reals $R = R^+ \cup R^- \cup \{0\}$ where $R^-$ is the set of negative reals.*

Clearly, $Z^+ \subset Q^+ \subset R^+$ and $Z \subset Q \subset R$, where '$\subset$' denotes the containment that is ' contained ' relation.

$Z_n = \{0, 1, 2, \ldots , n-1\}$ be the set of integers under multiplication or under addition modulo n. For examples $Z_2 = \{0, 1\}$. $1 + 1 \equiv 0 \pmod 2$, $1.1 \equiv 1 \pmod 2$. $Z_9 = \{0, 1, 2, \ldots , 8\}$, $3 + 6 \equiv 0 \pmod 9$, $3.3 \equiv 0 \pmod 9$, $2.8 \equiv 7 \pmod 9$, $6.2 \equiv 3 \pmod 9$.

This notation will be used and $Z_n$ will denote the set of finite integers modulo n.

## 1.2 Groupoids

In this section we recall the definition of groupoids and give some examples. Problems are given at the end of this section to make the reader familiar with the concept of groupoids.

**DEFINITION:** *Given an arbitrary set P a mapping of $P \times P$ into P is called a binary operation on P. Given such a mapping $\sigma: P \times P \to P$ we use it to define a product $*$ in P by declaring $a * b = c$ if $\sigma(a, b) = c$.*

**DEFINITION:** *A non empty set of elements G is said to form a groupoid if in G is defined a binary operation called the product denoted by $*$ such that $a * b \in G$ for all $a, b \in G$.*

It is important to mention here that the binary operation $*$ defined on the set G need not be associative that is $(a * b) * c \neq a * (b * c)$ in general for all $a, b, c \in G$, so we can say the groupoid (G, $*$) is a set on which is defined a non associative binary operation which is closed on G.

**DEFINITION:** *A groupoid G is said to be a commutative groupoid if for every a, b ∈ G we have a ∗ b = b ∗ a.*

**DEFINITION:** *A groupoid G is said to have an identity element e in G if a ∗ e = e ∗ a = a for all a ∈ G.*

We call the order of the groupoid G to be the number of distinct elements in it denoted by o(G) or |G|. If the number of elements in G is finite we say the groupoid G is of finite order or a finite groupoid otherwise we say G is an infinite groupoid.

*Example 1.2.1:* Let $G = \{a_1, a_2, a_3, a_4, a_0\}$. Define ∗ on G given by the following table:

| ∗ | $a_0$ | $a_1$ | $a_2$ | $a_3$ | $a_4$ |
|---|---|---|---|---|---|
| $a_0$ | $a_0$ | $a_4$ | $a_3$ | $a_2$ | $a_1$ |
| $a_1$ | $a_1$ | $a_0$ | $a_4$ | $a_3$ | $a_2$ |
| $a_2$ | $a_2$ | $a_1$ | $a_0$ | $a_4$ | $a_3$ |
| $a_3$ | $a_3$ | $a_2$ | $a_1$ | $a_0$ | $a_4$ |
| $a_4$ | $a_4$ | $a_3$ | $a_2$ | $a_1$ | $a_0$ |

Clearly (G, ∗) is a non commutative groupoid and does not contain an identity. The order of this groupoid is 5.

*Example 1.2.2:* Let (S, ∗) be a groupoid with 3 elements given by the following table:

| ∗ | $x_1$ | $x_2$ | $x_3$ |
|---|---|---|---|
| $x_1$ | $x_1$ | $x_3$ | $x_2$ |
| $x_2$ | $x_2$ | $x_1$ | $x_3$ |
| $x_3$ | $x_3$ | $x_2$ | $x_1$ |

This is a groupoid of order 3, which is non associative and non commutative.

*Example 1.2.3:* Consider the groupoid (P, x) where $P = \{p_0, p_1, p_2, p_3\}$ given by the following table:

| × | $p_0$ | $p_1$ | $p_2$ | $p_3$ |
|---|---|---|---|---|
| $p_0$ | $p_0$ | $p_2$ | $p_0$ | $p_2$ |
| $p_1$ | $p_3$ | $p_1$ | $p_3$ | $p_1$ |
| $p_2$ | $p_2$ | $p_0$ | $p_2$ | $p_0$ |
| $p_3$ | $p_1$ | $p_3$ | $p_1$ | $p_3$ |

This is a groupoid of order 4.

*Example 1.2.4:* Let Z be the set of integers define an operation '−' on Z that is usual subtraction; (Z, −) is a groupoid. This groupoid is of infinite order and is both non commutative and non associative.

**DEFINITION:** *Let (G, ∗) be a groupoid a proper subset H ⊂ G is a subgroupoid if (H, ∗) is itself a groupoid.*

***Example 1.2.5:*** Let R be the reals (R, −) is a groupoid where ' − ' is the usual subtraction on R.

Now $Z \subset R$ is a subgroupoid, as (Z, −) is a groupoid.

***Example 1.2.6:*** Let G be a groupoid given by the following table:

| * | $a_1$ | $a_2$ | $a_3$ | $a_4$ |
|---|---|---|---|---|
| $a_1$ | $a_1$ | $a_3$ | $a_1$ | $a_3$ |
| $a_2$ | $a_4$ | $a_2$ | $a_4$ | $a_2$ |
| $a_3$ | $a_3$ | $a_1$ | $a_3$ | $a_1$ |
| $a_4$ | $a_2$ | $a_4$ | $a_2$ | $a_4$ |

This has H = {$a_1$, $a_3$} and K = {$a_2$, $a_4$} to be subgroupoids of G.

***Example 1.2.7:*** G is a groupoid given by the following table:

| * | $a_0$ | $a_1$ | $a_2$ | $a_3$ | $a_4$ | $a_5$ | $a_6$ | $a_7$ | $a_8$ | $a_9$ | $a_{10}$ | $a_{11}$ |
|---|---|---|---|---|---|---|---|---|---|---|---|---|
| $a_0$ | $a_0$ | $a_3$ | $a_6$ | $a_9$ | $a_0$ | $a_3$ | $a_6$ | $a_9$ | $a_0$ | $a_3$ | $a_6$ | $a_9$ |
| $a_1$ | $a_1$ | $a_4$ | $a_7$ | $a_{10}$ | $a_1$ | $a_4$ | $a_7$ | $a_{10}$ | $a_1$ | $a_4$ | $a_7$ | $a_{10}$ |
| $a_2$ | $a_2$ | $a_5$ | $a_8$ | $a_{11}$ | $a_2$ | $a_5$ | $a_8$ | $a_{11}$ | $a_2$ | $a_5$ | $a_8$ | $a_{11}$ |
| $a_3$ | $a_3$ | $a_6$ | $a_9$ | $a_0$ | $a_3$ | $a_6$ | $a_9$ | $a_0$ | $a_3$ | $a_6$ | $a_9$ | $a_0$ |
| $a_4$ | $a_4$ | $a_7$ | $a_{10}$ | $a_1$ | $a_4$ | $a_7$ | $a_{10}$ | $a_1$ | $a_4$ | $a_7$ | $a_{10}$ | $a_1$ |
| $a_5$ | $a_5$ | $a_8$ | $a_{11}$ | $a_2$ | $a_5$ | $a_8$ | $a_{11}$ | $a_2$ | $a_5$ | $a_8$ | $a_{11}$ | $a_2$ |
| $a_6$ | $a_6$ | $a_9$ | $a_0$ | $a_3$ | $a_6$ | $a_9$ | $a_0$ | $a_3$ | $a_6$ | $a_9$ | $a_0$ | $a_3$ |
| $a_7$ | $a_7$ | $a_{10}$ | $a_1$ | $a_4$ | $a_7$ | $a_{10}$ | $a_1$ | $a_4$ | $a_7$ | $a_{10}$ | $a_1$ | $a_4$ |
| $a_8$ | $a_8$ | $a_{11}$ | $a_2$ | $a_5$ | $a_8$ | $a_{11}$ | $a_2$ | $a_5$ | $a_8$ | $a_{11}$ | $a_2$ | $a_5$ |
| $a_9$ | $a_9$ | $a_0$ | $a_3$ | $a_6$ | $a_9$ | $a_0$ | $a_3$ | $a_6$ | $a_9$ | $a_0$ | $a_3$ | $a_6$ |
| $a_{10}$ | $a_{10}$ | $a_1$ | $a_4$ | $a_7$ | $a_{10}$ | $a_1$ | $a_4$ | $a_7$ | $a_{10}$ | $a_1$ | $a_4$ | $a_7$ |
| $a_{11}$ | $a_{11}$ | $a_2$ | $a_5$ | $a_8$ | $a_{11}$ | $a_2$ | $a_5$ | $a_8$ | $a_{11}$ | $a_2$ | $a_5$ | $a_8$ |

Clearly, $H_1$ = {$a_0$, $a_3$, $a_6$, $a_9$}, $H_2$ = {$a_2$, $a_5$, $a_8$, $a_{11}$} and $H_3$ = {$a_1$, $a_4$, $a_7$, $a_{10}$} are the three subgroupoids of G.

**PROBLEM 1:** Give an example of an infinite commutative groupoid with identity.

**PROBLEM 2:** How many groupoids of order 3 exists?

**PROBLEM 3:** How many groupoids of order 4 exists?

**PROBLEM 4:** How many commutative groupoids of order 5 exists?

**PROBLEM 5:** Give a new operation ∗ on R the set of reals so that R is a commutative groupoid with identity.

**PROBLEM 6:** Does a groupoid of order 5 have subgroupoids?

**PROBLEM 7:** Find all subgroupoids for the following groupoid:

| *   | $x_0$ | $x_1$ | $x_2$ | $x_3$ | $x_4$ | $x_5$ | $x_6$ | $x_7$ |
|-----|-------|-------|-------|-------|-------|-------|-------|-------|
| $x_0$ | $x_0$ | $x_6$ | $x_4$ | $x_2$ | $x_0$ | $x_6$ | $x_4$ | $x_2$ |
| $x_1$ | $x_2$ | $x_0$ | $x_6$ | $x_4$ | $x_2$ | $x_0$ | $x_6$ | $x_4$ |
| $x_2$ | $x_4$ | $x_2$ | $x_0$ | $x_6$ | $x_4$ | $x_2$ | $x_0$ | $x_6$ |
| $x_3$ | $x_6$ | $x_4$ | $x_2$ | $x_0$ | $x_6$ | $x_4$ | $x_2$ | $x_0$ |
| $x_4$ | $x_0$ | $x_6$ | $x_4$ | $x_2$ | $x_0$ | $x_6$ | $x_4$ | $x_2$ |
| $x_5$ | $x_2$ | $x_0$ | $x_6$ | $x_4$ | $x_2$ | $x_0$ | $x_6$ | $x_4$ |
| $x_6$ | $x_4$ | $x_2$ | $x_0$ | $x_6$ | $x_4$ | $x_2$ | $x_0$ | $x_6$ |
| $x_7$ | $x_6$ | $x_4$ | $x_2$ | $x_0$ | $x_6$ | $x_4$ | $x_2$ | $x_0$ |

Is G a commutative groupoid? Can G have subgroupoids of order 2?

## 1.3 Definition of Semigroup with Examples

In this section, we just recall the definition of a semigroup since groupoids are the generalization of semigroups and as they are the basic tools to define a Smarandache groupoid. Hence we define semigroups and give some examples of them.

**DEFINITION:** *Let S be a non empty set S is said to be a semigroup if on S is defined a binary operation ' • ' such that*

1. *For all a, b ∈ S we have a • b ∈ S (closure).*
2. *For all a, b, c ∈ S we have a • (b • c) = (a • b) • c (associative law).*

*(S, •) is a semigroup.*

**DEFINITION:** *If in a semigroup (S, •) we have a • b = b • a for all a, b ∈ S we say S is a commutative semigroup.*

**DEFINITION:** *Let S be a semigroup, if the number of elements in a semigroup is finite we say the semigroup S is of finite order otherwise S is of infinite order.*

**DEFINITION:** *Let (S, •) be a semigroup, H be a proper subset of S. We say H is a subsemigroup if (H, •) is itself a semigroup.*

**DEFINITION:** *Let (S, •) be a semigroup if S contains an element e such that e • s = s • e = s for all s ∈ S we say S is a semigroup with identity e or S is a monoid.*

***Example 1.3.1:*** Let $Z_8$ = {0, 1, 2, ... , 7} be the set of integers modulo 8. $Z_8$ is a semigroup under multiplication mod 8, which is commutative and has 1 to be the identity.

***Example 1.3.2:*** Let $S_{2 \times 2} = \{(a_{ij}) \mid a_{ij} \in Z\}$, the set of all $2 \times 2$ matrices with entries from Z. $S_{2 \times 2}$ is an infinite semigroup under the matrix multiplication which is non commutative and has $\begin{pmatrix} 1 & 0 \\ 0 & 1 \end{pmatrix}$ to be its identity, that is $S_{2 \times 2}$ is a monoid.

***Example 1.3.3:*** Let $Z^+ = \{1, 2, \ldots, \infty\}$. $Z^+$ is a semigroup under addition. Clearly, $Z^+$ is only a semigroup and not a monoid.

***Example 1.3.4:*** Let $R_{2 \times 2} = \{(a_{ij}) \mid a_{ij} \in Z_2 = \{0, 1\}\}$, be the set of all $2 \times 2$ matrices with entries from the prime field $Z_2 = \{0, 1\}$. $R_{2 \times 2}$ is a semigroup under matrix multiplication. Clearly, $R_{2 \times 2}$ is a non commutative monoid and is of finite order.

Thus we have seen semigroups or monoids of finite order and infinite order both commutative and non commutative types. Now we define an ideal in a semigroup.

**DEFINITION:** *Let $(S, \bullet)$ be a semigroup, a non empty subset I of S is said to be a right ideal of S if I is a subsemigroup of S and for $s \in S$ and $a \in I$ we have $as \in I$.*
*Similarly, one can define left ideal in a semigroup. We say I is an ideal in a semigroup if I is simultaneously left and right ideal of S.*

*Note:* If S is a commutative semigroup we see the concept of right ideal and left ideal coincide.

***Example 1.3.5:*** Let $Z_9 = \{0, 1, 2, \ldots, 8\}$ be a semigroup under multiplication modulo 9. Clearly $G = \{0, 3\}$ is a subsemigroup of S.

***Example 1.3.6:*** Let $S_{3 \times 3} = \{(a_{ij}) \mid a_{ij} \in Z_3 = \{0, 1, 2\}\}$ be the set of all $3 \times 3$ matrices; under the operation matrix multiplication, $S_{3 \times 3}$ is a semigroup. Let $A_{3 \times 3} = \{\text{set of all upper triangular matrices with entries from } Z_3 = \{0, 1, 2\}\}$. Clearly, $A_{3 \times 3}$ is a subsemigroup of $S_{3 \times 3}$.

***Example 1.3.7:*** Let $Z_{12} = \{0, 1, 2, \ldots, 11\}$ be the semigroup under multiplication modulo 12. Clearly, $A = \{0, 6\}$ is an ideal of $Z_{12}$, also $B = \{0, 2, 4, 8, 10\}$ is an ideal of $Z_{12}$.

***Example 1.3.8:*** Let $Z^+ \cup \{0\}$ be the semigroup under multiplication. Clearly, $p(Z^+ \cup \{0\}) = \{0, p, 2p, \ldots\}$ is an ideal of $Z^+ \cup \{0\}$ for every integer $p \in Z^+$.

**DEFINITION:** *Let $(S, \bullet)$ and $(S', *)$ be any two semigroups. We say a map $\phi: (S, \bullet) \to (S', *)$ is a semigroup homomorphism if $\phi(s_1 \bullet s_2) = \phi(s_1) * \phi(s_2)$ for all $s_1, s_2 \in S$.*
*If $\phi$ is one to one and onto we say $\phi$ is an isomorphism on the semigroups.*

**DEFINITION:** *Let S be a semigroup with identity e. We say $s \in S$ is invertible or has an inverse in S if there exist a $s' \in S$ such that $ss' = s's = e$.*

<u>Remark:</u> If in a semigroup S with identity every element is invertible S becomes a group.

***Example 1.3.9:*** The inverse of 5 in $Z_6 = \{0, 1, 2, \ldots, 5\}$ under multiplication is itself for $5.5 \equiv 1 \pmod 6$ so 5 is invertible where as 2, 3, 4 and 0 are non invertible in $Z_6$.

**DEFINITION:** *Let (S, •) be a semigroup. We say an element $x \in S$ is an idempotent if $x \cdot x = x$.*

***Example 1.3.10:*** Let $Z_{10} = \{0, 1, 2, \ldots, 9\}$ be a semigroup under multiplication mod 10. Clearly $5^2 \equiv 5 \pmod{10}$.

**PROBLEM 1:** Find all the ideals in the semigroup $Z_{36} = \{0, 1, \ldots, 35\}$ under multiplication.

**PROBLEM 2:** Find all the subsemigroups which are not ideals in $R_{2\times 2}$ given in example 1.3.4.

**PROBLEM 3:** Find all the ideals in the semigroup $S_{3\times 3}$. $S_{3\times 3}$ given in example 1.3.6.

**PROBLEM 4:** Find all the ideals and subsemigroups of $Z_{128} = \{0, 1, 2, \ldots, 127\}$, $Z_{128}$ is a semigroup under multiplication modulo 128.

**PROBLEM 5:** Find a right ideal in $S_{2\times 2}$ given in example 1.3.2, which is not a left ideal.

**PROBLEM 6:** Find a left ideal in $S_{2\times 2}$ given in example 1.3.2, which is not a right ideal.

**PROBLEM 7:** Find a subsemigroup, which is not an ideal in $S_{2\times 2}$ given in example 1.3.2.

**PROBLEM 8:** Find a homomorphism between the semigroups $S = \{Z_{30}, \times\}$, the semigroup under multiplication modulo 30 and $R_{2\times 2} = \{(a_{ij}) \mid a_{ij} \in Z_2 = \{0, 1\}\}$ is a semigroup under matrix multiplication.

**PROBLEM 9:** How many elements are there in $S_{3\times 3}$ given in example 1.3.6.?

**PROBLEM 10:** Can $R_{2\times 2}$ given in example 1.3.4 be isomorphic to $Z_{16} = \{0, 1, 2, \ldots, 15\}$ semigroup under multiplication modulo 16? Justify your answer.

**PROBLEM 11:** Can $(Z^+, \times)$ be isomorphic with $(Z^+, +)$? Justify your answer.

**PROBLEM 12:** Find a homomorphism between the semigroups $Z_7 = \{0, 1, 2, \ldots, 6\}$ under multiplication modulo 7 and $Z_{14} = \{0, 1, 2, \ldots, 13\}$ under multiplication modulo 14.

**PROBLEM 13:** Find an invertible element in the semigroup $R_{2\times 2}$ given in example 1.3.4.

**PROBLEM 14:** Find a non invertible element in the semigroup $S_{3\times 3}$ given in example 1.3.6.

**PROBLEM 15:** Find all the invertible elements in $Z_{32} = \{0, 1, 2, \ldots, 31\}$; semigroup under multiplication modulo 32.

**PROBLEM 16:** Find all non invertible elements in $Z_{42} = \{0, 1, 2, \ldots, 41\}$, semigroup under multiplication modulo 42.

**PROBLEM 17:** Give an idempotent element in $(Z_{24}, \times)$, semigroup under multiplication modulo 24.

**PROBLEM 18:** Give an example of an idempotent element in the semigroup $S_{3\times 3}$ given in example 1.3.6.

**PROBLEM 19:** Give an example of an idempotent element in the semigroup $R_{2\times 2}$ given in example 1.3.4.

**PROBLEM 20:** Does an idempotent exist in $(Z^+, \times)$? Justify.

**PROBLEM 21:** Does an idempotent exist in $(Z^+, +)$? Justify.

## 1.4. Smarandache Groupoids

In this section we just recall the definition of Smarandache Groupoid (SG) studied in the year 2002, and explain this definition by examples so as to make the concept easy to grasp as the main aim of this book is the study of Smarandache groupoid using $Z_n$ and introduce some new notion in them. In Chapter 4 a complete work of SG and their properties are given and it is very recent (2002).

**DEFINITION:** *A Smarandache groupoid G is a groupoid which has a proper subset $S \subset G$ which is a semigroup under the operations of G.*

*Example 1.4.1:* Let $(G, *)$ be a groupoid on the modulo integers 6. $G = \{0, 1, 2, 3, 4, 5\}$ given by the following table:

| * | 0 | 1 | 2 | 3 | 4 | 5 |
|---|---|---|---|---|---|---|
| 0 | 0 | 3 | 0 | 3 | 0 | 3 |
| 1 | 1 | 4 | 1 | 4 | 1 | 4 |
| 2 | 2 | 5 | 2 | 5 | 2 | 5 |
| 3 | 3 | 0 | 3 | 0 | 3 | 0 |
| 4 | 4 | 1 | 4 | 1 | 4 | 1 |
| 5 | 5 | 2 | 5 | 2 | 5 | 2 |

Clearly, the following tables give proper subsets of G which are semigroups under *.

| * | 0 | 3 |
|---|---|---|
| 0 | 0 | 3 |
| 3 | 3 | 0 |

| * | 1 | 4 |
|---|---|---|
| 1 | 4 | 1 |
| 4 | 1 | 4 |

| * | 2 | 5 |
|---|---|---|
| 2 | 2 | 5 |
| 5 | 5 | 2 |

So (G, *) is a SG.

**DEFINITION:** *Let (G, *) be a SG. The number of elements in G is called the order of the SG. If the number of elements is finite, we say the SG is of finite order or finite otherwise it is infinite order or infinite.*

***Example 1.4.2:*** Let $Z_6$ = {0, 1, 2, 3, 4, 5} be a groupoid given by the following table:

| × | 0 | 1 | 2 | 3 | 4 | 5 |
|---|---|---|---|---|---|---|
| 0 | 0 | 5 | 4 | 3 | 2 | 1 |
| 1 | 2 | 1 | 0 | 5 | 4 | 3 |
| 2 | 4 | 3 | 2 | 1 | 0 | 5 |
| 3 | 0 | 5 | 4 | 3 | 2 | 1 |
| 4 | 2 | 1 | 0 | 5 | 4 | 3 |
| 5 | 4 | 3 | 2 | 1 | 0 | 5 |

Clearly, this is a SG as every singleton is a semigroup.

***Example 1.4.3:*** Let G be a SG given by the following table:

| * | $a_1$ | $a_2$ | $a_3$ | $a_4$ |
|---|---|---|---|---|
| $a_1$ | $a_1$ | $a_4$ | $a_3$ | $a_2$ |
| $a_2$ | $a_3$ | $a_2$ | $a_1$ | $a_4$ |
| $a_3$ | $a_1$ | $a_4$ | $a_3$ | $a_2$ |
| $a_4$ | $a_3$ | $a_2$ | $a_1$ | $a_4$ |

A= {$a_4$} is a semigroup as $a_4 * a_4 = a_4$. Hence, G is a SG.

**PROBLEM 1:** Give an example of a SG of order 8 and find all subsets which are semigroups.

**PROBLEM 2:** Does their exist a SG of order 3?

**PROBLEM 3:** How many SGs of order 3 exists?

**PROBLEM 4:** How many SGs of order 4 exists?

**PROBLEM 5:** Give an example of a SG of infinite order.

**PROBLEM 6:** Find an example of a SG of order 7.

**PROBLEM 7:** Can a SG of order 8 have subsets, which are semigroups of order 5?

## 1.5 Loops and its Properties

In this section we just recall the definition of loops and illustrate them with examples. A lot of study has been carried out on loops and special types of loops have been defined. Here we give the definition of a loop and illustrate them with examples, as groupoids are the generalization of loops and all loops are obviously groupoids and not conversely.

**DEFINITION:** *(L, •) is said to be a loop where L is a non empty set and '•' a binary operation, called the product defined on L satisfying the following conditions:*

1. *For a, b ∈ L we have a • b ∈ L (closure property).*
2. *There exist an element e ∈ L such that a • e = e • a = a for all a ∈ L (e is called the identity element of L).*
3. *For every ordered pair (a, b) ∈ L × L their exists a unique pair (x, y) ∈ L × L such that ax = b and ya = b.*

We say L is a commutative loop if a • b = b • a for all a, b ∈ L.

The number of elements in a loop L is the order of the loop denoted by o (L) or |L|, if |L| < ∞, it is finite otherwise infinite.

In this section, all the examples of the loops given are constructed by us.

*Example 1.5.1:* Let L = {e, $a_1$, $a_2$, $a_3$, $a_4$, $a_5$}. Define a binary operation ∗ given by the following table:

| ∗ | e | $a_1$ | $a_2$ | $a_3$ | $a_4$ | $a_5$ |
|---|---|---|---|---|---|---|
| e | e | $a_1$ | $a_2$ | $a_3$ | $a_4$ | $a_5$ |
| $a_1$ | $a_1$ | e | $a_3$ | $a_5$ | $a_2$ | $a_4$ |
| $a_2$ | $a_2$ | $a_5$ | e | $a_4$ | $a_1$ | $a_3$ |
| $a_3$ | $a_3$ | $a_4$ | $a_1$ | e | $a_5$ | $a_2$ |
| $a_4$ | $a_4$ | $a_3$ | $a_5$ | $a_2$ | e | $a_1$ |
| $a_5$ | $a_5$ | $a_2$ | $a_4$ | $a_1$ | $a_3$ | e |

This loop is non commutative.

*Example 1.5.2:* Let L = {e, $a_1$, $a_2$, ... , $a_7$} be the loop given by the following table:

| ∗ | e | $a_1$ | $a_2$ | $a_3$ | $a_4$ | $a_5$ | $a_6$ | $a_7$ |
|---|---|---|---|---|---|---|---|---|
| e | e | $a_1$ | $a_2$ | $a_3$ | $a_4$ | $a_5$ | $a_6$ | $a_7$ |
| $a_1$ | $a_1$ | e | $a_5$ | $a_2$ | $a_6$ | $a_3$ | $a_7$ | $a_4$ |
| $a_2$ | $a_2$ | $a_5$ | e | $a_6$ | $a_3$ | $a_7$ | $a_4$ | $a_1$ |
| $a_3$ | $a_3$ | $a_2$ | $a_6$ | e | $a_7$ | $a_4$ | $a_1$ | $a_5$ |
| $a_4$ | $a_4$ | $a_6$ | $a_3$ | $a_7$ | e | $a_1$ | $a_5$ | $a_2$ |
| $a_5$ | $a_5$ | $a_3$ | $a_7$ | $a_4$ | $a_1$ | e | $a_2$ | $a_6$ |
| $a_6$ | $a_6$ | $a_7$ | $a_4$ | $a_1$ | $a_5$ | $a_2$ | e | $a_3$ |
| $a_7$ | $a_7$ | $a_4$ | $a_1$ | $a_5$ | $a_2$ | $a_6$ | $a_3$ | e |

The loop (L, ∗) is a commutative one with order 8.

**DEFINITION:** *Let L be a loop. A nonempty subset P of L is called a sub loop of L if P itself is a loop under the operation of L.*

*Example 1.5.3:* For the loop L given in example 1.5.1. $\{e, a_1\}$, $\{e, a_2\}$, $\{e, a_3\}$, $\{e, a_4\}$ and $\{e, a_5\}$ are the sub loops of order 2.

For definition of Bol loop, Bruck loop, Moufang loop, alternative loop etc. please refer Bruck .R.H. A Survey of binary systems (1958).

*Example 1.5.4:* Let $L = \{e, a_1, a_2, a_3, a_4, a_5\}$; $\{L, *\}$ is a loop given by the following table:

| ∗     | e     | $a_1$ | $a_2$ | $a_3$ | $a_4$ | $a_5$ |
|-------|-------|-------|-------|-------|-------|-------|
| e     | e     | $a_1$ | $a_2$ | $a_3$ | $a_4$ | $a_5$ |
| $a_1$ | $a_1$ | e     | $a_5$ | $a_4$ | $a_3$ | $a_2$ |
| $a_2$ | $a_2$ | $a_3$ | e     | $a_1$ | $a_5$ | $a_4$ |
| $a_3$ | $a_3$ | $a_5$ | $a_4$ | e     | $a_2$ | $a_1$ |
| $a_4$ | $a_4$ | $a_2$ | $a_1$ | $a_5$ | e     | $a_3$ |
| $a_5$ | $a_5$ | $a_4$ | $a_3$ | $a_2$ | $a_1$ | e     |

The above is a non commutative loop of order 6.

**PROBLEM 1:** How many loops of order 4 exist?

**PROBLEM 2:** How many loops of order 5 exist?

**PROBLEM 3:** Can a loop of order 3 exist? Justify your answer.

**PROBLEM 4:** Give an example of a loop of order 7 and find all its sub loops.

**PROBLEM 5:** Can a loop L of finite order p, p a prime be generated by a single element? Justify your answer.

**PROBLEM 6:** Let $\{L, *\}$ be a loop given by the following table:

| ∗     | e     | $a_1$ | $a_2$ | $a_3$ | $a_4$ | $a_5$ |
|-------|-------|-------|-------|-------|-------|-------|
| e     | e     | $a_1$ | $a_2$ | $a_3$ | $a_4$ | $a_5$ |
| $a_1$ | $a_1$ | e     | $a_3$ | $a_5$ | $a_2$ | $a_4$ |
| $a_2$ | $a_2$ | $a_5$ | e     | $a_4$ | $a_1$ | $a_3$ |
| $a_3$ | $a_3$ | $a_4$ | $a_1$ | e     | $a_5$ | $a_2$ |
| $a_4$ | $a_4$ | $a_3$ | $a_5$ | $a_2$ | e     | $a_1$ |
| $a_5$ | $a_5$ | $a_2$ | $a_4$ | $a_1$ | $a_3$ | e     |

Does this loop have sub loops? Is it commutative?

**PROBLEM 7:** Give an example of a loop of order 13. Does this loop have proper sub loops?

## Supplementary Reading

1. R.H. Bruck, *A Survey of binary systems,* Berlin Springer Verlag, (1958).

2. Ivan Nivan and H.S.Zukerman, *Introduction to number theory,* Wiley Eastern Limited, (1989).

3. W.B. Vasantha Kandasamy, *Smarandache Groupoids,* http://www.gallup.unm.edu/~smarandache/Groupoids.pdf

# CHAPTER TWO
# GROUPOIDS AND ITS PROPERTIES

This chapter is completely devoted to the introduction to groupoids and the study of several properties and new concepts in them. We make this explicit by illustrative examples. Apart from this, several identities like Moufang, Bol, etc are defined on a groupoid and conditions are derived from them to satisfy these identities. It is certainly a loss to see when semigroups, which are nothing but associative groupoids are studied, groupoids which are more generalized concept of semigroups are not dealt well in general books.

In the first chapter, we define special identities in groupoids followed by study of sub structures in a groupoid and direct product of groupoids. We define some special properties of groupoids like conjugate subgroupoids and normal subgroupoids and obtain some interesting results about them.

## 2.1 Special Properties in Groupoids

In this section we introduce the notion of special identities like Bol identity, Moufang identity etc., to groupoids and give examples of each.

**DEFINITION:** *A groupoid G is said to be a Moufang groupoid if it satisfies the Moufang identity $(xy)(zx) = (x(yz))x$ for all $x, y, z$ in G.*

*Example 2.1.1:* Let $(G, *)$ be a groupoid given by the following table:

| * | $a_0$ | $a_1$ | $a_2$ | $a_3$ |
|---|---|---|---|---|
| $a_0$ | $a_0$ | $a_2$ | $a_0$ | $a_2$ |
| $a_1$ | $a_1$ | $a_3$ | $a_1$ | $a_3$ |
| $a_2$ | $a_2$ | $a_0$ | $a_2$ | $a_0$ |
| $a_3$ | $a_3$ | $a_1$ | $a_3$ | $a_1$ |

$(a_1 * a_3) * (a_2 * a_1) = (a_1 * (a_3 * a_2)) * a_1$
$(a_1 * a_3) * (a_2 * a_3) = a_3 * a_0 = a_3$.

Now $(a_1 * (a_3 * a_2)) a_1 = (a_1 * a_3) * a_1 = a_3 * a_1 = a_1$. Since $(a_1 * a_3) (a_2 * a_1) \neq (a_1 (a_3 * a_2)) * a_1$, we see G is not a Moufang groupoid.

**DEFINITION:** *A groupoid G is said to be a Bol groupoid if G satisfies the Bol identity $((xy) z) y = x ((yz) y)$ for all x, y, z in G.*

*Example 2.1.2:* Let (G, *) be a groupoid given by the following table:

| * | $a_0$ | $a_1$ | $a_2$ | $a_3$ | $a_4$ | $a_5$ |
|---|---|---|---|---|---|---|
| $a_0$ | $a_0$ | $a_3$ | $a_0$ | $a_3$ | $a_0$ | $a_3$ |
| $a_1$ | $a_2$ | $a_5$ | $a_2$ | $a_5$ | $a_2$ | $a_5$ |
| $a_2$ | $a_4$ | $a_1$ | $a_4$ | $a_1$ | $a_4$ | $a_1$ |
| $a_3$ | $a_0$ | $a_3$ | $a_0$ | $a_3$ | $a_0$ | $a_3$ |
| $a_4$ | $a_2$ | $a_5$ | $a_2$ | $a_5$ | $a_2$ | $a_5$ |
| $a_5$ | $a_4$ | $a_1$ | $a_4$ | $a_1$ | $a_4$ | $a_1$ |

It can be easily verified that the groupoid (G, *) is a Bol groupoid.

**DEFINITION:** *A groupoid G is said to be a P-groupoid if $(xy) x = x (yx)$ for all $x, y \in G$.*

**DEFINITION:** *A groupoid G is said to be right alternative if it satisfies the identity $(xy) y = x (yy)$ for all $x, y \in G$. Similarly we define G to be left alternative if $(xx) y = x (xy)$ for all $x, y \in G$.*

**DEFINITION:** *A groupoid G is alternative if it is both right and left alternative, simultaneously.*

*Example 2.1.3:* Consider the groupoid (G, *) given by the following table:

| * | $a_0$ | $a_1$ | $a_2$ | $a_3$ | $a_4$ | $a_5$ | $a_6$ | $a_7$ | $a_8$ | $a_9$ | $a_{10}$ | $a_{11}$ |
|---|---|---|---|---|---|---|---|---|---|---|---|---|
| $a_0$ | $a_0$ | $a_9$ | $a_6$ | $a_3$ | $a_0$ | $a_9$ | $a_6$ | $a_3$ | $a_0$ | $a_9$ | $a_6$ | $a_3$ |
| $a_1$ | $a_4$ | $a_1$ | $a_{10}$ | $a_7$ | $a_4$ | $a_1$ | $a_{10}$ | $a_7$ | $a_4$ | $a_1$ | $a_{10}$ | $a_7$ |
| $a_2$ | $a_8$ | $a_5$ | $a_2$ | $a_{11}$ | $a_8$ | $a_5$ | $a_2$ | $a_{11}$ | $a_8$ | $a_5$ | $a_2$ | $a_{11}$ |
| $a_3$ | $a_0$ | $a_9$ | $a_6$ | $a_3$ | $a_0$ | $a_9$ | $a_6$ | $a_3$ | $a_0$ | $a_9$ | $a_6$ | $a_3$ |
| $a_4$ | $a_4$ | $a_1$ | $a_{10}$ | $a_7$ | $a_4$ | $a_1$ | $a_{10}$ | $a_7$ | $a_4$ | $a_1$ | $a_{10}$ | $a_7$ |
| $a_5$ | $a_8$ | $a_5$ | $a_2$ | $a_{11}$ | $a_8$ | $a_5$ | $a_2$ | $a_{11}$ | $a_8$ | $a_5$ | $a_2$ | $a_{11}$ |
| $a_6$ | $a_0$ | $a_9$ | $a_6$ | $a_3$ | $a_0$ | $a_9$ | $a_6$ | $a_3$ | $a_0$ | $a_9$ | $a_6$ | $a_3$ |
| $a_7$ | $a_4$ | $a_1$ | $a_{10}$ | $a_7$ | $a_4$ | $a_1$ | $a_{10}$ | $a_7$ | $a_4$ | $a_1$ | $a_{10}$ | $a_7$ |
| $a_8$ | $a_8$ | $a_5$ | $a_2$ | $a_{11}$ | $a_8$ | $a_5$ | $a_2$ | $a_{11}$ | $a_8$ | $a_5$ | $a_2$ | $a_{11}$ |
| $a_9$ | $a_0$ | $a_9$ | $a_6$ | $a_3$ | $a_0$ | $a_9$ | $a_6$ | $a_3$ | $a_0$ | $a_9$ | $a_6$ | $a_3$ |
| $a_{10}$ | $a_4$ | $a_1$ | $a_{10}$ | $a_7$ | $a_4$ | $a_1$ | $a_{10}$ | $a_7$ | $a_4$ | $a_1$ | $a_{10}$ | $a_7$ |
| $a_{11}$ | $a_8$ | $a_5$ | $a_2$ | $a_{11}$ | $a_8$ | $a_5$ | $a_2$ | $a_{11}$ | $a_8$ | $a_5$ | $a_2$ | $a_{11}$ |

It is left for the reader to verify this groupoid is a P-groupoid.

*Example 2.1.4:* Let (G, *) be a groupoid given by the following table:

| * | $a_0$ | $a_1$ | $a_2$ | $a_3$ | $a_4$ | $a_5$ | $a_6$ | $a_7$ | $a_8$ | $a_9$ |
|---|---|---|---|---|---|---|---|---|---|---|
| $a_0$ | $a_0$ | $a_6$ | $a_2$ | $a_8$ | $a_4$ | $a_0$ | $a_6$ | $a_2$ | $a_8$ | $a_4$ |
| $a_1$ | $a_5$ | $a_1$ | $a_7$ | $a_3$ | $a_9$ | $a_5$ | $a_1$ | $a_7$ | $a_3$ | $a_9$ |
| $a_2$ | $a_0$ | $a_6$ | $a_2$ | $a_8$ | $a_4$ | $a_0$ | $a_6$ | $a_2$ | $a_8$ | $a_4$ |
| $a_3$ | $a_5$ | $a_1$ | $a_7$ | $a_3$ | $a_9$ | $a_5$ | $a_1$ | $a_7$ | $a_3$ | $a_9$ |
| $a_4$ | $a_0$ | $a_6$ | $a_2$ | $a_8$ | $a_4$ | $a_0$ | $a_6$ | $a_2$ | $a_8$ | $a_4$ |
| $a_5$ | $a_5$ | $a_1$ | $a_7$ | $a_3$ | $a_9$ | $a_5$ | $a_1$ | $a_7$ | $a_3$ | $a_9$ |
| $a_6$ | $a_0$ | $a_6$ | $a_2$ | $a_8$ | $a_4$ | $a_0$ | $a_6$ | $a_2$ | $a_8$ | $a_4$ |
| $a_7$ | $a_5$ | $a_1$ | $a_7$ | $a_3$ | $a_9$ | $a_5$ | $a_1$ | $a_7$ | $a_3$ | $a_9$ |
| $a_8$ | $a_0$ | $a_6$ | $a_2$ | $a_8$ | $a_4$ | $a_0$ | $a_6$ | $a_2$ | $a_8$ | $a_4$ |
| $a_9$ | $a_5$ | $a_1$ | $a_7$ | $a_3$ | $a_9$ | $a_5$ | $a_1$ | $a_7$ | $a_3$ | $a_9$ |

It is left for the reader to verify that (G, ∗) is a right alternative groupoid.

**PROBLEM 1:** Give an example of a groupoid which is Bol, Moufang and right alternative.

**PROBLEM 2:** Can a Bol groupoid be a Moufang groupoid? Justify your answer.

**PROBLEM 3:** Give an example of an alternative groupoid with identity.

**PROBLEM 4:** Can a groupoid G of prime order satisfy the Bol identity for all elements in G?

**PROBLEM 5:** Can a P-groupoid be a Bol-groupoid? Justify your answer?

**PROBLEM 6:** Does there exist a prime order groupoid, which is non commutative?

## 2.2 Substructures in Groupoids

In this section, we define properties on subsets of a groupoid like subgroupoid, conjugate subgroupoid, ideals etc. and derive some interesting results using them.

**DEFINITION:** *Let (G, ∗) be a groupoid. A proper subset H of G is said to be a subgroupoid of G if (H, ∗) is itself a groupoid.*

*Example 2.2.1:* Let (G, ∗) be a groupoid given by the following table. G = {$a_0$, $a_1$, ... , $a_8$}.

| * | $a_0$ | $a_1$ | $a_2$ | $a_3$ | $a_4$ | $a_5$ | $a_6$ | $a_7$ | $a_8$ |
|---|---|---|---|---|---|---|---|---|---|
| $a_0$ | $a_0$ | $a_3$ | $a_6$ | $a_0$ | $a_3$ | $a_6$ | $a_0$ | $a_3$ | $a_6$ |
| $a_1$ | $a_5$ | $a_8$ | $a_2$ | $a_5$ | $a_8$ | $a_2$ | $a_5$ | $a_8$ | $a_2$ |
| $a_2$ | $a_1$ | $a_4$ | $a_7$ | $a_1$ | $a_4$ | $a_7$ | $a_1$ | $a_4$ | $a_7$ |
| $a_3$ | $a_6$ | $a_0$ | $a_3$ | $a_6$ | $a_0$ | $a_3$ | $a_6$ | $a_0$ | $a_3$ |
| $a_4$ | $a_2$ | $a_5$ | $a_8$ | $a_2$ | $a_5$ | $a_8$ | $a_2$ | $a_5$ | $a_8$ |
| $a_5$ | $a_7$ | $a_1$ | $a_4$ | $a_7$ | $a_1$ | $a_4$ | $a_7$ | $a_1$ | $a_4$ |
| $a_6$ | $a_3$ | $a_6$ | $a_0$ | $a_3$ | $a_6$ | $a_0$ | $a_3$ | $a_6$ | $a_0$ |
| $a_7$ | $a_8$ | $a_2$ | $a_5$ | $a_8$ | $a_2$ | $a_5$ | $a_8$ | $a_2$ | $a_5$ |
| $a_8$ | $a_4$ | $a_7$ | $a_1$ | $a_4$ | $a_7$ | $a_1$ | $a_4$ | $a_7$ | $a_1$ |

The subgroupids are given by the following tables:

| * | $a_0$ | $a_3$ | $a_6$ |
|---|---|---|---|
| $a_0$ | $a_0$ | $a_0$ | $a_0$ |
| $a_3$ | $a_6$ | $a_6$ | $a_6$ |
| $a_6$ | $a_3$ | $a_3$ | $a_3$ |

| * | $a_1$ | $a_2$ | $a_4$ | $a_5$ | $a_7$ | $a_8$ |
|---|---|---|---|---|---|---|
| $a_1$ | $a_8$ | $a_2$ | $a_8$ | $a_2$ | $a_8$ | $a_2$ |
| $a_2$ | $a_4$ | $a_7$ | $a_4$ | $a_7$ | $a_4$ | $a_7$ |
| $a_4$ | $a_5$ | $a_8$ | $a_5$ | $a_8$ | $a_5$ | $a_8$ |
| $a_5$ | $a_1$ | $a_4$ | $a_1$ | $a_4$ | $a_1$ | $a_4$ |
| $a_7$ | $a_2$ | $a_5$ | $a_2$ | $a_5$ | $a_2$ | $a_5$ |
| $a_8$ | $a_7$ | $a_1$ | $a_7$ | $a_1$ | $a_7$ | $a_1$ |

*Example 2.2.2:* Let G = {$a_0$, $a_1$, $a_2$, $a_3$, $a_4$} be the groupoid is given by the following table:

| * | $a_0$ | $a_1$ | $a_2$ | $a_3$ | $a_4$ |
|---|---|---|---|---|---|
| $a_0$ | $a_0$ | $a_4$ | $a_3$ | $a_2$ | $a_1$ |
| $a_1$ | $a_2$ | $a_1$ | $a_0$ | $a_4$ | $a_3$ |
| $a_2$ | $a_4$ | $a_3$ | $a_2$ | $a_1$ | $a_0$ |
| $a_3$ | $a_1$ | $a_0$ | $a_4$ | $a_3$ | $a_2$ |
| $a_4$ | $a_3$ | $a_2$ | $a_1$ | $a_0$ | $a_4$ |

Clearly every singleton is a subgroupoid in G, as $a_i * a_i = a_i$ for i = 0, 1, 2, 3, 4.

**DEFINITION:** *A groupoid G is said to be an idempotent groupoid if $x^2 = x$ for all $x \in G$.*

The groupoid given in example 2.2.2 is an idempotent groupoid of order 5.

**DEFINITION:** *Let G be a groupoid. P a non empty proper subset of G, P is said to be a left ideal of the groupoid G if 1) P is a subgroupoid of G and 2) For all $x \in G$ and $a \in P$, $xa \in P$.*

*One can similarly define right ideal of the groupoid G. P is called an ideal if P is simultaneously a left and a right ideal of the groupoid G.*

*Example 2.2.3:* Consider the groupoid G given by the following table:

| * | $a_0$ | $a_1$ | $a_2$ | $a_3$ |
|---|---|---|---|---|
| $a_0$ | $a_0$ | $a_2$ | $a_0$ | $a_2$ |
| $a_1$ | $a_3$ | $a_1$ | $a_3$ | $a_1$ |
| $a_2$ | $a_2$ | $a_0$ | $a_2$ | $a_0$ |
| $a_3$ | $a_1$ | $a_3$ | $a_1$ | $a_3$ |

P = {$a_0$, $a_2$} and Q = {$a_1$, $a_3$} are right ideals of G. Clearly P and Q are not left ideals of G. In fact G has no left ideals only 2 right ideals.

*Example 2.2.4:* Let G be a groupoid given by the following table:

| * | $a_0$ | $a_1$ | $a_2$ | $a_3$ | $a_4$ | $a_5$ |
|---|---|---|---|---|---|---|
| $a_0$ | $a_0$ | $a_4$ | $a_2$ | $a_0$ | $a_4$ | $a_2$ |
| $a_1$ | $a_2$ | $a_0$ | $a_4$ | $a_2$ | $a_0$ | $a_4$ |
| $a_2$ | $a_4$ | $a_2$ | $a_0$ | $a_4$ | $a_2$ | $a_0$ |
| $a_3$ | $a_0$ | $a_4$ | $a_2$ | $a_0$ | $a_4$ | $a_2$ |
| $a_4$ | $a_2$ | $a_0$ | $a_4$ | $a_2$ | $a_0$ | $a_4$ |
| $a_5$ | $a_4$ | $a_2$ | $a_0$ | $a_4$ | $a_2$ | $a_0$ |

Clearly P = {$a_0$, $a_2$, $a_4$} is both a right and a left ideal of G; in fact an ideal of G.

*Example 2.2.5:* Let G be the groupoid given by the following table:

| * | $a_0$ | $a_1$ | $a_2$ | $a_3$ |
|---|---|---|---|---|
| $a_0$ | $a_0$ | $a_3$ | $a_2$ | $a_1$ |
| $a_1$ | $a_2$ | $a_1$ | $a_0$ | $a_3$ |
| $a_2$ | $a_0$ | $a_3$ | $a_2$ | $a_1$ |
| $a_3$ | $a_2$ | $a_1$ | $a_0$ | $a_3$ |

Clearly, G has P = {$a_0$, $a_2$} and Q = {$a_1$, $a_3$} to be the only left ideals of G and G has no right ideals. Now we proceed to define normal subgroupoids of G.

**DEFINITION:** *Let G be a groupoid A subgroupoid V of G is said to be a normal subgroupoid of G if*

1. *aV = Va*
2. *(Vx)y = V(xy)*
3. *y(xV) = (yx)V*

*for all x, y, a ∈ V.*

*Example 2.2.6:* Consider the groupoid given in example 2.2.4. Clearly P = {$a_0$, $a_2$, $a_4$} is a normal subgroupoid of G.

**DEFINITION:** *A groupoid G is said to be simple if it has no non trivial normal subgroupoids.*

*Example 2.2.7:* The groupoid G given by the following table is simple.

| * | $a_0$ | $a_1$ | $a_2$ | $a_3$ | $a_4$ | $a_5$ | $a_6$ |
|---|---|---|---|---|---|---|---|
| $a_0$ | $a_0$ | $a_4$ | $a_1$ | $a_5$ | $a_2$ | $a_6$ | $a_3$ |
| $a_1$ | $a_3$ | $a_0$ | $a_4$ | $a_1$ | $a_5$ | $a_2$ | $a_6$ |
| $a_2$ | $a_6$ | $a_3$ | $a_0$ | $a_4$ | $a_1$ | $a_5$ | $a_2$ |
| $a_3$ | $a_2$ | $a_6$ | $a_3$ | $a_0$ | $a_4$ | $a_1$ | $a_5$ |
| $a_4$ | $a_5$ | $a_2$ | $a_6$ | $a_3$ | $a_0$ | $a_4$ | $a_1$ |
| $a_5$ | $a_1$ | $a_5$ | $a_2$ | $a_6$ | $a_3$ | $a_0$ | $a_4$ |
| $a_6$ | $a_4$ | $a_1$ | $a_5$ | $a_2$ | $a_6$ | $a_3$ | $a_0$ |

It is left for the reader to verify $(G, *) = \{a_0, a_1, a_2, \ldots, a_6, *\}$ has no normal subgroupoids. Hence, G is simple.

**DEFINITION:** *A groupoid G is normal if*

1. $xG = Gx$
2. $G(xy) = (Gx)y$
3. $y(xG) = (yx)G$ *for all $x, y \in G$.*

*Example 2.2.8:* Let $G = \{a_0, a_1, a_2, \ldots, a_6\}$ be the groupoid given by the following table:

| *   | $a_0$ | $a_1$ | $a_2$ | $a_3$ | $a_4$ | $a_5$ | $a_6$ |
|-----|-------|-------|-------|-------|-------|-------|-------|
| $a_0$ | $a_0$ | $a_4$ | $a_1$ | $a_5$ | $a_2$ | $a_6$ | $a_3$ |
| $a_1$ | $a_3$ | $a_0$ | $a_4$ | $a_1$ | $a_5$ | $a_2$ | $a_6$ |
| $a_2$ | $a_6$ | $a_3$ | $a_0$ | $a_4$ | $a_1$ | $a_5$ | $a_2$ |
| $a_3$ | $a_2$ | $a_6$ | $a_3$ | $a_0$ | $a_4$ | $a_1$ | $a_5$ |
| $a_4$ | $a_5$ | $a_2$ | $a_6$ | $a_3$ | $a_0$ | $a_4$ | $a_1$ |
| $a_5$ | $a_1$ | $a_5$ | $a_2$ | $a_6$ | $a_3$ | $a_0$ | $a_4$ |
| $a_6$ | $a_4$ | $a_1$ | $a_5$ | $a_2$ | $a_6$ | $a_3$ | $a_0$ |

This groupoid is a normal groupoid.

This has no proper subgroupoid, which is normal.

Not all groupoids in general are normal by an example.

*Example 2.2.9:* Let $G = \{a_0, a_1, \ldots, a_9\}$ be the groupoid given by the following table:

| *   | $a_0$ | $a_1$ | $a_2$ | $a_3$ | $a_4$ | $a_5$ | $a_6$ | $a_7$ | $a_8$ | $a_9$ |
|-----|-------|-------|-------|-------|-------|-------|-------|-------|-------|-------|
| $a_0$ | $a_0$ | $a_2$ | $a_4$ | $a_6$ | $a_8$ | $a_0$ | $a_2$ | $a_4$ | $a_6$ | $a_8$ |
| $a_1$ | $a_1$ | $a_3$ | $a_5$ | $a_7$ | $a_9$ | $a_1$ | $a_3$ | $a_5$ | $a_7$ | $a_9$ |
| $a_2$ | $a_2$ | $a_4$ | $a_6$ | $a_8$ | $a_0$ | $a_2$ | $a_4$ | $a_6$ | $a_8$ | $a_0$ |
| $a_3$ | $a_3$ | $a_5$ | $a_7$ | $a_9$ | $a_1$ | $a_3$ | $a_5$ | $a_7$ | $a_9$ | $a_1$ |
| $a_4$ | $a_4$ | $a_6$ | $a_8$ | $a_0$ | $a_2$ | $a_4$ | $a_6$ | $a_8$ | $a_0$ | $a_2$ |
| $a_5$ | $a_5$ | $a_7$ | $a_9$ | $a_1$ | $a_3$ | $a_5$ | $a_7$ | $a_9$ | $a_1$ | $a_3$ |
| $a_6$ | $a_6$ | $a_8$ | $a_0$ | $a_2$ | $a_4$ | $a_6$ | $a_8$ | $a_0$ | $a_2$ | $a_4$ |
| $a_7$ | $a_7$ | $a_9$ | $a_1$ | $a_3$ | $a_5$ | $a_7$ | $a_9$ | $a_1$ | $a_3$ | $a_9$ |
| $a_8$ | $a_8$ | $a_0$ | $a_2$ | $a_4$ | $a_6$ | $a_8$ | $a_0$ | $a_2$ | $a_4$ | $a_6$ |
| $a_9$ | $a_9$ | $a_1$ | $a_3$ | $a_5$ | $a_7$ | $a_9$ | $a_1$ | $a_3$ | $a_5$ | $a_7$ |

This groupoid is not normal. Clearly $a_1 \bullet (a_0\ a_1\ a_2\ a_3\ a_4\ a_5\ a_6\ a_7\ a_8\ a_9) = a_1G = \{a_1, a_3, a_5, a_7, a_9\}$ where as $Ga_1 = G$. So G is not a normal groupoid.

**DEFINITION:** *Let G be a groupoid H and K be two proper subgroupoids of G, with $H \cap K = \phi$. We say H is conjugate with K if there exists a $x \in H$ such that $H = x K$ or $Kx$ ('or' in the mutually exclusive sense).*

***Example 2.2.10:*** Let G be a groupoid; G = {$a_0, a_1, \ldots, a_{10}, a_{11}$} defined by the following table:

| * | $a_0$ | $a_1$ | $a_2$ | $a_3$ | $a_4$ | $a_5$ | $a_6$ | $a_7$ | $a_8$ | $a_9$ | $a_{10}$ | $a_{11}$ |
|---|---|---|---|---|---|---|---|---|---|---|---|---|
| $a_0$ | $a_0$ | $a_3$ | $a_6$ | $a_9$ | $a_0$ | $a_3$ | $a_6$ | $a_9$ | $a_0$ | $a_3$ | $a_6$ | $a_9$ |
| $a_1$ | $a_1$ | $a_4$ | $a_7$ | $a_{10}$ | $a_1$ | $a_4$ | $a_7$ | $a_{10}$ | $a_1$ | $a_4$ | $a_7$ | $a_{10}$ |
| $a_2$ | $a_2$ | $a_5$ | $a_8$ | $a_{11}$ | $a_2$ | $a_5$ | $a_8$ | $a_{11}$ | $a_2$ | $a_5$ | $a_8$ | $a_{11}$ |
| $a_3$ | $a_3$ | $a_6$ | $a_9$ | $a_0$ | $a_3$ | $a_6$ | $a_9$ | $a_0$ | $a_3$ | $a_6$ | $a_9$ | $a_0$ |
| $a_4$ | $a_4$ | $a_7$ | $a_{10}$ | $a_1$ | $a_4$ | $a_7$ | $a_{10}$ | $a_1$ | $a_4$ | $a_7$ | $a_{10}$ | $a_1$ |
| $a_5$ | $a_5$ | $a_8$ | $a_{11}$ | $a_2$ | $a_5$ | $a_8$ | $a_{11}$ | $a_2$ | $a_5$ | $a_8$ | $a_{11}$ | $a_2$ |
| $a_6$ | $a_6$ | $a_9$ | $a_0$ | $a_3$ | $a_6$ | $a_9$ | $a_0$ | $a_3$ | $a_6$ | $a_9$ | $a_0$ | $a_3$ |
| $a_7$ | $a_7$ | $a_{10}$ | $a_1$ | $a_4$ | $a_7$ | $a_{10}$ | $a_1$ | $a_4$ | $a_7$ | $a_{10}$ | $a_1$ | $a_4$ |
| $a_8$ | $a_8$ | $a_{11}$ | $a_2$ | $a_5$ | $a_8$ | $a_{11}$ | $a_2$ | $a_5$ | $a_8$ | $a_{11}$ | $a_2$ | $a_5$ |
| $a_9$ | $a_9$ | $a_0$ | $a_3$ | $a_6$ | $a_9$ | $a_0$ | $a_3$ | $a_6$ | $a_9$ | $a_0$ | $a_3$ | $a_3$ |
| $a_{10}$ | $a_{10}$ | $a_1$ | $a_4$ | $a_7$ | $a_{10}$ | $a_1$ | $a_4$ | $a_7$ | $a_{10}$ | $a_1$ | $a_4$ | $a_7$ |
| $a_{11}$ | $a_{11}$ | $a_2$ | $a_5$ | $a_8$ | $a_{11}$ | $a_2$ | $a_5$ | $a_8$ | $a_{11}$ | $a_2$ | $a_5$ | $a_8$ |

Take K = {$a_0, a_3, a_6, a_9$} and H = {$a_2, a_5, a_8, a_{11}$}, H $\cap$ K = $\phi$.

It can be verified the subgroupoids H and K are conjugate with each other.

To get more normal subgroupoids and ideals, that is to make groupoids have richer structures we define direct products of groupoids.

**DEFINITION:** *Let $(G_1, \theta_1), (G_2, \theta_2), \ldots, (G_n, \theta_n)$ be n groupoids with $\theta_i$ binary operations defined on each $G_i$, i = 1, 2, 3, ..., n. The direct product of $G_1, \ldots, G_n$ denoted by $G = G_1 \times \ldots \times G_n = \{(g_1, \ldots, g_n) \mid g_i \in G_i\}$ by component wise multiplication on G, G becomes a groupoid.*

*For if $g = (g_1, \ldots, g_n)$ and $h = (h_1, \ldots, h_n)$ then $g \bullet h = \{(g_1\theta_1 h_1, g_2\theta_2 h_2, \ldots, g_n\theta_n h_n)\}$. Clearly, $gh \in G$. Hence G is a groupoid.*

This direct product helps us to construct more groupoids with desired properties. For example by using several simple groupoids, we can get a groupoid, which may not be simple. Likewise, we can get normal groupoids and normal subgroupoids using the direct product.

***Example 2.2.11:*** Let $(G_1, *)$ and $(G_2, o)$ be 2 groupoids given by the following tables:

| * | $a_0$ | $a_1$ | $a_2$ |
|---|---|---|---|
| $a_0$ | $a_0$ | $a_2$ | $a_1$ |
| $a_1$ | $a_1$ | $a_0$ | $a_2$ |
| $a_2$ | $a_2$ | $a_1$ | $a_0$ |

| o | $x_0$ | $x_1$ | $x_2$ | $x_3$ |
|---|---|---|---|---|
| $x_0$ | $x_0$ | $x_1$ | $x_2$ | $x_3$ |
| $x_1$ | $x_3$ | $x_0$ | $x_1$ | $x_2$ |
| $x_2$ | $x_2$ | $x_3$ | $x_0$ | $x_1$ |
| $x_3$ | $x_1$ | $x_2$ | $x_3$ | $x_0$ |

$G_1 \times G_2 = \{(a_i, x_j) \mid a_i \in G_1, i = 0, 1, 2\ x_j \in G_2, j = 0, 1, 2, 3\}$. $G_1 \times G_2$ is a groupoid of order 12. It is left for the reader to find the subgroupoids, if any, normal subgroupoids if any, ideals if any (right, left) and so on.

**PROBLEM 1:** Can a groupoid of order 11 have normal subgroupoids?

**PROBLEM 2:** Will all groupoids of prime order be normal groupoids?

**PROBLEM 3:** Give an example of a groupoid of order 16 where all of its subgroupoids are normal subgroupoids.

**PROBLEM 4:** Give an example of a groupoid, which has no proper subgroupoids.

**PROBLEM 5:** Does their exist a simple groupoid of order 24?

**PROBLEM 6:** Give an example of a groupoid in which all subgroupoids are ideals.

**PROBLEM 7:** Give an example of a groupoid, which has no ideals.

**PROBLEM 8:** Can we get a normal subgroupoid using two normal groupoids?

**PROBLEM 9:** Using two groupoids which have no right ideals can we find right ideals in the direct product?

**PROBLEM 10:** Do the existence of an ideal in a groupoid G guarantee the existence of an ideal in the direct product $G \times G_1$ for an arbitrary groupoid $G_1$?

**PROBLEM 11:** Is the direct product of 2 normal groupoids, a normal groupoid? Justify!

**PROBLEM 12:** Is the direct product of simple groupoids $G_1$, $G_2$, $G_3$ make G, ($G = G_1 \times G_2 \times G_3$) also simple? Justify.

**PROBLEM 13:** Does the direct product of a simple groupoid $G_1$ with a groupoid $G_2$, which has normal subgroupoids, be simple? That is will $G_1 \times G_2$ be simple?

**PROBLEM 14:** Can problem 11, 12 and 13 be true even if the groupoids do not posses identity?

**PROBLEM 15:** Give an example of a groupoid in which every pair of subgroupoids is conjugate with each other.

**PROBLEM 16:** Does their exist a groupoid which has no pair of subgroupoids which are conjugate?

**PROBLEM 17:** Can a simple groupoid have conjugate subgroupoids?

## 2.3 Some Special Properties of a Groupoid

In this section, we define some special properties of groupoids like inner commutative, left (right) identity, left (right) zero divisor, center of the groupoid, conjugate pair etc. We illustrate these properties by several interesting examples. This is done mainly as, we do not have any textbook on groupoids.

**DEFINITION:** *Let G be a groupoid, we say G is inner commutative if every subgroupoid of G is commutative.*

***Example 2.3.1:*** Consider the groupoid $G = \{x_0, x_1, ... , x_{11}\}$ given by the following table:

| * | $x_0$ | $x_1$ | $x_2$ | $x_3$ | $x_4$ | $x_5$ | $x_6$ | $x_7$ | $x_8$ | $x_9$ | $x_{10}$ | $x_{11}$ |
|---|---|---|---|---|---|---|---|---|---|---|---|---|
| $x_0$ | $x_0$ | $x_4$ | $x_8$ | $x_0$ | $x_4$ | $x_8$ | $x_0$ | $x_4$ | $x_8$ | $x_0$ | $x_4$ | $x_8$ |
| $x_1$ | $x_1$ | $x_5$ | $x_9$ | $x_1$ | $x_5$ | $x_9$ | $x_1$ | $x_5$ | $x_9$ | $x_1$ | $x_5$ | $x_9$ |
| $x_2$ | $x_2$ | $x_6$ | $x_{10}$ | $x_2$ | $x_6$ | $x_{10}$ | $x_2$ | $x_6$ | $x_{10}$ | $x_2$ | $x_6$ | $x_{10}$ |
| $x_3$ | $x_3$ | $x_7$ | $x_{11}$ | $x_3$ | $x_7$ | $x_{11}$ | $x_3$ | $x_7$ | $x_{11}$ | $x_3$ | $x_7$ | $x_{11}$ |
| $x_4$ | $x_4$ | $x_8$ | $x_0$ | $x_4$ | $x_8$ | $x_0$ | $x_4$ | $x_8$ | $x_0$ | $x_4$ | $x_8$ | $x_0$ |
| $x_5$ | $x_5$ | $x_9$ | $x_1$ | $x_5$ | $x_9$ | $x_1$ | $x_5$ | $x_9$ | $x_1$ | $x_5$ | $x_9$ | $x_1$ |
| $x_6$ | $x_6$ | $x_{10}$ | $x_2$ | $x_6$ | $x_{10}$ | $x_2$ | $x_6$ | $x_{10}$ | $x_2$ | $x_6$ | $x_{10}$ | $x_2$ |
| $x_7$ | $x_7$ | $x_{11}$ | $x_3$ | $x_7$ | $x_{11}$ | $x_3$ | $x_7$ | $x_{11}$ | $x_3$ | $x_7$ | $x_{11}$ | $x_3$ |
| $x_8$ | $x_8$ | $x_0$ | $x_4$ | $x_8$ | $x_0$ | $x_4$ | $x_8$ | $x_0$ | $x_4$ | $x_8$ | $x_0$ | $x_4$ |
| $x_9$ | $x_9$ | $x_1$ | $x_5$ | $x_9$ | $x_1$ | $x_5$ | $x_9$ | $x_1$ | $x_5$ | $x_9$ | $x_1$ | $x_5$ |
| $x_{10}$ | $x_{10}$ | $x_2$ | $x_6$ | $x_{10}$ | $x_2$ | $x_6$ | $x_{10}$ | $x_2$ | $x_6$ | $x_{10}$ | $x_2$ | $x_6$ |
| $x_{11}$ | $x_{11}$ | $x_3$ | $x_7$ | $x_{11}$ | $x_3$ | $x_7$ | $x_{11}$ | $x_3$ | $x_7$ | $x_{11}$ | $x_3$ | $x_7$ |

The subgroupoids of G are $P_1 = \{x_0, x_4, x_8\}$, $P_2 = \{x_1, x_5, x_9\}$, $P_3 = \{x_2, x_6, x_{10}\}$ and $P_4 = \{x_3, x_7, x_{11}\}$. Clearly this groupoid is inner commutative.

**DEFINITION:** *Let G be a groupoid we say an element $e \in G$ is a left identity if $ea = a$ for all $a \in G$. Similarly we can define right identity of the groupoid G, if $e \in G$ happens to be simultaneously both right and left identity we say the groupoid G has an identity.*

***Example 2.3.2:*** In groupoid given in example 2.3.1 we see $a * e = a$ for $e = x_0$ and for all $a \in G$ thus $x_0$ is the right identity of G.

**DEFINITION:** *Let G be a groupoid. We say a in G has right zero divisor if $a * b = 0$ for some $b \neq 0$ in G and a in G has left zero divisor if $b * a = 0$. We say G has zero divisors if $a \bullet b = 0$ and $b * a = 0$ for $a, b \in G \setminus \{0\}$ A zero divisor in G can be left or right divisor.*

**DEFINITION:** *Let G be a groupoid. The center of the groupoid $C(G) = \{x \in G \mid ax = xa$ for all $a \in G\}$.*

**DEFINITION:** *Let G be a groupoid. We say $a, b \in G$ is a conjugate pair if $a = bx$ (or xa for some $x \in G$) and $b = ay$ (or ya for some $y \in G$).*

*Example 2.3.3:* Let G be a groupoid given by the following table:

| • | 0 | 1 | 2 | 3 | 4 | 5 | 6 |
|---|---|---|---|---|---|---|---|
| 0 | 0 | 3 | 6 | 2 | 5 | 1 | 4 |
| 1 | 2 | 5 | 1 | 4 | 0 | 3 | 6 |
| 2 | 4 | 0 | 3 | 6 | 2 | 5 | 1 |
| 3 | 6 | 2 | 5 | 1 | 4 | 0 | 3 |
| 4 | 1 | 4 | 0 | 3 | 6 | 2 | 5 |
| 5 | 3 | 6 | 2 | 5 | 1 | 4 | 0 |
| 6 | 5 | 1 | 4 | 0 | 3 | 6 | 2 |

| | |
|---|---|
| $1 * 4 = 0$ but $4 * 1 = 4$ | $4 * 2 = 0$ but $2 * 4 = 2$ |
| $2 * 1 = 0$ but $1 * 2 = 1$ | $5 * 6 = 0$ but $6 * 5 = 6$ |
| $3 * 5 = 0$ but $5 * 3 = 5$ | $6 * 3 = 0$ but $3 * 6 = 3$ |

Thus we see from this example $a \cdot b = 0$ need not imply $b \cdot a = 0$. That is why we can have an element in a groupoid to be right zero divisor or a left zero divisor.

*Example 2.3.4:* Let G be a groupoid given by the following table:

| * | $a_0$ | $a_1$ | $a_2$ |
|---|---|---|---|
| $a_0$ | $a_0$ | $a_2$ | $a_1$ |
| $a_1$ | $a_1$ | $a_0$ | $a_2$ |
| $a_2$ | $a_2$ | $a_1$ | $a_0$ |

The center of G is $C(G) = \phi$.

*Example 2.3.5:* Let G be a groupoid given by the following table:

| * | 0 | 1 | 2 |
|---|---|---|---|
| 0 | 0 | 1 | 2 |
| 1 | 2 | 0 | 1 |
| 2 | 1 | 2 | 0 |

Clearly (1, 2) is a conjugate pair as $2.0 = 1$ and $1.0 = 2$, (0,1) is also a conjugate pair as $0 = 1.1$, $1 = 0.1$; 0 is conjugate to 1.

**DEFINITION:** *A groupoid G is said to be strictly non commutative if $a \cdot b \neq b * a$ for any a, b $\in$ G where $a \neq b$.*

**DEFINITION:** *Let G be a groupoid of order n. An element a in G is said to be right conjugate with b in G if we can find x, y $\in$ G such that $a \cdot x = b$ and $b \cdot y = a$ ($x * a = b$ and $y * b = a$). Similarly, we define left conjugate.*

<u>**PROBLEM 1:**</u> Give an example of a groupoid in which every element has only left zero divisor.

**PROBLEM 2:** Construct a groupoid of order 25 in which every element has a right conjugate.

**PROBLEM 3:** Does their exist a groupoid of order 19 which is inner commutative but not commutative?

**PROBLEM 4:** Can a groupoid of order 17 have a non empty center?

**PROBLEM 5:** Give an example of a groupoid of order 24 in which every element is either a left or a right zero divisor.

**PROBLEM 6:** Construct a groupoid of order 18, which has only right identity.

**PROBLEM 7:** Can a non prime order groupoid have left identity? Justify your answer by an example!

## 2.4 Infinite Groupoids and its Properties

In this section, we introduce several interesting classes of groupoids of infinite order constructed using integers or rationals or reals. We obtain some important properties satisfied by them. It is pertinent to note that we are not in a position to get infinite Smarandache groupoids, which remains as an open problem.

**DEFINITION:** *Let $Z^+$ be the set of integers. Define an operation $*$ on $Z^+$ by $x * y = mx + ny$ where $m, n \in Z^+$, $m < \infty$ and $n < \infty$ $(m, n) = 1$ and $m \neq n$. Clearly $\{Z^+, *, (m, n)\}$ is a groupoid denoted by $Z^+ (m, n)$. We have for varying m and n get infinite number of groupoids of infinite order denoted by $\mathbb{Z}^+$.*

***Example 2.4.1:*** $Z^+ (1, 3) = \{4, 5, 6, 7, \ldots \infty\}$ is a non commutative and non associative groupoid of infinite order with no identity.

**THEOREM 2.4.1:** *No groupoid in the class $\mathbb{Z}^+$ is commutative.*

*Proof:* It is left for the reader to verify this result.

**THEOREM 2.4.2:** *No groupoid in the class $\mathbb{Z}^+$ is associative.*

*Proof:* Simple number theoretic methods will yield the result.

**COROLLARY 2.4.3:** *No groupoid in this class $\mathbb{Z}^+$ is a semigroup.*

*Proof:* Follows from theorem 2.4.2.

**THEOREM 2.4.4:** *No groupoid in this class $\mathbb{Z}^+$ is a P-groupoid.*

*Proof:* Let $Z^+(m, n)$ be a groupoid in $Z^+$. If $Z^+(m, n)$ is said to be a P-groupoid it must satisfy $(x * y) * x = x * (y * x)$ that is $m^2x + mny + nx = mx + mny + n^2x$. This is impossible for any $m \neq n$, $(m, n) = 1$ and $m, n \in Z^+$. Hence the claim.

**THEOREM 2.4.5:** *Groupoids $Z^+(1, 2m)$, $m \in Z^+$ in $Z^+$ always have subgroupoids of infinite order.*

*Proof:* Follows by simple number theoretic techniques.

**THEOREM 2.4.6:** *The groupoids $Z^+(1, p)$ where p is a prime can be partitioned as conjugate subgroupoids.*

*Proof:* Simple number theoretic methods will yield the result. We can now construct still many more classes of groupoids using $Z^+$ by making $m \neq n$ but m need not be relatively prime with n; still another new class by making $m = n$. Similarly by replacing $Z^+$ by $Q^+$ or $R^+$ also, we can construct different classes of infinite groupoids and similar results can be got using $Q^+$ and $R^+$.

**PROBLEM 1:** Find all the subgroupoids in $Z^+(3, 6)$.

**PROBLEM 2:** Does $Z^+(3, 11)$ have conjugate subgroupoids?

**PROBLEM 3:** Find subgroupoids of $Q^+(11/3, 7/3)$.

**PROBLEM 4:** Can $R^+(\sqrt{2}, \sqrt{3})$ have subgroupoids, which are conjugate with each other?

**PROBLEM 5:** Is $Q(7, 11)$ a P-groupoid?

**PROBLEM 6:** Can $R(2, 4)$ be a Moufang groupoid?

**PROBLEM 7:** Can $Q(7, 49)$ have conjugate subgroupoids?

**PROBLEM 8:** Find $(t, u)$ so that $Q(t, u)$ has conjugate subgroupoids which partition $Q(t, u)$.

**PROBLEM 9:** If $Z(m, n)$ has conjugate subgroupoids that partition $Z(m, n)$ find m and n.

## Supplementary Reading

1. R.H.Bruck, *A Survey of Binary Systems,* Springer Verlag, (1958).

2. W.B.Vasantha Kandasamy, *On Ordered Groupoids and its Groupoid Rings,* Jour. of Maths and Comp. Sci., Vol 9, 145-147, (1996).

**CHAPTER THREE**

# NEW CLASSES OF GROUPOIDS USING $Z_n$

In this chapter, we introduce four new classes of groupoids using $Z_n$. We study their properties. These new classes are different in their own way and they enjoy nice properties. To the best of the authors knowledge this chapter is the first attempt to construct non abstract classes of finite groupoids using $Z_n$. Throughout this chapter $Z_n$ denotes the set of integers modulo n that is $Z_n = \{0, 1, 2, 3, ... , n\}$; $n \geq 3$, $n < \infty$.

## 3.1 Definition of the Class of Groupoids Z (n)

Here we define a new class of groupoids denoted by Z (n) using $Z_n$ and study their various properties.

**DEFINITION:** *Let $Z_n = \{0, 1, 2, ... , n – 1\}$ $n \geq 3$. For a, b $\in Z_n \setminus \{0\}$ define a binary operation $*$ on $Z_n$ as follows. $a * b = ta + ub \pmod{n}$ where t,u are 2 distinct elements in $Z_n \setminus \{0\}$ and (t, u) =1 here ' + ' is the usual addition of 2 integers and ' ta ' means the product of the two integers t and a. We denote this groupoid by $\{Z_n, (t, u), *\}$ or in short by $Z_n$ (t, u).*

It is interesting to note that for varying t, u $\in Z_n \setminus \{0\}$ with (t, u) = 1 we get a collection of groupoids for a fixed integer n. This collection of groupoids is denoted by Z (n) that is Z (n) = $\{Z_n, (t, u), * \mid$ for integers t, u $\in Z_n \setminus \{0\}$ such that (t, u) = 1$\}$. Clearly every groupoid in this class is of order n.

*Example 3.1.1:* Using $Z_3 = \{0, 1, 2\}$. The groupoid $\{Z_3, (1, 2), *\} = (Z_3 (1, 2))$ is given by the following table:

| * | 0 | 1 | 2 |
|---|---|---|---|
| 0 | 0 | 2 | 1 |
| 1 | 1 | 0 | 2 |
| 2 | 2 | 1 | 0 |

Clearly this groupoid is non associative and non commutative and its order is 3.

***Example 3.1.2:*** Let $Z_3 = \{0, 1, 2\}$, $Z_3(2, 1) \in Z(3)$ is a groupoid given by the following table:

| * | 0 | 1 | 2 |
|---|---|---|---|
| 0 | 0 | 1 | 2 |
| 1 | 2 | 0 | 1 |
| 2 | 1 | 2 | 0 |

This groupoid is also non commutative and non associative. Hence, in $Z(3)$ we have only 2 groupoids of order 3 using $Z_3$.

***Example 3.1.3:*** Let $Z_7 = \{0, 1, 2, \ldots, 6\}$. The groupoid $Z_7(3, 4)$ is given by the following table:

| * | 0 | 1 | 2 | 3 | 4 | 5 | 6 |
|---|---|---|---|---|---|---|---|
| 0 | 0 | 4 | 1 | 5 | 2 | 6 | 3 |
| 1 | 3 | 0 | 4 | 1 | 5 | 2 | 6 |
| 2 | 6 | 3 | 0 | 4 | 1 | 5 | 2 |
| 3 | 2 | 6 | 3 | 0 | 4 | 1 | 5 |
| 4 | 5 | 2 | 6 | 3 | 0 | 4 | 1 |
| 5 | 1 | 5 | 2 | 6 | 3 | 0 | 4 |
| 6 | 4 | 1 | 5 | 2 | 6 | 3 | 0 |

This is a finite groupoid of order 7, which is non commutative.

**THEOREM 3.1.1:** *Let $Z_n = \{0, 1, 2, \ldots, n\}$. A groupoid in $Z(n)$ is a semigroup if and only if $t^2 \equiv t \pmod{n}$ and $u^2 \equiv u \pmod{n}$ for $t, u \in Z_n \setminus \{0\}$ and $(t, u) = 1$.*

*Proof:* Given $Z_n = \{0, 1, 2, \ldots, n-1\}$. Let $t, u \in Z_n \setminus \{0\}$ with $t^2 \equiv t \pmod{n}$, $u^2 \equiv u \pmod{n}$ and $(t, u) = 1$ to show $Z_n(t, u)$ is a semigroup it is sufficient if we show $(a * b) * c = a * (b * c)$ for all $a, b, c \in Z_n$.

$a * (b * c) = a * (tb + uc) \pmod{n} = ta + utb + u^2c \pmod{n}$. Now $(a * b) * c = (ta + ub) * c \pmod{n} = t^2a + utb + uc \pmod{n}$. Using $t^2 \equiv t \pmod{n}$ and $u^2 \equiv u \pmod{n}$ we see $a * (b * c) = (a * b) * c$. Hence $Z_n(t, u)$ is a semigroup.

On the other hand if $Z_n(t, u)$ is a semigroup to show $t^2 \equiv t \pmod{n}$ and $u^2 \equiv u \pmod{n}$. Given $Z_n(t, u)$ is a semigroup so $a * (b * c) \equiv (a * b) * c \pmod{n}$ for all $a, b, c \in Z_n$.

Now $a * (b * c) = (a * b) * c$ implies $ta + tub + u^2c = t^2a + tub + uc$. So $(t - t^2)a + (u^2 - u)c \equiv o \pmod{n}$ this is possible for all $a, c \in Z_n$ only if $t^2 \equiv t \pmod{n}$ and $u^2 \equiv u \pmod{n}$. Hence the claim.

**COROLLARY 3.1.2:** *If in the above theorem when; n is a prime $Z(n)$ never contains a semigroup.*

*Proof:* Since $t^2 \equiv t \pmod{n}$ or $u^2 \equiv u \pmod{n}$ can never occur for $t, u \in Z_n \setminus \{0, 1\}$ when n is a prime, Z (n) has no semigroups.

**THEOREM 3.1.3:** *The groupoid $Z_n$ (t, u) is an idempotent groupoid if and only if $t + u \equiv 1$ (mod n).*

*Proof:* If $Z_n$ (t, u) is to be an idempotent groupoid then $a^2 \equiv a \pmod{n}$ for all $a \in Z_n$. Now $a^2 = a * a = ta + ua \equiv a \pmod{n}$ that is $(t + u - 1) a \equiv 0 \pmod{n}$. This is possible if and only if $t + u \equiv 1 \pmod{n}$.

**THEOREM 3.1.4:** *No groupoid in Z (n) has {0} as an ideal.*

*Proof:* Obvious by the very definition of groupoids in Z (n).

*Example 3.1.4:* Consider the groupoid $Z_4$ (2, 3) and $Z_4$ (3, 2) in Z (4) given by the following tables:

Table for $Z_4$ (2, 3)

| * | 0 | 1 | 2 | 3 |
|---|---|---|---|---|
| 0 | 0 | 3 | 2 | 1 |
| 1 | 2 | 1 | 0 | 3 |
| 2 | 0 | 3 | 2 | 1 |
| 3 | 2 | 1 | 0 | 3 |

P ={0, 2} and Q = {1, 3} are the left ideals of $Z_4$ (2, 3). These two left ideals are not right ideals.

Table for $Z_4$ (3, 2)

| * | 0 | 1 | 2 | 3 |
|---|---|---|---|---|
| 0 | 0 | 2 | 0 | 2 |
| 1 | 3 | 1 | 3 | 1 |
| 2 | 2 | 0 | 2 | 0 |
| 3 | 1 | 3 | 1 | 3 |

Clearly P = {0, 2} and Q = {1, 3} are right ideal of $Z_4$ (3, 2) and are not left ideals of $Z_4$ (3, 2).

Consequent of the example we have the following theorem.

**THEOREM 3.1.5:** *P is a left ideal of $Z_n$ (t, u) if and only if P is a right ideal of $Z_n$ (u, t).*

*Proof:* Simple number theoretic computations and the simple translation of rows to columns gives a transformation of $Z_n$ (t, u) to $Z_n$ (u, t) and this will yield the result.

*Example 3.1.5:* Let $Z_{10}$ (3, 7) be the groupoid given by the following table:

| * | 0 | 1 | 2 | 3 | 4 | 5 | 6 | 7 | 8 | 9 |
|---|---|---|---|---|---|---|---|---|---|---|
| 0 | 0 | 7 | 4 | 1 | 8 | 5 | 2 | 9 | 6 | 3 |
| 1 | 3 | 0 | 7 | 4 | 1 | 8 | 5 | 2 | 9 | 6 |
| 2 | 6 | 3 | 0 | 7 | 4 | 1 | 8 | 5 | 2 | 9 |
| 3 | 9 | 6 | 3 | 0 | 7 | 4 | 1 | 8 | 5 | 2 |
| 4 | 2 | 9 | 6 | 3 | 0 | 7 | 4 | 1 | 8 | 5 |
| 5 | 5 | 2 | 9 | 6 | 3 | 0 | 7 | 4 | 1 | 8 |
| 6 | 8 | 5 | 2 | 9 | 6 | 3 | 0 | 7 | 4 | 1 |
| 7 | 1 | 8 | 5 | 2 | 9 | 6 | 3 | 0 | 7 | 4 |
| 8 | 4 | 1 | 8 | 5 | 2 | 9 | 6 | 3 | 0 | 7 |
| 9 | 7 | 4 | 1 | 8 | 5 | 2 | 9 | 6 | 3 | 0 |

This groupoid has no left or right ideals in it.

**THEOREM 3.1.6:** *Let $Z_n$ (t, u) be a groupoid. If n = t + u where both t and u are primes then $Z_n$ (t, u) is simple.*

*Proof:* It is left for the reader to verify this using number theoretic techniques as every row and every column in the table have distinct n elements.

**COROLLARY 3.1.7:** *If $Z_p$ (t, u) is a groupoid such that t + u = p, (t, u) = 1 then also $Z_p$ (t, u) is simple.*

*Proof:* Left for the reader to verify.

**PROBLEM 1:** Find the number of groupoids in Z (17).

**PROBLEM 2:** How many groupoids of order 7 exist?

**PROBLEM 3:** How many groupoids in Z (14) are semigroups?

**PROBLEM 4:** Is $Z_9$ (3, 8) an idempotent groupoid?

**PROBLEM 5:** Does $Z_{11}$ (4, 7) satisfy Bol identity?

**PROBLEM 6:** Does there exists a Moufang groupoid in Z (15)?

**PROBLEM 7:** Find P-groupoids in Z (12).

**PROBLEM 8:** Find all the subgroupoids of $Z_3$ (7, 14).

**PROBLEM 9:** Find all right ideals of $Z_{21}$ (3,7).

**PROBLEM 10:** Does $Z_{21}$ (3, 7) have ideals?

**PROBLEM 11:** Can $Z_{11}$ (2, 8) have ideals?

**PROBLEM 12:** Find all the left ideals of $Z_{16}$ (3, 9).

**PROBLEM 13:** Does Z (19) contain any simple groupoid?

**PROBLEM 14:** Find normal groupoids of $Z_{16}$ (3, 9).

**PROBLEM 15:** Does Z (16) have simple groupoids?

**PROBLEM 16:** Characterize all groupoids in Z (15), which has no ideals right or left.

**PROBLEM 17:** Does there exist a groupoid in the class Z (18) in which (xx) y = x (xy) for all x, y ∈ $Z_{18}$?

**PROBLEM 18:** Can $Z_{11}$ (3, 4) be written as a partition of conjugate groupoids? Justify your answer.

**PROBLEM 19:** Find the conjugate groupoids of $Z_{15}$ (9, 8).

## 3.2 New Class of Groupoids $Z^*(n)$

In this section, we introduce another new class of groupoids $Z^*(n)$, which is a generalization of the class of groupoids Z (n) and the class of groupoids Z (n) is a proper subclass of $Z^*(n)$. We study these groupoids and obtain some interesting results about them.

**DEFINITION:** *Let $Z_n$ = {0, 1, 2, ... , n –1} n ≥ 3, n < ∞. Define ∗ a closed binary operation on $Z_n$ as follows. For any a, b ∈ $Z_n$ define a ∗ b = at + bu (mod n) where (t, u) need not always be relatively prime but t ≠ u and t, u ∈ $Z_n$ \ {0}.*

Clearly {$Z_n$, (t, u), ∗} is a groupoid of order n. For varying values of t and u we get a class of groupoids for the fixed integer n. This new class of groupoids is denoted by $Z^*(n)$. $Z^*(n)$ is an extended class of Z (n) that is Z (n) ⊂ $Z^*(n)$. We denote the groupoids in the class $Z^*(n)$ by $Z_n^*$(t, u). Even if ' ∗ ' is omitted by looking at (t, u) ≠ 1 we can easily guess that $Z_n$ (t, u) is in the class $Z^*(n)$. In this section we only concentrate on groupoids in $Z^*(n)$ and not in Z (n).

*Example 3.2.1:* Let $Z_5$ (2, 4) ∈ $Z^*(5)$ be the groupoid given by the following table:

| ∗ | 0 | 1 | 2 | 3 | 4 |
|---|---|---|---|---|---|
| 0 | 0 | 4 | 3 | 2 | 1 |
| 1 | 2 | 1 | 0 | 4 | 3 |
| 2 | 4 | 3 | 2 | 1 | 0 |
| 3 | 1 | 0 | 4 | 3 | 2 |
| 4 | 3 | 2 | 1 | 0 | 4 |

The subgroupoids of $Z_5$ (2, 4) are $A_1 = \{1\}$, $A_2 = \{2\}$, $A_3 = \{3\}$ and $A_4 = \{4\}$. Thus singletons are the only subgroupoids of the groupoid $Z_5$ (2, 4).

***Example 3.2.2:*** The groupoid $Z_6$ (2, 4) $\in Z^*(6)$ is given by the following table:

| * | 0 | 1 | 2 | 3 | 4 | 5 |
|---|---|---|---|---|---|---|
| 0 | 0 | 4 | 2 | 0 | 4 | 2 |
| 1 | 2 | 0 | 4 | 2 | 0 | 4 |
| 2 | 4 | 2 | 0 | 4 | 2 | 0 |
| 3 | 0 | 4 | 2 | 0 | 4 | 2 |
| 4 | 2 | 0 | 4 | 2 | 0 | 4 |
| 5 | 4 | 2 | 0 | 4 | 2 | 0 |

Clearly A= {0, 2, 4} is a subgroupoid of $Z_6$ (2, 4).

**THEOREM 3.2.1:** *The number of groupoids in $Z^*(n)$ is $(n - 1)(n - 2)$.*

*Proof:* Left for the reader to verify using simple number theoretic techniques.

**THEOREM 3.2.2:** *The number of groupoids in the class Z (n) is bounded by $(n - 1)(n - 2)$.*

*Proof:* Follows from the fact $Z(n) \subset Z^*(n)$ and $|Z(n)| < |Z^*(n)| = (n-1)(n-2)$. Hence the claim.

***Example 3.2.3:*** Consider the groupoid $Z_8$ (2, 6) in $Z^*(8)$ given by the following table:

| * | 0 | 1 | 2 | 3 | 4 | 5 | 6 | 7 |
|---|---|---|---|---|---|---|---|---|
| 0 | 0 | 6 | 4 | 2 | 0 | 6 | 4 | 2 |
| 1 | 2 | 0 | 6 | 4 | 2 | 0 | 6 | 4 |
| 2 | 4 | 2 | 0 | 6 | 4 | 2 | 0 | 6 |
| 3 | 6 | 4 | 2 | 0 | 6 | 4 | 2 | 0 |
| 4 | 0 | 6 | 4 | 2 | 0 | 6 | 4 | 2 |
| 5 | 2 | 0 | 6 | 4 | 2 | 0 | 6 | 4 |
| 6 | 4 | 2 | 0 | 6 | 4 | 2 | 0 | 6 |
| 7 | 6 | 4 | 2 | 0 | 6 | 4 | 2 | 0 |

Clearly {0, 2, 4, 6} is a subgroupoid of $Z_8$ (2, 6).

**THEOREM 3.2.3:** *Let $Z_n$ (t, u) be a groupoid in $Z^*(n)$ such that (t, u) = t, n = 2m, t / 2m and t + u = 2m. Then $Z_n$ (t, u) has subgroupoids of order 2m / t or n / t.*

*Proof:* Given n is even and t + u = n so that u = n – t. Thus $Z_n$ (t, u) = $Z_n$ (t, n – t). Now using the fact $t \cdot Z_n = \left\{0, t, 2t, 3t, \ldots, \left(\dfrac{n}{t} - 1\right)t\right\}$ that is t . $Z_n$ has only n / t elements and these n / t elements from a subgroupoid. Hence $Z_n$ (t, n – t) where (t, n – t) = t has only subgroupoids of order n / t.

***Example 3.2.4:*** Let the groupoid $Z_{12}$ (2, 10) in $Z^*(12)$ be given by the following table:

| * | 0 | 1 | 2 | 3 | 4 | 5 | 6 | 7 | 8 | 9 | 10 | 11 |
|---|---|---|---|---|---|---|---|---|---|---|----|----|
| 0 | 0 | 10 | 8 | 6 | 4 | 2 | 0 | 10 | 8 | 6 | 4 | 2 |
| 1 | 2 | 0 | 10 | 8 | 6 | 4 | 2 | 0 | 10 | 8 | 6 | 4 |
| 2 | 4 | 2 | 0 | 10 | 8 | 6 | 4 | 2 | 0 | 10 | 8 | 6 |
| 3 | 6 | 4 | 2 | 0 | 10 | 8 | 6 | 4 | 2 | 0 | 10 | 8 |
| 4 | 8 | 6 | 4 | 2 | 0 | 10 | 8 | 6 | 4 | 2 | 0 | 10 |
| 5 | 10 | 8 | 6 | 4 | 2 | 0 | 10 | 8 | 6 | 4 | 2 | 0 |
| 6 | 0 | 10 | 8 | 6 | 4 | 2 | 0 | 10 | 8 | 6 | 4 | 2 |
| 7 | 2 | 0 | 10 | 8 | 6 | 4 | 2 | 0 | 10 | 8 | 6 | 4 |
| 8 | 4 | 2 | 0 | 10 | 8 | 6 | 4 | 2 | 0 | 10 | 8 | 6 |
| 9 | 6 | 4 | 2 | 0 | 10 | 8 | 6 | 4 | 2 | 0 | 10 | 8 |
| 10 | 8 | 6 | 4 | 2 | 0 | 10 | 8 | 6 | 4 | 2 | 0 | 10 |
| 11 | 10 | 8 | 6 | 4 | 2 | 0 | 10 | 8 | 6 | 4 | 2 | 0 |

This has a subgroupoid P = {0, 2, 4, 6, 8, 10} of order 6 that is 12 / 2.

***Example 3.2.5:*** Consider the groupoid $Z_{12}$ (3, 9) the subgroupoid of $Z_{12}$ (3, 9) is {0, 3, 6, 9} which is order 4 that is 12 / 3.

***Example 3.2.6:*** Consider the groupoid $Z_{12}$ (4, 8). The subgroupoid of $Z_{12}$ (4, 8) is P = {0, 4, 8} of order 3 that is 12 / 4. Thus the class $Z^*(12)$ when t + u = 12, (t, u) = t gives 3 distinct groupoids having subgroupoids of order 6, 4 and 3.

**THEOREM 3.2.4:** *Let $Z_n$ (t, u) $\in Z^*(n)$ be a groupoid such that n is even t + u = n with (t, u) = t. Then $Z_n$ (t, u) has only one subgroupoid of order n / t and it is a normal subgroupoid of $Z_n$ (t, u).*

*Proof:* This is left for the reader to verify.

***Example 3.2.7:*** Let $Z_{12}$ (10, 8) be the groupoid given by the following table:

|   | 0 | 1 | 2 | 3 | 4 | 5 | 6 | 7 | 8 | 9 | 10 | 11 |
|---|---|---|---|---|---|---|---|---|---|---|----|----|
| 0 | 0 | 8 | 4 | 0 | 8 | 4 | 0 | 8 | 4 | 0 | 8 | 4 |
| 1 | 10 | 6 | 2 | 10 | 6 | 2 | 10 | 6 | 2 | 10 | 6 | 2 |
| 2 | 8 | 4 | 0 | 8 | 4 | 0 | 8 | 4 | 0 | 8 | 4 | 0 |
| 3 | 6 | 2 | 10 | 6 | 2 | 10 | 6 | 2 | 10 | 6 | 2 | 10 |
| 4 | 4 | 0 | 8 | 4 | 0 | 8 | 4 | 0 | 8 | 4 | 0 | 8 |
| 5 | 2 | 10 | 6 | 2 | 10 | 6 | 2 | 10 | 6 | 2 | 10 | 6 |
| 6 | 0 | 8 | 4 | 0 | 8 | 4 | 0 | 8 | 4 | 0 | 8 | 4 |
| 7 | 10 | 6 | 2 | 10 | 6 | 2 | 10 | 6 | 2 | 10 | 6 | 2 |
| 8 | 8 | 4 | 0 | 8 | 4 | 0 | 8 | 4 | 0 | 8 | 4 | 0 |
| 9 | 6 | 2 | 10 | 6 | 2 | 10 | 6 | 2 | 10 | 6 | 2 | 10 |
| 10 | 4 | 0 | 8 | 4 | 0 | 8 | 4 | 0 | 8 | 4 | 0 | 8 |
| 11 | 2 | 10 | 6 | 2 | 10 | 6 | 2 | 10 | 6 | 2 | 10 | 6 |

The subgroupoids of $Z_{12}$ (10, 8) are A = {0, 4, 8} and B = {2, 6, 10}. This has 2 subgroupoids of order 3.

*Example 3.2.8:* Let $Z_{10}$ (8, 4) $\in Z^*(10)$ be the groupoid of order 10.

The only subgroupoid of $Z_{10}$ (8, 4) is {0, 4, 2, 8, 6} which is of order 5 and (8, 4) = 4. So we see in our theorem 3.2.4 that only when (t, u) = t with t + u = n we need have subgroupoids but other wise also we have subgroupoids which is evident from the following table given for the groupoid $Z_{10}$ (8, 4).

| * | 0 | 1 | 2 | 3 | 4 | 5 | 6 | 7 | 8 | 9 |
|---|---|---|---|---|---|---|---|---|---|---|
| 0 | 0 | 4 | 8 | 2 | 6 | 0 | 4 | 8 | 2 | 6 |
| 1 | 8 | 2 | 6 | 0 | 4 | 8 | 2 | 6 | 0 | 4 |
| 2 | 6 | 0 | 4 | 8 | 2 | 6 | 0 | 4 | 8 | 2 |
| 3 | 4 | 8 | 2 | 6 | 0 | 4 | 8 | 2 | 6 | 0 |
| 4 | 2 | 6 | 0 | 4 | 8 | 2 | 6 | 0 | 4 | 8 |
| 5 | 0 | 4 | 8 | 2 | 6 | 0 | 4 | 8 | 2 | 6 |
| 6 | 8 | 2 | 6 | 0 | 4 | 8 | 2 | 6 | 0 | 4 |
| 7 | 6 | 0 | 4 | 8 | 2 | 6 | 0 | 4 | 8 | 2 |
| 8 | 4 | 8 | 2 | 6 | 0 | 4 | 8 | 2 | 6 | 0 |
| 9 | 2 | 6 | 0 | 4 | 8 | 2 | 6 | 0 | 4 | 8 |

**PROBLEM 1:** Show Z (3) = $Z^*$(3).

**PROBLEM 2:** Show Z (4) = $Z^*$(4).

**PROBLEM 3:** Prove Z (8) $\subsetneq Z^*$(8).

**PROBLEM 4:** Find all the groupoids in $Z^*$(9).

**PROBLEM 5:** Find the groupoids in $Z^*$(12) \ Z (12) (that is groupoids in $Z^*$(12) which are not in Z (12)).

**PROBLEM 6:** Find the subgroupoid of $Z_{36}$ (3, 33). What is its order?

**PROBLEM 7:** Does the groupoid $Z_{36}$ (7, 29) have subgroupoids?

**PROBLEM 8:** Find all the subgroupoids of these groupoids $Z_8$ (2, 6), $Z_8$ (4, 6) and $Z_8$ (4, 6).

**PROBLEM 9:** Find all the subgroupoids of the groupoids $Z_{12}$ (6, 8), $Z_{12}$ (10, 6), $Z_{12}$ (10, 8) and $Z_{12}$ (8, 4).

**PROBLEM 10:** Find normal subgroupoids in $Z_{24}$(2, 22).

**PROBLEM 11:** Does the groupoid $Z_{11}$ (3, 7) have normal subgroupoids?

## 3.3 On New Class of Groupoids $Z^{**}(n)$

In this section we introduce yet another new class of groupoids, which is a generalization of $Z^*(n)$ using $Z_n$ and study their properties and obtain some interesting results about them.

**DEFINITION:** *Let $Z_n = \{0, 1, 2, ... , n – 1\}$ $n \geq 3$, $n < \infty$. Define $*$ on $Z_n$ as $a * b = ta + ub$ (mod n) where t and u $\in Z_n \setminus \{0\}$ and t can also equal u. For a fixed n and for varying t and u we get a class of groupoids of order n which we denote by $Z^{**}(n)$.*

For a fixed integer n we have $(n – 1)^2$ groupoids in the class $Z^{**}(n)$.

*Example 3.3.1:* Let $Z_5 (2, 2) \in Z^{**}(5)$ be the groupoid given by the following table:

| * | 0 | 1 | 2 | 3 | 4 |
|---|---|---|---|---|---|
| 0 | 0 | 2 | 4 | 1 | 3 |
| 1 | 2 | 4 | 1 | 3 | 0 |
| 2 | 4 | 1 | 3 | 0 | 2 |
| 3 | 1 | 3 | 0 | 2 | 4 |
| 4 | 3 | 0 | 2 | 4 | 1 |

This groupoid is commutative which is non associative. This has no subgroupoids.

*Example 3.3.2:* Let $Z_7 (3, 3) \in Z^{**}(7)$ be the groupoid given by the following table:

| * | 0 | 1 | 2 | 3 | 4 | 5 | 6 |
|---|---|---|---|---|---|---|---|
| 0 | 0 | 3 | 6 | 2 | 5 | 1 | 4 |
| 1 | 3 | 6 | 2 | 5 | 1 | 4 | 0 |
| 2 | 6 | 2 | 5 | 1 | 4 | 0 | 3 |
| 3 | 2 | 5 | 1 | 4 | 0 | 3 | 6 |
| 4 | 5 | 1 | 4 | 0 | 3 | 6 | 2 |
| 5 | 1 | 4 | 0 | 3 | 6 | 2 | 5 |
| 6 | 4 | 0 | 3 | 6 | 2 | 5 | 1 |

This groupoid is also commutative and non associative having no subgroupoids.

*Example 3.3.3:* Let $Z_6 (2, 2) \in Z^{**}(6)$ be the groupoid given by the following table:

| * | 0 | 1 | 2 | 3 | 4 | 5 |
|---|---|---|---|---|---|---|
| 0 | 0 | 2 | 4 | 0 | 2 | 4 |
| 1 | 2 | 4 | 0 | 2 | 4 | 0 |
| 2 | 4 | 0 | 2 | 4 | 0 | 2 |
| 3 | 0 | 2 | 4 | 0 | 2 | 4 |
| 4 | 2 | 4 | 0 | 2 | 4 | 0 |
| 5 | 4 | 0 | 2 | 4 | 0 | 2 |

This is a commutative groupoid and has {0, 2, 4} as subgroupoid so we consider that if n is not a prime we may have subgroupoids.

***Example 3.3.4:*** $Z_6$ (5, 5) ∈ $Z^{**}$(6) is given by the following table:

| * | 0 | 1 | 2 | 3 | 4 | 5 |
|---|---|---|---|---|---|---|
| 0 | 0 | 5 | 4 | 3 | 2 | 1 |
| 1 | 5 | 4 | 3 | 2 | 1 | 0 |
| 2 | 4 | 3 | 2 | 1 | 0 | 5 |
| 3 | 3 | 2 | 1 | 0 | 5 | 4 |
| 4 | 2 | 1 | 0 | 5 | 4 | 3 |
| 5 | 1 | 0 | 5 | 4 | 3 | 2 |

This has {4} and {2} as subgroupoids.

**THEOREM 3.3.1:** *The groupoids $Z_n$ (t, t) ∈ $Z^{**}$(n) for t < n are commutative groupoids.*

*Proof:* Let a, b ∈ $Z_n$ now consider a * b = ta + tb (mod n) and b * a = tb + ta (mod n). Clearly a * b ≡ b * a (mod n) for every a, b ∈ $Z_n$. So $Z_n$ (t, t) is a commutative groupoids.

**THEOREM 3.3.2:** *The groupoids $Z_p$ (t, t) ∈ $Z^{**}$(5), p a prime, t < p are normal groupoids.*

*Proof:* Clearly in the groupoid $Z_n$ (t, t) we have $aZ_n$ (t, t) = $Z_n$ (t, t) a and ([$Z_n$ (t, t)] x) y = [$Z_n$ (t, t)] (xy) for all x, y ∈ $Z_n$ (t, t), y (x$Z_n$ (t, t)) = (yx) [$Z_n$ (t, t)] for all a, x, y ∈ $Z_n$ (t, t). The verification is simple and left for the reader.

**THEOREM 3.3.3:** *The groupoids $Z_n$ (t, t) ∈ $Z^{**}$(n) are P-groupoids.*

*Proof:* Let $Z_n$ (t, t) be a groupoid. If $Z_n$ (t, t) is to be a P-groupoid we must have (x * y) * x = x * (y * x). Now (x * y) *x = (tx + ty) * x = $t^2$x + $t^2$y + tx (mod n) and x * (y * x) = x * (ty + tx) = tx + $t^2$y + $t^2$x.

Clearly; (x * y) * x = x * (y * x) for all x, y ∈ $Z_n$ (t, t). Hence the claim.

**THEOREM 3.3.4:** *Let $Z_n$ (t, t) ∈ $Z^{**}$(n) be a groupoid, n a prime 1 < t < n. Clearly $Z_n$ (t, t) is not an alternative groupoid if p is a prime.*

*Proof:* Now $Z_n$ (t, t) where n is a prime and 1 < t < n, to show the groupoid is alternative we have to prove (xy) y = x (yy) for all x, y ∈ $Z_n$ (t, t). Since $t^2$ ≡ t (mod p) is never possible if t ≠ 1. Consider (x * y) * y = (tx + ty) * y = $t^2$x + $t^2$y + ty. Now x * (y * y) = x * (ty + ty) = $t^2$x + $t^2$y + $t^2$y. So, $t^2$x + $t^2$y + ty ≠ tx + $t^2$y + $t^2$y if n is a prime. Hence, $Z_n$ (t, t) is not an alternative groupoid.

**COROLLARY 3.3.5:** *Let $Z_n$ (t, t) be a groupoid where n is not a prime. $Z_n$(t, t) is an alternative groupoid if and only if $t^2$ ≡ t (mod n).*

*Proof:* Obvious from the fact that if $Z_n$ (t, t) is alternative then $t^2x + t^2y + ty \equiv tx + t^2y + t^2y$ (mod n) that is if and only if $t^2 \equiv t$ (mod n).

**PROBLEM 1:** Does $Z_{18}$ (9, 9) have any normal subgroupoid?

**PROBLEM 2:** Can $Z_{18}$ (9, 9) have two subgroupoids, which are conjugate subgroupoids?

**PROBLEM 3:** Prove $Z_7$ (4, 4) is an idempotent groupoid?

**PROBLEM 4:** Can $Z_9$ (t, t) be an idempotent groupoid if $t \neq 5$?

**PROBLEM 5:** Find all the subgroupoids of $Z_{32}$ (7, 7).

**PROBLEM 6:** Find all subgroupoids in $Z_{19}$ (3, 3).

**PROBLEM 7:** Does $Z_{32}$ (8, 8) have normal subgroupoid?

**PROBLEM 8:** Can $Z_{32}$ (16, 16) have conjugate normal subgroupoid?

**PROBLEM 9:** Can $Z_{21}$ (7, 7) be a P-groupoid?

**PROBLEM 10:** Is $Z_{25}$ (17, 17) a normal groupoid?

## 3.4 On Groupoids $Z^{***}(n)$

In this chapter, we introduce the most generalized class of groupoids using $Z_n$. We have $Z(n) \subset Z^*(n) \subset Z^{**}(n) \subset Z^{***}(n)$. Thus, $Z^{***}(n)$ is the groupoid, which has many more interesting properties than the other classes. We now proceed on to define this new class of groupoids.

**DEFINITION:** *Let $Z_n = \{0, 1, 2, \ldots, n-1\}$ $n \geq 3$, $n < \infty$. Define $*$ on $Z_n$ as follows. $a * b = ta + ub$ (mod n) where $t, u \in Z_n$. Here t or u can also be zero.*

For a fixed n if we collect the set of all groupoids and denote it by $Z^{***}(n)$. Clearly $Z^{***}(n)$ contains $n(n-1)$ groupoids each of order n.

*Example 3.4.1:* The groupoid $Z_6$ (3, 0) $\in Z^{***}(6)$ be given by the following table:

| * | 0 | 1 | 2 | 3 | 4 | 5 |
|---|---|---|---|---|---|---|
| 0 | 0 | 0 | 0 | 0 | 0 | 0 |
| 1 | 3 | 3 | 3 | 3 | 3 | 3 |
| 2 | 0 | 0 | 0 | 0 | 0 | 0 |
| 3 | 3 | 3 | 3 | 3 | 3 | 3 |
| 4 | 0 | 0 | 0 | 0 | 0 | 0 |
| 5 | 3 | 3 | 3 | 3 | 3 | 3 |

Clearly this is a non commutative groupoid having {0, 3} as its subgroupoid.

***Example 3.4.2:*** Consider the groupoid $Z_6$ (0, 2) in $Z^{***}(6)$ given by the following table:

| * | 0 | 1 | 2 | 3 | 4 | 5 |
|---|---|---|---|---|---|---|
| 0 | 0 | 2 | 4 | 0 | 2 | 4 |
| 1 | 0 | 2 | 4 | 0 | 2 | 4 |
| 2 | 0 | 2 | 4 | 0 | 2 | 4 |
| 3 | 0 | 2 | 4 | 0 | 2 | 4 |
| 4 | 0 | 2 | 4 | 0 | 2 | 4 |
| 5 | 0 | 2 | 4 | 0 | 2 | 4 |

This is a non commutative groupoid with {0, 2, 4} as its subgroupoid.

**THEOREM 3.4.1:** *The groupoids $Z_n$ (0, t) $\in Z^{***}(n)$ is a P-groupoid and alternative groupoid if and only if $t^2 \equiv t \mod n$.*

*Proof:* This is left as an exercise for the reader to prove the result.

Here we do not prove many of the results as the main aim of this book is to display Smarandache groupoid notions and give central importance only to them so several proofs in groupoids in chapter 2 and chapter 3 are left for reader to prove using number theoretic notions.

So we propose the proof of these results as problem at the end of the section for the reader to prove them.

**PROBLEM 1:** Does the groupoid $Z_{11}$ (3, 0) satisfy the Moufang identity?

**PROBLEM 2:** Can $Z_{19}$ (0, 2) have a proper subgroupoids?

**PROBLEM 3:** Find all the subgroupoids in $Z_{24}$ (0, 6).

**PROBLEM 4:** Is $Z_7$ (0, 4) a normal groupoid?

**PROBLEM 5:** Is $Z_n$ (t, 0) a Bol groupoid?

**PROBLEM 6:** Will $Z_{20}$ (0, 4) be an alternative groupoid?

**PROBLEM 7:** Can $Z_n$ (t, 0) satisfy Moufang identity?

**PROBLEM 8:** Will $Z_{13}$ (0, 3) be an alternative groupoid?

**PROBLEM 9:** Find all subgroupoids and ideals if any of $Z_{12}$ (4, 0).

## 3.5 Groupoids with Identity Using $Z_n$

We saw that all the four new classes of groupoids built using $Z_n$ viz, $Z(n)$, $Z^*(n)$, $Z^{**}(n)$ and $Z^{***}(n)$ had no groupoid (which are not semigroups) with identity. Hence we have proposed in Chapter 7 an open problem about the possible construction of a groupoid which is not a semigroup using $Z_n$. The only way to obtain groupoids with identity is by adjoining an element with $Z_n$.

**DEFINITION:** *Let $Z_n = \{0, 1, 2, \ldots, n-1\}$ $n \geq 3$, $n < \infty$. Let $G = Z_n \cup \{e\}$ where $e \notin Z_n$. Define operation $*$ on $G$ by $a_i * a_i = e$ for all $a_i \in Z_n$ and $a_i * e = e * a_i = a_i$ for all $a_i \in Z_n$.*

*Finally $a_i * a_j = ta_i + ua_j \pmod{n}$ $t, u \in Z_n$. It can be checked $(G, *)$ is a groupoid with identity $e$ and the order of $G$ is $n + 1$.*

***Example 3.5.1:*** Let $G = \{Z_3(2, 2) \cup e\}$ the groupoid $(G, *)$ is given by the following table:

| * | e | 0 | 1 | 2 |
|---|---|---|---|---|
| e | e | 0 | 1 | 2 |
| 0 | 0 | e | 2 | 1 |
| 1 | 1 | 2 | e | 0 |
| 2 | 2 | 1 | 0 | e |

This is a commutative groupoid with identity.

***Example 3.5.2:*** Let $G = \{Z_n(t, u) \cup e\}$ the groupoid $(G, *)$ is a groupoid with identity then every pair $\{m, e\}$ by the very definition is a subgroupoid of G.

**THEOREM 3.5.1:** *Let $G = (Z_n(t, u) \cup \{e\})$ $e \notin Z_n$, be a groupoid with operation $*$ defined by $m * m = e$ and $m * e = e * m = m$ for all $m \in G$. Then $(G, *)$ is a groupoid with unit and has non trivial subgroups.*

*Proof:* It is clear from the fact that every pair $\{m, e\}$ is a cyclic group of order 2; Hence the claim.

**THEOREM 3.5.2:** *Let $G = (Z_n(m, m-1) \cup \{e\})$ if $n > 3$, n odd and t is a positive integer from $Z_n$ such that $(m, n) = 1$ and $(m - 1, n) = 1$ with $m < n$ and $*$ is an operation such that*

1. *$e * i = i * e = i$ for all $i \in G$*
2. *$i * i = e$ for all $i \in G$*
3. *$i * j = [mj - (m - 1)i] \bmod n$*

*for all $i, j \in Z_n$, $i \neq j$, then $G = (Z_n(m, m-1) \cup \{e\})$ is a loop which is a groupoid.*

*Proof:* Since simple number theoretic methods will yield the results, It is left for the reader to prove that $(G, *)$ is a loop.

Thus, we do not know whether there exists groupoids with identity using $Z_n$ other than the once mentioned here got by adjoining an extra element e.

Thus we see we can have groupoids with identity given by $G = (Z_n \cup \{e\})$.

**PROBLEM 1:** Find how many groupoid exist using $Z_5 \cup \{e\}$. (' $*$ ' defined as in the definition).

**PROBLEM 2:** Find the subgroupoids of (G, $*$) where $G = Z_6 \cup \{*\}$.

**PROBLEM 3:** Is (G, $*$) where $G = Z_4 (3, 2) \cup \{e\}$ a Moufang groupoid?

**PROBLEM 4:** Show (G, $*$) where $G = \{Z_7 (6, 5) \cup \{e\}\}$ is not an idempotent groupoid.

**PROBLEM 5:** Is (G, $*$) when $G = (Z_8 (7, 6) \cup \{e\})$ a P-groupoid?

**PROBLEM 6:** How many groupoids exist using $Z_{13} \cup \{e\}$?

**PROBLEM 7:** Is the groupoid given in problem 3, a right alternative groupoid?

## Supplementary Reading

1. R.H.Bruck, *A Survey of Binary Systems,* Springer Verlag, (1958).

2. W.B.Vasantha Kandasamy, *On Ordered Groupoids and its Groupoids Rings,* Jour. of Maths and Comp. Sci., Vol 9, 145-147, (1996).

3. W.B.Vasantha Kandasamy, *New classes of finite groupoids using $Z_n$*, Varahmihir Journal of Mathematical Sciences, Vol 1, 135-143, (2001).

# CHAPTER FOUR
# SMARANDACHE GROUPOIDS

In this chapter, we show the notion of SG, define and study several of its new properties. Here we work on new Smarandache substructures, identities on SGs and study them. We have given many examples of SG to make the concepts defined about them explicit.

## 4.1 Smarandache Groupoids

**DEFINITION:** *A Smarandache groupoid G is a groupoid which has a proper subset S, $S \subset G$ such that S under the operations of G is a semigroup.*

*Example 4.1.1:* Let $(G, *)$ be a groupoid given by the following table:

| * | 0 | 1 | 2 | 3 | 4 | 5 |
|---|---|---|---|---|---|---|
| 0 | 0 | 3 | 0 | 3 | 0 | 3 |
| 1 | 1 | 4 | 1 | 4 | 1 | 4 |
| 2 | 2 | 5 | 2 | 5 | 2 | 5 |
| 3 | 3 | 0 | 3 | 0 | 3 | 0 |
| 4 | 4 | 1 | 4 | 1 | 4 | 1 |
| 5 | 5 | 2 | 5 | 2 | 5 | 2 |

Clearly, $S_1 = \{0, 3\}$, $S_2 = \{1, 4\}$ and $S_3 = \{2, 5\}$ are proper subsets of G which are semigroups of G.

So $(G, *)$ is a SG.

*Example 4.1.2:* Let G be a groupoid on $Z_{10} = \{0, 1, 2, \ldots, 9\}$. Define operation $*$ on $Z_{10}$ for $1, 5 \in Z_{10} \setminus \{0\}$ by $a * b = 1a + 5b \pmod{10}$ for all $a, b \in Z_{10}$.

The groupoid is given by the following table:

| * | 0 | 1 | 2 | 3 | 4 | 5 | 6 | 7 | 8 | 9 |
|---|---|---|---|---|---|---|---|---|---|---|
| 0 | 0 | 5 | 0 | 5 | 0 | 5 | 0 | 5 | 0 | 5 |
| 1 | 1 | 6 | 1 | 6 | 1 | 6 | 1 | 6 | 1 | 6 |
| 2 | 2 | 7 | 2 | 7 | 2 | 7 | 2 | 7 | 2 | 7 |
| 3 | 3 | 8 | 3 | 8 | 3 | 8 | 3 | 8 | 3 | 8 |
| 4 | 4 | 9 | 4 | 9 | 4 | 9 | 4 | 9 | 4 | 9 |
| 5 | 5 | 0 | 5 | 0 | 5 | 0 | 5 | 0 | 5 | 0 |
| 6 | 6 | 1 | 6 | 1 | 6 | 1 | 6 | 1 | 6 | 1 |
| 7 | 7 | 2 | 7 | 2 | 7 | 2 | 7 | 2 | 7 | 2 |
| 8 | 8 | 3 | 8 | 3 | 8 | 3 | 8 | 3 | 8 | 3 |
| 9 | 9 | 4 | 9 | 4 | 9 | 4 | 9 | 4 | 9 | 4 |

Clearly, $(G, *) = Z_{10}(1, 5)$ is a SG as $S = \{0, 5\}$ is a semigroup of $Z_{10}$.

***Example 4.1.3:*** Let $Z^+ \cup \{0\}$ be the set of integers. Define ' – ' the difference between any two elements as a binary operation on $Z^+ \cup \{0\}$. $Z^+ \cup \{0\}$ is a SG as every set $\{0, n\}$; $n \in Z^+$ is a semigroup of $Z^+ \cup \{0\}$.

**DEFINITION:** *Let G be a Smarandache groupoid (SG) if the number of elements in G is finite we say G is a finite SG, otherwise G is said to be an infinite SG.*

**DEFINITION:** *Let G be a Smarandache groupoid. G is said to be a Smarandache commutative groupoid if there is a proper subset, which is a semigroup, is a commutative semigroup.*

***Example 4.1.4:*** Let $G = \{a_1, a_2, a_3\}$ be a groupoid given by the following table:

| * | $a_1$ | $a_2$ | $a_3$ |
|---|---|---|---|
| $a_1$ | $a_1$ | $a_3$ | $a_2$ |
| $a_2$ | $a_2$ | $a_1$ | $a_3$ |
| $a_3$ | $a_3$ | $a_2$ | $a_1$ |

Clearly, G is a groupoid and G is not commutative as $a_2 * a_3 \neq a_3 * a_2$ but G is a Smarandache commutative groupoid as $\{a_1\} = S$ is a commutative semigroup of G. G is not a commutative groupoid but G is a Smarandache commutative groupoid.

***Example 4.1.5:*** Let G be a groupoid given by the following table:

| * | $a_o$ | $a_1$ | $a_2$ | $a_3$ |
|---|---|---|---|---|
| $a_0$ | $a_0$ | $a_3$ | $a_2$ | $a_1$ |
| $a_1$ | $a_2$ | $a_1$ | $a_o$ | $a_3$ |
| $a_2$ | $a_0$ | $a_3$ | $a_2$ | $a_1$ |
| $a_3$ | $a_2$ | $a_1$ | $a_0$ | $a_3$ |

Clearly, G is a non-commutative groupoid but G is a Smarandache commutative groupoid as $P_1 = \{a_0\}$, $P_2 = \{a_1\}$, $P_3 = \{a_3\}$ and $P_4 = \{a_2\}$ are the only proper subsets which are semigroups of G and all of them are commutative.

This following example is a non-commutative groupoid, which is a Smarandache commutative groupoid.

*Example 4.1.6:* Let G be a groupoid given by the following table:

| * | 0 | 1 | 2 | 3 | 4 |
|---|---|---|---|---|---|
| 0 | 0 | 4 | 3 | 2 | 1 |
| 1 | 2 | 1 | 0 | 4 | 3 |
| 2 | 4 | 3 | 2 | 1 | 0 |
| 3 | 1 | 0 | 4 | 3 | 2 |
| 4 | 3 | 2 | 1 | 0 | 4 |

Now {1}, {0}, {2}, {3} and {4} are the only proper subsets, which are semigroups. Clearly G is a non-commutative groupoid but G is Smarandache commutative.

*Example 4.1.7:* Let G be a groupoid given by the following table:

| * | 0 | 1 | 2 | 3 |
|---|---|---|---|---|
| 0 | 0 | 2 | 0 | 2 |
| 1 | 2 | 0 | 2 | 0 |
| 2 | 0 | 2 | 0 | 2 |
| 3 | 2 | 0 | 2 | 0 |

This groupoid G is commutative and (G, ∗) is a SG as the proper subset $A = \{0, 2\} \subset G$ is a semigroup. Since G is a commutative groupoid G is a Smarandache commutative groupoid.

In view of the above examples we have the following theorem.

**THEOREM 4.1.1:** *Let G be a commutative groupoid if G is a Smarandache groupoid then G is a Smarandache commutative groupoid. Conversely, if G is a Smarandache commutative groupoid G need not in general be a commutative groupoid.*

*Proof:* Given G is a commutative groupoid and G is a SG so every proper subset which is a semigroup is commutative. Hence, G is a Smarandache commutative groupoid.

To prove the converse we give the following example.

Consider the groupoid (G, ∗) given in example 4.1.6, clearly (G, ∗) is a non commutative groupoid, but (G, ∗) is a Smarandache commutative groupoid. Hence the claim.

**PROBLEM 1:** Construct a groupoid of order 17, which is not commutative but Smarandache commutative.

**PROBLEM 2:** How many SGs of order 5 exist?

**PROBLEM 3:** Give an example of a Smarandache commutative groupoid of order 24 in which every semigroup is commutative.

**PROBLEM 4:** Find a Smarandache commutative groupoid of order 12.

**PROBLEM 5:** Find all semigroups of the groupoid $Z_8$ (2, 7). Is $Z_8$ (2, 7) a SG?

**PROBLEM 6:** Is $Z_{11}$ (3, 7) a SG? Justify your answer.

**PROBLEM 7:** Find all semigroups in $Z_{18}$ (2, 9).

## 4.2 Substructures in Smarandache Groupoids

In this section we extend the research to the concepts of Smarandache subgroupoids, Smarandache left (right) ideal, Smarandache ideal, Smarandache seminormal groupoid, Smarandache normal groupoid, Smarandache semiconjugate groupoid, Smarandache conjugate groupoid, Smarandache inner commutative groupoid and obtain some interesting results about them. Further, we give several examples of each definition to make the definition easy for understanding.

**DEFINITION:** *Let $(G, *)$ be a Smarandache groupoid (SG). A non-empty subset H of G is said to be a Smarandache subgroupoid if H contains a proper subset $K \subset H$ such that K is a semigroup under the operation $*$.*

**THEOREM 4.2.1:** *Every subgroupoid of a Smarandache groupoid S need not in general be a Smarandache subgroupoid of S.*

*Proof*: Let $S = Z_6 = \{0, 1, 2, 3, 4, 5\}$ set of modulo integer 6.

The groupoid $(S,*)$ is given by the following table:

| * | 0 | 1 | 2 | 3 | 4 | 5 |
|---|---|---|---|---|---|---|
| 0 | 0 | 5 | 4 | 3 | 2 | 1 |
| 1 | 4 | 3 | 2 | 1 | 0 | 5 |
| 2 | 2 | 1 | 0 | 5 | 4 | 3 |
| 3 | 0 | 5 | 4 | 3 | 2 | 1 |
| 4 | 4 | 3 | 2 | 1 | 0 | 5 |
| 5 | 2 | 1 | 0 | 5 | 4 | 3 |

The only subgroupoids of S are $A_1 = \{0, 3\}$, $A_2 = \{0, 2, 4\}$ and $A_3 = \{1, 3, 5\}$. Clearly, $A_2$ has no non-trivial semigroup as $\{0\}$ is a trivial semigroup. Thus, $A_2$ is a subgroupoid, which is not a Smarandache subgroupoid of S.

*Remark:* Just as in a group G we call the identity element of G as a trivial subgroup, we in case of groupoids built using the finite integers $Z_n$ call $\{0\}$ the trivial subgroupoid of $Z_n$ so the singleton $\{0\}$ is a trivial semigroup of $Z_n$.

**THEOREM 4.2.2:** *If a groupoid G contains a Smarandache subgroupoid then G is a SG.*

*Proof*: Let G be a groupoid. $H \subset G$ be a Smarandache subgroupoid that is H contains a proper subset $P \subset H$ such that P is a semigroup under the same operations of G. So $P \subset G$ is a semigroup. Hence the claim.

**DEFINITION:** *A Smarandache left ideal A of the Smarandache groupoid G satisfies the following conditions:*

1. *A is a Smarandache subgroupoid*
2. *$x \in G$ and $a \in A$ then $xa \in A$.*

*Similarly, we can define Smarandache right ideal. If A is both a Smarandache right ideal and Smarandache left ideal simultaneously then we say A is a Smarandache ideal.*

*Example 4.2.1:* Let $Z_6$ (4, 5) be a SG given by the following table:

| * | 0 | 1 | 2 | 3 | 4 | 5 |
|---|---|---|---|---|---|---|
| 0 | 0 | 5 | 4 | 3 | 2 | 1 |
| 1 | 4 | 3 | 2 | 1 | 0 | 5 |
| 2 | 2 | 1 | 0 | 5 | 4 | 3 |
| 3 | 0 | 5 | 4 | 3 | 2 | 1 |
| 4 | 4 | 3 | 2 | 1 | 0 | 5 |
| 5 | 2 | 1 | 0 | 5 | 4 | 3 |

Now the set $A = \{1, 3, 5\}$ is a Smarandache left ideal of $Z_6$ (4, 5) and is not a Smarandache right ideal; the subset $\{3\}$ is a semigroup of A.

*Example 4.2.2:* Let $G = Z_6$ (2, 4) be the SG given by the following table:

| * | 0 | 1 | 2 | 3 | 4 | 5 |
|---|---|---|---|---|---|---|
| 0 | 0 | 4 | 2 | 0 | 4 | 2 |
| 1 | 2 | 0 | 4 | 2 | 0 | 4 |
| 2 | 4 | 2 | 0 | 4 | 2 | 0 |
| 3 | 0 | 4 | 2 | 0 | 4 | 2 |
| 4 | 2 | 0 | 4 | 2 | 0 | 4 |
| 5 | 4 | 2 | 0 | 4 | 2 | 0 |

Clearly, G is a SG as P = {0, 3} ⊂ G is a semigroup under ∗. Now Q = {0, 2, 4} is an ideal of G but Q is not a Smarandache ideal of G for Q is not even a Smarandache subgroupoid.

We take {0} as a trivial Smarandache ideal only whenever {0} happens to be an ideal for in every groupoid $Z_n$, {0} is not in general an ideal of $Z_n$ as in the case of rings.

So whenever {0} happens to be an ideal we can consider it as a trivial one.

This example leads us to the following theorem.

**THEOREM 4.2.3:** *Let G be a groupoid. If A is a Smarandache ideal of G then A is an ideal of G. Conversely an ideal of G in general is not a Smarandache ideal of G even if G is a SG.*

*Proof*: By the very definition of Smarandache ideal A in a groupoid G; A will be an ideal of the groupoid G. To prove the converse. Consider the groupoid given in example 4.2.2. Clearly, it is proved in example 4.2.2. G = $Z_6$ is a SG and Q = {0, 2, 4} is an ideal which is not a Smarandache ideal of G. Hence the claim.

**DEFINITION:** *Let G be a SG. V be a Smarandache subgroupoid of G. We say V is a Smarandache seminormal groupoid if*

1. *aV = X for all a ∈ G*
2. *Va = Y for all a ∈ G*

*where either X or Y is a Smarandache subgroupoid of G but X and Y are both subgroupoids.*

*Example 4.2.3:* Let G = $Z_6$ (4, 5) be a groupoid given in example 4.2.1. Clearly, G is a SG as P = {3} is a semigroup. Take A = {1, 3, 5}, A is also a Smarandache subgroupoid of G.

Now aA = A is a SG and Aa = {0, 2, 4}, but, {0, 2, 4} is not a Smarandache subgroupoid of G. Hence A is a Smarandache seminormal groupoid of G.

**DEFINITION:** *Let G be a SG and V be a Smarandache subgroupoid of G. V is said to be a Smarandache normal groupoid if aV = X and Va = Y for all a ∈ G where both X and Y are Smarandache subgroupoids of G.*

**THEOREM 4.2.4:** *Every Smarandache normal groupoid is a Smarandache seminormal groupoid and not conversely.*

*Proof:* By the very definition of Smarandache normal groupoid and Smarandache seminormal groupoid we see every Smarandache normal groupoid is a Smarandache seminormal groupoid. To prove the converse we consider the groupoid given in example 4.2.1.

We have proved A = {1, 3, 5} ⊂ $Z_6$ (4, 5) is a Smarandache seminormal groupoid as aA = A but A is not a Smarandache normal groupoid as Aa = {0, 2, 4} but {0, 2, 4} is not a Smarandache subgroupoid of $Z_6$ (4, 5) so A is not a Smarandache normal groupoid. Hence a claim.

Now we will see an example of a Smarandache normal groupoid.

***Example 4.2.4:*** Let $G = Z_8$ (2, 6) be the groupoid given by the following table:

| * | 0 | 1 | 2 | 3 | 4 | 5 | 6 | 7 |
|---|---|---|---|---|---|---|---|---|
| 0 | 0 | 6 | 4 | 2 | 0 | 6 | 4 | 2 |
| 1 | 2 | 0 | 6 | 4 | 2 | 0 | 6 | 4 |
| 2 | 4 | 2 | 0 | 6 | 4 | 2 | 0 | 6 |
| 3 | 6 | 4 | 2 | 0 | 6 | 4 | 2 | 0 |
| 4 | 0 | 6 | 4 | 2 | 0 | 6 | 4 | 2 |
| 5 | 2 | 0 | 6 | 4 | 2 | 0 | 6 | 4 |
| 6 | 4 | 2 | 0 | 6 | 4 | 2 | 0 | 6 |
| 7 | 6 | 4 | 2 | 0 | 6 | 4 | 2 | 0 |

Clearly $G = Z_8$ (2, 6) is a SG for {0, 4} is a semigroup of G. Consider A = {0, 2, 4, 6} $\subset$ G, it is easily verified A is Smarandache subgroupoid. Clearly aA = A for all a $\in$ G and Aa = A for all a $\in$ G. Hence A is a Smarandache normal groupoid of G.

Now we proceed on to see when are two subgroupoids Smarandache semiconjugate with each other.

**DEFINITION:** *Let G be a SG. H and P be any two subgroupoids of G. We say H and P are Smarandache semiconjugate subgroupoids of G if*

1. *H and P are Smarandache subgroupoids of G*
2. *H = xP or Px or*
3. *P = xH or Hx for some x $\in$ G.*

***Example 4.2.5:*** Let $G = Z_8$ (2, 4) be a groupoid given in Example 4.2.4:

G is a SG as S = {0, 4} is a semigroup of G. Take P = {0, 3, 2, 4, 6} clearly P is a Smarandache subgroupoid of G. Take Q = {0, 2, 4, 6}, Q is also a Smarandache subgroupoid of G. Now 7P = Q. Hence P and Q are Smarandache semiconjugate subgroupoids of G.

**DEFINITION:** *Let G be a Smarandache groupoid. H and P be subgroupoids of G. We say H and P are Smarandache conjugate subgroupoids of G if*

1. *H and P are Smarandache subgroupoids of G*
2. *H = xP or Px and*
3. *P = xH or Hx.*

*Remark:* If one of H = xP (or Px) or P = xH (or Hx) occurs it is Smarandache semiconjugate but for it to be Smarandache conjugate we need both H = xP (or Px) and P = xH (or Hx) should occur.

***Example 4.2.6:*** Let $G = Z_{12}(1, 3)$ be the groupoid given by the following table:

| * | 0 | 1 | 2 | 3 | 4 | 5 | 6 | 7 | 8 | 9 | 10 | 11 |
|---|---|---|---|---|---|---|---|---|---|---|----|----|
| 0 | 0 | 3 | 6 | 9 | 0 | 3 | 6 | 9 | 0 | 3 | 6  | 9  |
| 1 | 1 | 4 | 7 | 10| 1 | 4 | 7 | 10| 1 | 4 | 7  | 10 |
| 2 | 2 | 5 | 8 | 11| 2 | 5 | 8 | 11| 2 | 5 | 8  | 11 |
| 3 | 3 | 6 | 9 | 0 | 3 | 6 | 9 | 0 | 3 | 6 | 9  | 0  |
| 4 | 4 | 7 | 10| 1 | 4 | 7 | 10| 1 | 4 | 7 | 10 | 1  |
| 5 | 5 | 8 | 11| 2 | 5 | 8 | 11| 2 | 5 | 8 | 11 | 2  |
| 6 | 6 | 9 | 0 | 3 | 6 | 9 | 0 | 3 | 6 | 9 | 0  | 3  |
| 7 | 7 | 10| 1 | 4 | 7 | 10| 1 | 4 | 7 | 10| 1  | 4  |
| 8 | 8 | 11| 2 | 5 | 8 | 11| 2 | 5 | 8 | 11| 2  | 5  |
| 9 | 9 | 0 | 3 | 6 | 9 | 0 | 3 | 6 | 9 | 0 | 3  | 6  |
| 10| 10| 1 | 4 | 7 | 10| 1 | 4 | 7 | 10| 1 | 4  | 7  |
| 11| 11| 2 | 5 | 8 | 11| 2 | 5 | 8 | 11| 2 | 5  | 8  |

G is a SG as the set $S = \{0, 6\}$ is a semigroup of G. $A_1 = \{0, 3, 6, 9\}$ and $A_2 = \{2, 5, 8, 11\}$ are two Smarandache subgroupoids of G for $P_1 = \{0, 6\} \subset A_1$ and $P_2 = \{8\} \subset A_2$ where $P_1$ and $P_2$ semigroups of $A_1$ and $A_2$ respectively.

Now $A_1 = 3A_2 = 3\{2, 5, 8, 11\} = \{0, 3, 6, 9\}$. Further $A_2 = 2A_1 = 2\{0, 3, 6, 9\} = \{2, 5, 8, 11\}$. Hence, $A_1$ and $A_2$ are Smarandache conjugate subgroupoids of G.

**THEOREM 4.2.5:** *Let G be a SG. If P and K be any two Smarandache subgroupoids of G which are Smarandache conjugate then they are Smarandache semiconjugate. But Smarandache semiconjugate subgroupoids need not in general be Smarandache conjugate subgroupoids of G.*

*Proof*: By the very definition of Smarandache semiconjugate subgroupoids and Smarandache conjugate subgroupoids we see every Smarandache conjugate subgroupoid is Smarandache semiconjugate subgroupoid. It is left for the reader as an exercise to construct an example of a Smarandache semiconjugate subgroupoids, which are not Smarandache conjugate subgroupoids.

**DEFINITION:** *Let G be a SG we say G is a Smarandache inner commutative groupoid if every semigroup contained in every Smarandache subgroupoid of G is commutative.*

***Example 4.2.7:*** Let G be the groupoid given by the following table:

| * | 0 | 1 | 2 | 3 |
|---|---|---|---|---|
| 0 | 0 | 3 | 2 | 1 |
| 1 | 2 | 1 | 0 | 3 |
| 2 | 0 | 3 | 2 | 1 |
| 3 | 2 | 1 | 0 | 3 |

Clearly, G is a SG as $S = \{1\}$ and $P = \{2\}$ are proper subsets of G which are semigroups under the operation $*$. Now $A_1 = \{0, 2\}$ and $A_2 = \{1, 3\}$ are subgroupoids of G.

In fact $A_1$ and $A_2$ are Smarandache subgroupoids of G but $A_1$ and $A_2$ are not commutative subgroupoids, but G is Smarandache inner commutative as the semigroups contained in them are commutative.

*Example 4.2.8:* Let $G = Z_5 (3, 3)$ be a groupoid given by the following table:

| * | 0 | 1 | 2 | 3 | 4 |
|---|---|---|---|---|---|
| 0 | 0 | 3 | 1 | 4 | 2 |
| 1 | 3 | 1 | 4 | 2 | 0 |
| 2 | 1 | 4 | 2 | 0 | 3 |
| 3 | 4 | 2 | 0 | 3 | 1 |
| 4 | 2 | 0 | 3 | 1 | 4 |

G is a commutative groupoid. G is a SG as {1}, {2}, {3} and {4} are subsets which are semigroups. So G is a Smarandache inner commutative groupoid.

*Example 4.2.9:* Let G be the groupoid given by the following table:

| * | a | b | c | d |
|---|---|---|---|---|
| a | a | c | a | c |
| b | a | c | a | c |
| c | a | c | a | c |
| d | a | c | a | c |

G is both non-commutative and non-associative for $a * b \neq b * a$ as $a * b = c$ and $b * a = a$ also $a * (b * d) \neq (a * b) * d$. For $a * (b * d) = a$ and $(a * b) * d = c$.

G is a SG as $P = \{a, c\} \subset G$ is a semigroup under the operation $*$. This groupoid has subgroupoids, which are non-commutative, but G is a inner commutative groupoid.

**THEOREM 4.2.6:** *Every Smarandache commutative groupoid is a Smarandache inner commutative groupoid and not conversely.*

*Proof:* By the very definitions of Smarandache inner commutative and Smarandache commutative we see every Smarandache inner commutative groupoid is a Smarandache commutative groupoid.

To prove the converse we construct the following example.

Let $Z_2 = \{0, 1\}$ be the prime field of characteristic 2. Consider the set of all $2 \times 2$ matrices with entries from the field $Z_2 = \{0, 1\}$ and denote it by $M_{2 \times 2}$.

$M_{2 \times 2}$ contains exactly 16 elements. By defining a non-associative and non-commutative operation 'o' we make $M_{2 \times 2}$ a Smarandache commutative groupoid which is not a Smarandache inner-commutative groupoid.

The proof is as follows:

$$M_{2\times 2} = \left\{ \begin{array}{cccc} \begin{pmatrix} 0 & 0 \\ 0 & 0 \end{pmatrix} \begin{pmatrix} 0 & 1 \\ 0 & 0 \end{pmatrix} \begin{pmatrix} 1 & 0 \\ 0 & 0 \end{pmatrix} \begin{pmatrix} 0 & 0 \\ 1 & 0 \end{pmatrix} \\ \begin{pmatrix} 0 & 0 \\ 0 & 1 \end{pmatrix} \begin{pmatrix} 1 & 1 \\ 0 & 0 \end{pmatrix} \begin{pmatrix} 1 & 0 \\ 1 & 0 \end{pmatrix} \begin{pmatrix} 0 & 1 \\ 0 & 1 \end{pmatrix} \\ \begin{pmatrix} 0 & 0 \\ 1 & 1 \end{pmatrix} \begin{pmatrix} 1 & 0 \\ 0 & 1 \end{pmatrix} \begin{pmatrix} 0 & 1 \\ 1 & 0 \end{pmatrix} \begin{pmatrix} 1 & 1 \\ 0 & 1 \end{pmatrix} \\ \begin{pmatrix} 1 & 1 \\ 1 & 0 \end{pmatrix} \begin{pmatrix} 0 & 1 \\ 1 & 1 \end{pmatrix} \begin{pmatrix} 1 & 0 \\ 1 & 1 \end{pmatrix} \begin{pmatrix} 1 & 1 \\ 1 & 1 \end{pmatrix} \end{array} \right\}.$$

$M_{2\times 2}$ is made into a groupoid by defining multiplication ' o ' as follows. If

$$A = \begin{pmatrix} a_1 & a_2 \\ a_3 & a_4 \end{pmatrix} \text{ and } B = \begin{pmatrix} b_1 & b_2 \\ b_3 & b_4 \end{pmatrix},$$

the operation 'o' defined on $M_{2\times 2}$ is

$$A \circ B = \begin{pmatrix} a_1 & a_2 \\ a_3 & a_4 \end{pmatrix} \circ \begin{pmatrix} b_1 & b_2 \\ b_3 & b_4 \end{pmatrix}$$

$$A \circ B = \begin{bmatrix} a_1 b_3 + a_2 b_1 \pmod 2 & a_1 b_4 + a_2 b_2 \pmod 2 \\ a_3 b_3 + a_4 b_1 \pmod 2 & a_3 b_4 + a_4 b_2 \pmod 2 \end{bmatrix}.$$

Clearly ($M_{2\times 2}$, o) is a groupoid. It is left for the reader to verify ' o ' is non-associative on $M_{2\times 2}$. $M_{2\times 2}$ is a SG for the set

$$A = \left\{ \begin{pmatrix} 0 & 0 \\ 0 & 0 \end{pmatrix} \begin{pmatrix} 1 & 0 \\ 0 & 0 \end{pmatrix}, \circ \right\},$$

is a semigroup as

$$\begin{pmatrix} 1 & 0 \\ 0 & 0 \end{pmatrix} \circ \begin{pmatrix} 1 & 0 \\ 0 & 0 \end{pmatrix} = \begin{pmatrix} 0 & 0 \\ 0 & 0 \end{pmatrix}.$$

Now consider

$$A_1 = \left\{ \begin{pmatrix} 0 & 0 \\ 0 & 0 \end{pmatrix} \begin{pmatrix} 1 & 0 \\ 0 & 0 \end{pmatrix} \begin{pmatrix} 0 & 1 \\ 0 & 0 \end{pmatrix} \begin{pmatrix} 1 & 1 \\ 0 & 0 \end{pmatrix}, \text{'o'} \right\},$$

is a Smarandache subgroupoid. But $A_1$ is non-commutative Smarandache subgroupoid for $A_1$ contains the non-commutative semigroup S.

$$S = \left\{ \begin{pmatrix} 0 & 0 \\ 0 & 0 \end{pmatrix}, \begin{pmatrix} 1 & 1 \\ 0 & 0 \end{pmatrix}, \begin{pmatrix} 1 & 0 \\ 0 & 0 \end{pmatrix}, o \right\}$$

such that

$$\begin{pmatrix} 1 & 0 \\ 0 & 0 \end{pmatrix} o \begin{pmatrix} 1 & 1 \\ 0 & 0 \end{pmatrix} = \begin{pmatrix} 0 & 0 \\ 0 & 0 \end{pmatrix}$$

and

$$\begin{pmatrix} 1 & 1 \\ 0 & 0 \end{pmatrix} \begin{pmatrix} 1 & 0 \\ 0 & 0 \end{pmatrix} = \begin{pmatrix} 1 & 0 \\ 0 & 0 \end{pmatrix}.$$

So ($M_{2 \times 2}$, o) is a Smarandache commutative groupoid but not a Smarandache inner commutative groupoid.

**PROBLEM 1:** Is $Z_{20}$ (3, 8) a SG? Does $Z_{20}$ (3, 8) have semigroups of order 5?

**PROBLEM 2:** Find all the subsets, which are semigroups of the groupoid $Z_{15}$ (7, 3).

**PROBLEM 3:** Can $Z_{11}$ (7, 4) have subsets, which are commutative semigroups?

**PROBLEM 4:** Find Smarandache right ideals of $Z_{15}$ (7, 8).

**PROBLEM 5:** Does $Z_{15}$ (7, 8) have Smarandache left ideals and Smarandache ideals?

**PROBLEM 6:** Find all Smarandache subgroupoids of $Z_{11}$ (3, 8).

**PROBLEM 7:** Find two Smarandache semiconjugate subgroupoids of $Z_{24}$ (3, 11).

**PROBLEM 8:** Find two Smarandache conjugate subgroupoids of $Z_{30}$ (5, 27).

**PROBLEM 9:** Does $Z_{22}$ (9, 11) have Smarandache ideals?

**PROBLEM 10:** Find Smarandache seminormal groupoids of the groupoid $Z_{22}$ (8, 13).

**PROBLEM 11:** Can $Z_{22}$ (8, 15) have Smarandache normal groupoids?

**PROBLEM 12:** Find all the Smarandache seminormal groupoids which are Smarandache semiconjugate in any groupoid G.

**PROBLEM 13:** Does there exist SGs in which all Smarandache normal subgroupoids are

a) Smarandache semiconjugate with each other
b) Smarandache conjugate with each other?

**PROBLEM 14:** Is $Z_{18}$ (3, 14) Smarandache commutative? Justify.

**PROBLEM 15:** Can the groupoid $Z_{19}$ (4, 11) be Smarandache inner commutative? Justify your answer.

### 4.3 Identities in Smarandache Groupoids

In this section we introduce some identities in SGs called Bol identity and Moufang identity, P-identity and alternative identity and illustrate them with examples.

**DEFINITION:** *A SG (G, ∗) is said to be a Smarandache Moufang groupoid if there exist $H \subset G$ such that H is a Smarandache subgroupoid satisfying the Moufang identity (xy)(zx) = (x (yz)) x for all x, y, z in H.*

**DEFINITION:** *Let S be a SG. If every Smarandache subgroupoid H of S satisfies the Moufang identity for all x, y, z in H then S is a Smarandache strong Moufang groupoid.*

*Example 4.3.1:* Let G = $Z_{10}$ (5, 6) be a groupoid given by the following table:

| ∗ | 0 | 1 | 2 | 3 | 4 | 5 | 6 | 7 | 8 | 9 |
|---|---|---|---|---|---|---|---|---|---|---|
| 0 | 0 | 6 | 2 | 8 | 4 | 0 | 6 | 2 | 8 | 4 |
| 1 | 5 | 1 | 7 | 3 | 9 | 5 | 1 | 7 | 3 | 9 |
| 2 | 0 | 6 | 2 | 8 | 4 | 0 | 6 | 2 | 8 | 4 |
| 3 | 5 | 1 | 7 | 3 | 9 | 5 | 1 | 7 | 3 | 9 |
| 4 | 0 | 6 | 2 | 8 | 4 | 0 | 6 | 2 | 8 | 4 |
| 5 | 5 | 1 | 7 | 3 | 9 | 5 | 1 | 7 | 3 | 9 |
| 6 | 0 | 6 | 2 | 8 | 4 | 0 | 6 | 2 | 8 | 4 |
| 7 | 5 | 1 | 7 | 3 | 9 | 5 | 1 | 7 | 3 | 9 |
| 8 | 0 | 6 | 2 | 8 | 4 | 0 | 6 | 2 | 8 | 4 |
| 9 | 5 | 1 | 7 | 3 | 9 | 5 | 1 | 7 | 3 | 9 |

Clearly, $Z_{10}$ (5, 6) is a SG as S = {2} is a semi-group. It can be verified $Z_{10}$ (5, 6) is a Smarandache strong Moufang groupoid as the Moufang identity (x ∗ y) ∗ (z ∗ x) = [x ∗ (y ∗ z)] ∗ x is true for

$$(x ∗ y) ∗ (z ∗ x) = (5x + 6y) ∗ (5z + 6x)$$
$$= 5(5x + 6y) + 6(5z + 6x)$$
$$= x.$$

Now $[x ∗ (y ∗ z)] ∗ x = (x ∗ [5y + 6z]) ∗ x$
$$= [5x + 6(5y + 6z)] ∗ x$$
$$= (5x + 6z) ∗ x$$
$$= 5(5x + 6z) + 6x$$

$$= 5x + 6x$$
$$= x.$$

Since x, y, z are arbitrary elements from $Z_{10}$ (5,6) we have this groupoid to be a Smarandache strong Moufang groupoid. Since every Smarandache, subgroupoid of G will also satisfy the Moufang identity.

In view of this example, we give the following nice theorem.

**THEOREM 4.3.1:** *Let S be a SG if the Moufang identity $(x * y) * (z * x) = (x * (y * z)) * x$ is true for all $x, y, z \in S$ then S is a Smarandache strong Moufang groupoid.*

*Proof:* If $(x * y) * (z * x) = (x * (y * z)) * x$ for all $x, y, z \in S$ then clearly we see every Smarandache subgroupoid of S satisfies the Moufang identity hence S is a Smarandache strong Moufang groupoid.

Now we see an example of a groupoid, which is only a Smarandache Moufang groupoid, and not a Smarandache strong Moufang groupoid.

*Example 4.3.2:* Let $G = Z_{12}$ (3, 9) be a groupoid given by the following table:

| *  | 0 | 1 | 2 | 3 | 4 | 5 | 6 | 7 | 8 | 9 | 10 | 11 |
|----|---|---|---|---|---|---|---|---|---|---|----|----|
| 0  | 0 | 9 | 6 | 3 | 0 | 9 | 6 | 3 | 0 | 9 | 6  | 3  |
| 1  | 3 | 0 | 9 | 6 | 3 | 0 | 9 | 6 | 3 | 0 | 9  | 6  |
| 2  | 6 | 3 | 0 | 9 | 6 | 3 | 0 | 9 | 6 | 3 | 0  | 9  |
| 3  | 9 | 6 | 3 | 0 | 9 | 6 | 3 | 0 | 9 | 6 | 3  | 0  |
| 4  | 0 | 9 | 6 | 3 | 0 | 9 | 6 | 3 | 0 | 9 | 6  | 3  |
| 5  | 3 | 0 | 9 | 6 | 3 | 0 | 9 | 6 | 3 | 0 | 9  | 6  |
| 6  | 6 | 3 | 0 | 9 | 6 | 3 | 0 | 9 | 6 | 3 | 0  | 9  |
| 7  | 9 | 6 | 3 | 0 | 9 | 6 | 3 | 0 | 9 | 6 | 3  | 0  |
| 8  | 0 | 9 | 6 | 3 | 0 | 9 | 6 | 3 | 0 | 9 | 6  | 3  |
| 9  | 3 | 0 | 9 | 6 | 3 | 0 | 9 | 6 | 3 | 0 | 9  | 6  |
| 10 | 6 | 3 | 0 | 9 | 6 | 3 | 0 | 9 | 6 | 3 | 0  | 9  |
| 11 | 9 | 6 | 3 | 0 | 9 | 6 | 3 | 0 | 9 | 6 | 3  | 0  |

Clearly, $A_1 = \{0, 4, 8\}$ and $A_2 = \{0, 3, 6, 9\}$ are Smarandache subgroupoids of G as $P_1 = \{0, 4\} \subset A_1$ is a semigroup and $P_2 = \{0, 6\} \subset A_2$ is a semigroup of G. $A_2$ does not satisfy Moufang identity but $A_1$ satisfies the Moufang identity. So G is only a Smarandache Moufang groupoid but not a Smarandache strong Moufang groupoid.

**THEOREM 4.3.2:** *Every Smarandache strong Moufang groupoid is a Smarandache Moufang groupoid and not conversely.*

*Proof*: By the very definition of the Smarandache strong Moufang groupoid and the Smarandache Moufang groupoid it is clear that every Smarandache strong Moufang groupoid is a Smarandache Moufang groupoid.

To prove the converse, consider the example 4.3.2 clearly that groupoid is proved to be Smarandache Moufang groupoid but it is not Smarandache strong Moufang groupoid for $x, y, z \in Z_{12}(3, 9)$ we see $(x * y) * (z * x) \neq (x * (y * z)) * x$ for

$(x * y) * (z * x)$ = $6x + 3y + 3z \pmod{12}$

Now $(x * (y * z)) * x$ = $6x + 9y + 3z \pmod{12}$.

Thus, $Z_{12}(3, 9)$ satisfies the Moufang identity for the subgroupoid $\{0, 6\}$.

***Example 4.3.3:*** Let $G = Z_{12}(3, 4)$ be a groupoid given by the following table:

| * | 0 | 1 | 2 | 3 | 4 | 5 | 6 | 7 | 8 | 9 | 10 | 11 |
|---|---|---|---|---|---|---|---|---|---|---|----|----|
| 0 | 0 | 4 | 8 | 0 | 4 | 8 | 0 | 4 | 8 | 0 | 4  | 8  |
| 1 | 3 | 7 | 11| 3 | 7 | 11| 3 | 7 | 11| 3 | 7  | 11 |
| 2 | 6 | 10| 2 | 6 | 10| 2 | 6 | 10| 2 | 6 | 10 | 2  |
| 3 | 9 | 1 | 5 | 9 | 1 | 5 | 9 | 1 | 5 | 9 | 1  | 5  |
| 4 | 0 | 4 | 8 | 0 | 4 | 8 | 0 | 4 | 8 | 0 | 4  | 8  |
| 5 | 3 | 7 | 11| 3 | 7 | 11| 3 | 7 | 11| 3 | 7  | 11 |
| 6 | 6 | 10| 2 | 6 | 10| 2 | 6 | 10| 2 | 6 | 10 | 2  |
| 7 | 9 | 1 | 5 | 9 | 1 | 5 | 9 | 1 | 5 | 9 | 1  | 5  |
| 8 | 0 | 4 | 8 | 0 | 4 | 8 | 0 | 4 | 8 | 0 | 4  | 8  |
| 9 | 3 | 7 | 11| 3 | 7 | 11| 3 | 7 | 11| 3 | 7  | 11 |
| 10| 6 | 10| 2 | 6 | 10| 2 | 6 | 10| 2 | 6 | 10 | 2  |
| 11| 9 | 1 | 5 | 9 | 1 | 5 | 9 | 1 | 5 | 9 | 1  | 5  |

This groupoid is a SG as it contains semigroups viz. $\{2\}$, $\{4\}$ and $\{10\}$. Now $Z_{12}(3, 4)$ is a Smarandache strong Bol groupoid for $((x * y) * z) * y = x * [(y * z) * y]$ for all $x, y, z \in Z_{12}(3, 4)$.

Consider,

$$((x * y) * z) * y = [(3x + 4y) * z] * y$$
$$= [9x + 12y + 4z] * y$$
$$= [9x + 4z] * y$$
$$= 27x + 12z + 4y$$
$$= 3x + 4y.$$

Now

$$x * [(y * z) * y] = x * [(3y + 4z) * y]$$
$$= x * [9y + 12z + 4y]$$
$$= x * y$$
$$= 3x + 4y.$$

Thus $Z_{12}(3, 4)$ is a Smarandache strong Bol groupoid as every triple satisfies the Bol identity so every subgroupoid of $Z_{12}(3, 4)$ which is a Smarandache subgroupoid of G will satisfy the Bol identity.

**DEFINITION:** *Let G be a groupoid, G is said to be a Smarandache Bol groupoid if G has a subgroupoid H of G such that H is a Smarandache subgroupoid and satisfies the Bol identity:*

*((x * y) * z) * y = x * ((y * z) * y) for all x, y, z ∈ H.*

**DEFINITION:** *Let G be a groupoid, G is said to be a Smarandache strong Bol groupoid if every Smarandache subgroupoid of G satisfies the Bol identity.*

**THEOREM 4.3.3:** *Let G be a SG if every triple satisfies the Bol identity then G is a Smarandache strong Bol groupoid.*

*Proof:* If ((x * y) * z) * y = x * ((y * z) * y) for all x, y, z ∈ G then we see G is a Smarandache strong Bol groupoid as every Smarandache subgroupoid of G will satisfy the Bol identity.

Now we give an example of a Smarandache Bol groupoid S.

*Example 4.3.4:* Let $Z_4$ (2, 3) be the groupoid given by the following table:

| * | 0 | 1 | 2 | 3 |
|---|---|---|---|---|
| 0 | 0 | 3 | 2 | 1 |
| 1 | 2 | 1 | 0 | 3 |
| 2 | 0 | 3 | 2 | 1 |
| 3 | 2 | 1 | 0 | 3 |

Clearly, $Z_4$ (2, 3) is a SG as the set {2} ⊂ $Z_4$ (2, 3) is a semigroup. Now A = {0, 2} is a Smarandache subgroupoid of $Z_4$ (2, 3). In fact A = {0, 2} satisfies the Bol identity (left for the reader to verify). But $Z_4$ (2, 3) does not satisfy Bol identity for all x, y, z ∈ $Z_4$ (2,3). For ((x * y) * z) * y ≠ x * [(y * z) * y],

Consider

| ((x * y) * z) * y | = | [(2x + 3y) * z] * y |
|---|---|---|
| | = | (4x + 6y + 3z) * y |
| | = | (2y + 3z) * y |
| | = | 4x + 6z + 3y |
| | = | 2z + 3y. |

Now  x * [(y * z) * y]  =  x * [(2y + 3z) * y]
                       =  x * [4y + 6z + 3y]
                       =  x * [2z + 3y]
                       =  2x + 6z + 9y
                       =  2x + 2z + y.

Clearly 2z + 3y ≠ 2x + 2z + y for all choices of x, y and z. So $Z_4$ (2,3) is not a Smarandache strong Bol groupoid.

It is only a Smarandache Bol groupoid. In view of these examples, we have a nice theorem.

**THEOREM 4.3.4:** *Every Smarandache strong Bol groupoid is a Smarandache Bol groupoid and not conversely.*

*Proof*: Now by the very definition of Smarandache strong Bol groupoid and Smarandache Bol groupoid we see every Smarandache strong Bol groupoid is a Smarandache Bol groupoid.

To prove the converse, consider the groupoid $Z_4$ (2, 3) given in example 4.3.4. It is not a Smarandache strong Bol groupoid but it is clearly a Smarandache Bol groupoid. Hence the claim.

Now we proceed on to define Smarandache P- groupoids.

**DEFINITION:** *Let G be a SG. G is said to be a Smarandache P – groupoid if G contains a proper subset $H \subset G$ where H is a Smarandache subgroupoid of G and satisfies the identity $(x * y) * x = x * (y * x)$ for all $x, y \in H$.*

**DEFINITION:** *Let G be a SG. G is said to be a Smarandache strong P-groupoid if every Smarandache subgroupoid of G is a Smarandache P-groupoid of G.*

We illustrate this by the following example.

*Example 4.3.5:* Consider the groupoid $Z_6$ (4, 3) given by the following table:

| * | 0 | 1 | 2 | 3 | 4 | 5 |
|---|---|---|---|---|---|---|
| 0 | 0 | 3 | 0 | 3 | 0 | 3 |
| 1 | 4 | 1 | 4 | 1 | 4 | 1 |
| 2 | 2 | 5 | 2 | 5 | 2 | 5 |
| 3 | 0 | 3 | 0 | 3 | 0 | 3 |
| 4 | 4 | 1 | 4 | 1 | 4 | 1 |
| 5 | 2 | 5 | 2 | 5 | 2 | 5 |

Clearly $Z_6$ (4, 3) is a SG as every singleton set i.e. {2}, {1}, {3}, {4} are subsets which are semigroups of $Z_6$ (4, 3). Now $Z_6$ (4, 3) is a Smarandache strong P-groupoid for $(x * y) * x = x * (y * x)$ for all $x, y \in Z_6$ (4, 3).

Consider $(x * y) * x$ = $(4x + 3y) * x$  
= $16x + 12y + 5x$  
= $x$.

Now $x * (y * x)$ = $x * [4y + 3x]$  
= $4x + 12y + 9x$  
= $x$.

This is true for all $x, y \in Z_6$ (4, 3). Hence, $Z_6$ (4, 3) is a Smarandache strong P-groupoid.

***Example 4.3.6:*** Consider the groupoid $Z_4$ (2, 3) given by the following table:

| * | 0 | 1 | 2 | 3 |
|---|---|---|---|---|
| 0 | 0 | 3 | 2 | 1 |
| 1 | 2 | 1 | 0 | 3 |
| 2 | 0 | 3 | 2 | 1 |
| 3 | 2 | 1 | 0 | 3 |

Clearly, $Z_4$ (2, 3) is a SG as {2} and {3} are semigroups. Now is $Z_4$ (2, 3) a Smarandache strong P-groupoid? That is; Is $(x * y) * x = x * (y * x)$?

Consider $(x * y) * x$ = $[2x + 3y] * x$
 = $4x + 6y + 3x$
 = $2y + 3x$.

Now $x * [y * x]$ = $x * [2y + 3x]$
 = $2x + 6y + 9x$
 = $3x + 2y$.

Thus, $Z_4$ (2, 3) is a Smarandache strong P-groupoid.

***Example 4.3.7:*** Consider the groupoid $Z_6$ (3, 5) given by the following table:

| * | 0 | 1 | 2 | 3 | 4 | 5 |
|---|---|---|---|---|---|---|
| 0 | 0 | 5 | 4 | 3 | 2 | 1 |
| 1 | 3 | 2 | 1 | 0 | 5 | 4 |
| 2 | 0 | 5 | 4 | 3 | 2 | 1 |
| 3 | 3 | 2 | 1 | 0 | 5 | 4 |
| 4 | 0 | 5 | 4 | 3 | 2 | 1 |
| 5 | 3 | 2 | 1 | 0 | 5 | 4 |

Now $Z_6$ (3, 5) is a SG as A = {0, 3} is a semigroup of $Z_6$ (3, 5). The only Smarandache subgroupoid of $Z_6$ (3,5) is A = {0, 3} and A satisfies P-groupoid identity.

Hence $Z_6$ (3, 5) is a Smarandache strong P-groupoid but not all elements of $Z_6$(3, 5) satisfy P-groupoid identity.

This example leads us to the following theorem.

**THEOREM 4.3.5:** *Let G be a SG. If G is a Smarandache strong P-groupoid then every pair in G need not satisfy the P-groupoid identity.*

*Proof:* We prove this only by an example. Clearly the groupoid in example 4.3.7 is a Smarandache strong P-groupoid but it can be verified that every pair in G does not satisfy the P-groupoid identity.

**THEOREM 4.3.6:** *Every Smarandache strong P-groupoid is a Smarandache P-groupoid but not conversely.*

*Proof:* By the very definition of Smarandache strong P-groupoid and Smarandache P-groupoid we see every Smarandache strong P-groupoid is a Smarandache P-groupoid. To prove the converse, consider the example. Consider the groupoid $Z_{12}$ (5, 10) given by the following table:

| *  | 0  | 1  | 2  | 3  | 4  | 5  | 6  | 7  | 8  | 9  | 10 | 11 |
|----|----|----|----|----|----|----|----|----|----|----|----|----|
| 0  | 0  | 10 | 8  | 6  | 4  | 2  | 0  | 10 | 8  | 6  | 4  | 2  |
| 1  | 5  | 3  | 1  | 11 | 9  | 7  | 5  | 3  | 1  | 11 | 9  | 7  |
| 2  | 10 | 8  | 6  | 4  | 2  | 0  | 10 | 8  | 6  | 4  | 2  | 0  |
| 3  | 3  | 1  | 11 | 9  | 7  | 5  | 3  | 1  | 11 | 9  | 7  | 5  |
| 4  | 8  | 6  | 4  | 2  | 0  | 10 | 8  | 6  | 4  | 2  | 0  | 10 |
| 5  | 1  | 11 | 9  | 7  | 5  | 3  | 1  | 11 | 9  | 7  | 5  | 3  |
| 6  | 6  | 4  | 2  | 0  | 10 | 8  | 6  | 4  | 2  | 0  | 10 | 8  |
| 7  | 11 | 9  | 7  | 5  | 3  | 1  | 11 | 9  | 7  | 5  | 3  | 1  |
| 8  | 4  | 2  | 0  | 10 | 8  | 6  | 4  | 2  | 0  | 10 | 8  | 6  |
| 9  | 9  | 7  | 5  | 3  | 1  | 11 | 9  | 7  | 5  | 3  | 1  | 11 |
| 10 | 2  | 0  | 10 | 8  | 6  | 4  | 2  | 0  | 10 | 8  | 6  | 4  |
| 11 | 7  | 5  | 3  | 1  | 11 | 9  | 7  | 5  | 3  | 1  | 11 | 9  |

Now {6} is a semigroup and A = {0, 6} is a Smarandache subgroupoid and the elements in A satisfy the identity $(x * y) * x = x * (y * x)$.

For $(x * y) * x$          =     $(5x + 10y) * x$
                               =     $25x + 50y + 10x$
                               =     $11x + 2y$.

Now $x * (y * x)$          =     $x * [5y + 10x]$
                               =     $5x + 50y + 100x$
                               =     $9x + 2y$.

Now the set {0, 6} satisfies the identity. Thus $Z_{12}$ (5, 10) is only a Smarandache P-groupoid and it is not a Smarandache strong P-groupoid.

Now we proceed on to define right (left) alternative identities in SGs.

**DEFINITION:** *Let G be a SG, G is said to be a Smarandache right alternative groupoid if a proper subset H of G where H is a Smarandache subgroupoid of G satisfies the right alternative identity; (xy) y = x (yy) for all x, y ∈ H.*

*On similar lines we define Smarandache left alternative groupoid if it satisfies the left alternative identity that is (xx) y = x (xy) for all x, y ∈ H. A SG, G is said to be Smarandache alternative groupoid if a proper subset H of G which is a Smarandache subgroupoid, simultaneously satisfies both the right and left alternative identities.*

**DEFINITION:** *Let G be SG if every Smarandache subgroupoid of G satisfies the right alternative identity, we call G a Smarandache strong right alternative groupoid. Similarly, Smarandache strong left alternative groupoid is defined.*

*G is said to be Smarandache strong alternative groupoid if every Smarandache subgroupoid of G satisfies both the right and left alternative identities simultaneously.*

We illustrate this by the following example:

*Example 4.3.8:* Let $Z_{14}$ (7, 8) be a groupoid given by the following table:

| * | 0 | 1 | 2 | 3 | 4 | 5 | 6 | 7 | 8 | 9 | 10 | 11 | 12 | 13 |
|---|---|---|---|---|---|---|---|---|---|---|----|----|----|----|
| 0 | 0 | 8 | 2 | 10 | 4 | 12 | 6 | 0 | 8 | 2 | 10 | 4 | 12 | 6 |
| 1 | 7 | 1 | 9 | 3 | 11 | 5 | 13 | 7 | 1 | 9 | 3 | 11 | 5 | 13 |
| 2 | 0 | 8 | 2 | 10 | 4 | 12 | 6 | 0 | 8 | 2 | 10 | 4 | 12 | 6 |
| 3 | 7 | 1 | 9 | 3 | 11 | 5 | 13 | 7 | 1 | 9 | 3 | 11 | 5 | 13 |
| 4 | 0 | 8 | 2 | 10 | 4 | 12 | 6 | 0 | 8 | 2 | 10 | 4 | 12 | 6 |
| 5 | 7 | 1 | 9 | 3 | 11 | 5 | 13 | 7 | 1 | 9 | 3 | 11 | 5 | 13 |
| 6 | 0 | 8 | 2 | 10 | 4 | 12 | 6 | 0 | 8 | 2 | 10 | 4 | 12 | 6 |
| 7 | 7 | 1 | 9 | 3 | 11 | 5 | 13 | 7 | 1 | 9 | 3 | 11 | 5 | 13 |
| 8 | 0 | 8 | 2 | 10 | 4 | 12 | 6 | 0 | 8 | 2 | 10 | 4 | 12 | 6 |
| 9 | 7 | 1 | 9 | 3 | 11 | 5 | 13 | 7 | 1 | 9 | 3 | 11 | 5 | 13 |
| 10 | 0 | 8 | 2 | 10 | 4 | 12 | 6 | 0 | 8 | 2 | 10 | 4 | 12 | 6 |
| 11 | 7 | 1 | 9 | 3 | 11 | 5 | 13 | 7 | 1 | 9 | 3 | 11 | 5 | 13 |
| 12 | 0 | 8 | 2 | 10 | 4 | 12 | 6 | 0 | 8 | 2 | 10 | 4 | 12 | 6 |
| 13 | 7 | 1 | 9 | 3 | 11 | 5 | 13 | 7 | 1 | 9 | 3 | 11 | 5 | 13 |

Clearly, $Z_{14}$ (7, 8) is a SG as {4} is a semigroup.

Now the subgroupoids of $Z_{14}$ (7, 8) are

| * | 0 | 2 |
|---|---|---|
| 0 | 0 | 2 |
| 2 | 0 | 2 |

| * | 0 | 4 |
|---|---|---|
| 0 | 0 | 4 |
| 4 | 0 | 4 |

| * | 0 | 8 |
|---|---|---|
| 0 | 0 | 8 |
| 8 | 0 | 8 |

All these subgroupoids are Smarandache subgroupoids. This groupoid is both right alternative and left alternative. Further $Z_{14}$ (7, 8) is a Smarandache strong alternative groupoid.

Now we give an example of a Smarandache alternative groupoid, which is a Smarandache strong alternative groupoid where only SG satisfies the identity but all pairs do not satisfy the alternative identity.

*Example 4.3.9:* Let $Z_{12}$ (1,6) be a groupoid given by the following table:

| * | 0 | 1 | 2 | 3 | 4 | 5 | 6 | 7 | 8 | 9 | 10 | 11 |
|---|---|---|---|---|---|---|---|---|---|---|----|----|
| 0 | 0 | 6 | 0 | 6 | 0 | 6 | 0 | 6 | 0 | 6 | 0 | 6 |
| 1 | 1 | 7 | 1 | 7 | 1 | 7 | 1 | 7 | 1 | 7 | 1 | 7 |
| 2 | 2 | 8 | 2 | 8 | 2 | 8 | 2 | 8 | 2 | 8 | 2 | 8 |
| 3 | 3 | 9 | 3 | 9 | 3 | 9 | 3 | 9 | 3 | 9 | 3 | 9 |
| 4 | 4 | 10 | 4 | 10 | 4 | 10 | 4 | 10 | 4 | 10 | 4 | 10 |
| 5 | 5 | 11 | 5 | 11 | 5 | 11 | 5 | 11 | 5 | 11 | 5 | 11 |
| 6 | 6 | 0 | 6 | 0 | 6 | 0 | 6 | 0 | 6 | 0 | 6 | 0 |
| 7 | 1 | 7 | 1 | 7 | 1 | 7 | 1 | 7 | 1 | 7 | 1 | 7 |
| 8 | 8 | 2 | 8 | 2 | 8 | 2 | 8 | 2 | 8 | 2 | 8 | 2 |
| 9 | 9 | 3 | 9 | 3 | 9 | 3 | 9 | 3 | 9 | 3 | 9 | 3 |
| 10 | 10 | 4 | 10 | 4 | 10 | 4 | 10 | 4 | 10 | 4 | 10 | 4 |
| 11 | 11 | 5 | 11 | 5 | 11 | 5 | 11 | 5 | 11 | 5 | 11 | 5 |

Now $Z_{12}$ (16) is a SG as the singleton sets {2}, {4} and {8} are semigroups. Now A = (10, 4) is the Smarandache subgroupoid given by the table:

| * | 4 | 10 |
|---|---|----|
| 4 | 4 | 4 |
| 10 | 10 | 10 |

Now it is left for the reader to verify that both the right and left alternative identities are satisfied by A, but not by all elements of G. Hence, $Z_{12}$ (1,6) is a Smarandache alternative groupoid which is evident from the fact the subgroupoid $A_2 = \{1, 7\}$ does not satisfy the alternative identities. Thus two examples give no doubt only a Smarandache strong groupoid but are formulated differently.

**THEOREM 4.3.7:** *Every Smarandache strong alternative groupoid is a Smarandache alternative groupoid and not conversely.*

*Proof*: Clearly, by the very definitions we see every Smarandache strong alternative groupoid is a Smarandache alternative groupoid. To prove the converse we give an example.

From example 4.3.9 $Z_{12}$ (1, 6) is only a Smarandache alternative groupoid, which is not a Smarandache strong alternative groupoid.

**PROBLEM 1:** Give an example of a SG of order 10, which is a Smarandache Moufang groupoid.

**PROBLEM 2:** Can a groupoid of order 11 be a Smarandache strong Bol groupoid?

**PROBLEM 3:** Check whether $Z_{13}$ (3, 9) is a Smarandache alternative groupoid.

**PROBLEM 4:** Is $Z_7$ (3, 4) a Smarandache P-groupoid?

**PROBLEM 5:** Can $Z_8$ (5, 3) be a Smarandache Bol groupoid?

**PROBLEM 6:** Can a groupoid satisfy Bol identity, Moufang identity simultaneously? (or) Can a groupoid be both Smarandache strong Moufang groupoid and Smarandache strong Bol groupoid? Illustrate them by example.

**PROBLEM 7:** Is $Z_{30}$ (7, 3) a Smarandache strong alternative groupoid?

**PROBLEM 8:** Give an example of a Smarandache strong alternative groupoid, which is not a Smarandache strong P-groupoid.

**PROBLEM 9:** Find a Smarandache strong Bol groupoid G in which only every Smarandache subgroupoid satisfies the Bol identity but every triple in G does not satisfy the Bol identity.

**PROBLEM 10:** Give an example of a Smarandache Moufang groupoid, which is not a Smarandache strong Moufang groupoid.

**PROBLEM 11:** Give an example of Smarandache Bol groupoid, which is not a Smarandache strong Bol groupoid.

**PROBLEM 12:** Is $Z_{19}$ (3, 12) a Smarandache Bol groupoid?

**PROBLEM 13:** Is $Z_{19}$ (4, 15) a Smarandache Bol groupoid?

**PROBLEM 14:** Can $Z_{19}$ (17, 2) be a Smarandache Bol groupoid?

**PROBLEM 15:** Compare the groupoids in Problems 12, 13 and 14 and derive some interesting results about them.

**PROBLEM 16:** Obtain conditions on t and u so that $Z_n$ (t, u) is a Smarandache alternative groupoid.

**PROBLEM 17:** Prove $Z_8$ (3, 6) cannot be Smarandache strong P-groupoid.

**PROBLEM 18:** Can $Z_9$ (2,4) be a Smarandache strong right alternative groupoid?

**PROBLEM 19:** Prove $Z_{12}$ (3, 6) cannot be a Smarandache strong Moufang groupoid.

## 4.4 More Properties on Smarandache Groupoids

In this section we introduce more properties about SGs like Smarandache direct product of groupoids, Smarandache homomorphism of groupoids and give some examples of them.

**DEFINITION:** *Let $G_1, G_2, \ldots, G_n$ be n-groupoids. We say $G = G_1 \times G_2 \times \ldots \times G_n$ is a Smarandache direct product of groupoids if G has a proper subset H of G which is a semigroup under the operations of G. It is important to note that each $G_i$ need not be a SG* *for in that case G will obviously be a SG. Here we take any set of groupoids and find their direct product.*

***Example 4.4.1:*** Let $G_1 = Z_{12}(8, 4)$ and $G_2 = Z_4(2, 3)$ be two groupoids. Now $G = G_1 \times G_2$ is a direct product. It is left for the reader to verify G is a Smarandache direct product.

***Example 4.4.2:*** Let $G = G_1 \times G_1 \times G_1$ where $G_1 = Z_{12}(8, 4)$. Prove G is a Smarandache direct product.

**DEFINITION:** *Let $(G, *)$ and $(G', o)$ be any two SGs. A map $\phi$ from $(G, *)$ to $(G', o)$ is said to be a SG homomorphism if $\phi : A \to A'$ is a semigroup homomorphism where $A \subset G$ and $A' \subset G'$ are semigroups of G and G' respectively. We have $\phi(a * b) = \phi(a) \circ \phi(b)$ for all a, b ∈ A. We say a SG homomorphism is a SG isomorphism if $\phi$ restricted to the semigroups is a semigroup isomorphism.*

*Remark 1*: If G and G' be two SGs which are isomorphic we need not have |G| = |G'|.

*Remark 2*: Two SGs may be homomorphic for some semigroups in them and be isomorphic for some other set of semigroups. Yet, these SGs become isomorphic.

***Example 4.4.3:*** Let $Z_4(2, 3)$ and $Z_6(4, 5)$ be two given SGs given by the following tables:

Table for $Z_4(2, 3)$

| * | 0 | 1 | 2 | 3 |
|---|---|---|---|---|
| 0 | 0 | 3 | 2 | 1 |
| 1 | 2 | 1 | 0 | 3 |
| 2 | 0 | 3 | 2 | 1 |
| 3 | 2 | 1 | 0 | 3 |

Clearly, $A = \{3\} \subset Z_4(2,3)$ is a SG. The table for groupoid $Z_6(4, 5)$ is given below.

| * | 0 | 1 | 2 | 3 | 4 | 5 |
|---|---|---|---|---|---|---|
| 0 | 0 | 5 | 4 | 3 | 2 | 1 |
| 1 | 4 | 3 | 2 | 1 | 0 | 5 |
| 2 | 2 | 1 | 0 | 5 | 4 | 3 |
| 3 | 0 | 5 | 4 | 3 | 2 | 1 |
| 4 | 4 | 3 | 2 | 1 | 0 | 5 |
| 5 | 2 | 1 | 0 | 5 | 4 | 3 |

Clearly, $A' = \{3\} \subset Z_6(4, 5)$ is a SG. We see $Z_6(4, 5) \cong Z_4(2, 3)$ under $\phi$ which maps 3 to 3. Though the groupoids are of order 4 and 6 yet they are Smarandache isomorphic.

**PROBLEM 1:** Give an example of a direct product of groupoids, which is not a Smarandache direct product.

**PROBLEM 2:** Can $Z_7(3, 4)$ be SG isomorphic with $Z_{11}(3, 4)$? Justify your answer.

**PROBLEM 3:** Can $Z_{13}$ (3, 2) be SG isomorphic with $Z_{13}$ (6, 3)?

**PROBLEM 4:** Give a SG homomorphism between $Z_{19}$ (3, 7) and $Z_{19}$ (3, 6).

**PROBLEM 5:** Is $Z_n$ (t, u) SG isomorphic with $Z_n$(u, t)?

**PROBLEM 6:** Find conditions on t and u in problem 5 so that for all n, $Z_n$ (t, u) is SG isomorphic with $Z_n$ (u, t).

**PROBLEM 7:** Let G = $Z_3$ (1, 2) × $Z_7$ (3, 4) × $Z_{11}$ (4, 5). Prove G is a Smarandache direct product. Find all the semigroups of G.

**PROBLEM 8:** Does there exist a SG isomorphism between $Z_8$ (4,6) and $Z_9$ (3,6)?

**PROBLEM 9:** Construct a SG homomorphism between $Z_6$ (2, 3) and $Z_{12}$ (2, 3).

## 4.5 Smarandache Groupoids with Identity

In this section we introduce finite SGs with identity and leave it as an open problem in chapter seven of finding infinite SGs with identity. When we say SGs with identity we only mean the identity should simultaneously be both left and right identity.

*Example 4.5.1:* The groupoid (G, ∗) given by the following table is a SG.

| ∗ | e | $a_0$ | $a_1$ | $a_2$ | $a_3$ | $a_4$ |
|---|---|---|---|---|---|---|
| e | e | $a_0$ | $a_1$ | $a_2$ | $a_3$ | $a_4$ |
| $a_0$ | $a_0$ | e | $a_2$ | $a_4$ | $a_1$ | $a_3$ |
| $a_1$ | $a_1$ | $a_2$ | e | $a_1$ | $a_3$ | $a_0$ |
| $a_2$ | $a_2$ | $a_4$ | $a_1$ | e | $a_0$ | $a_2$ |
| $a_3$ | $a_3$ | $a_1$ | $a_3$ | $a_0$ | e | $a_4$ |
| $a_4$ | $a_4$ | $a_3$ | $a_0$ | $a_2$ | $a_4$ | e |

(G, ∗) is a SG with identity e. Every pair {e, $a_i$}; i = 1, 2, 3, 4 are semigroups in G. Now consider the following two groupoids built using the set G = {e, $a_0$, $a_1$, $a_2$, $a_3$, $a_4$}.

*Example 4.5.2:* Consider the groupoid {G, o} given by the following table:

| o | e | $a_0$ | $a_1$ | $a_2$ | $a_3$ | $a_4$ |
|---|---|---|---|---|---|---|
| e | e | $a_0$ | $a_1$ | $a_2$ | $a_3$ | $a_4$ |
| $a_0$ | $a_0$ | e | $a_4$ | $a_3$ | $a_2$ | $a_1$ |
| $a_1$ | $a_1$ | $a_2$ | e | $a_0$ | $a_4$ | $a_3$ |
| $a_2$ | $a_2$ | $a_4$ | $a_3$ | e | $a_1$ | $a_0$ |
| $a_3$ | $a_3$ | $a_1$ | $a_0$ | $a_4$ | e | $a_2$ |
| $a_4$ | $a_4$ | $a_3$ | $a_2$ | $a_1$ | $a_0$ | e |

Clearly, this is a SG with identity 'e' and every pair $\{a_i, e\}$ for $i = 0, 1, 2, 3, 4$ are semigroups. It is still interesting to see that (G, o) is a loop. Now using the same set G we define $*$ on G as follows.

*Example 4.5.3:* The groupoid $(G, *)$ is given by the following table:

| * | e | $a_0$ | $a_1$ | $a_2$ | $a_3$ | $a_4$ |
|---|---|---|---|---|---|---|
| e | e | $a_0$ | $a_1$ | $a_2$ | $a_3$ | $a_4$ |
| $a_0$ | $a_0$ | e | $a_1$ | $a_2$ | $a_3$ | $a_4$ |
| $a_1$ | $a_1$ | $a_2$ | e | $a_4$ | $a_0$ | $a_1$ |
| $a_2$ | $a_2$ | $a_4$ | $a_0$ | e | $a_2$ | $a_3$ |
| $a_3$ | $a_3$ | $a_1$ | $a_2$ | $a_3$ | e | $a_0$ |
| $a_4$ | $a_4$ | $a_3$ | $a_4$ | $a_0$ | $a_1$ | e |

This is also a SG with identity but this is not a loop. Thus we can get SGs with identity only by adjoining ' e ' with them; some SGs with identity happens to be loops. Here we ask the reader to find examples of SGs of infinite order using the set of positive integers $Z^+$. Clearly any operation $*$ on $Z^+$ such that $a * b = ta + ub$ where t and u are distinct elements of $Z^+$ does not make $(Z^+, *)$ into a SG for we are not in a position to find a proper subset which is a semigroup. Thus finding SGs of infinite order seems to be a very difficult task. Hence we have proposed open problems in Chapter 7 about SGs of infinite order.

**PROBLEM 1:** How many SGs with identity of order three and four exist?

**PROBLEM 2:** Give an example of a SG with identity of order 7.

**PROBLEM 3:** Give an example of a SG of order 5, which is non commutative.

**PROBLEM 4:** Find an example of a SG of order 6, which is commutative.

**PROBLEM 5:** Prove if $G = Z^+ \cup \{0\}$ or $Q+ \cup \{0\}$ or $R+ \cup \{0\}$ is a SG with operation ' – ' difference.

## Supplementary Reading

# CHAPTER FIVE
# SMARANDACHE GROUPOIDS USING $Z_n$

In this chapter we introduce conditions for the 4 new classes of groupoids Z (n), $Z^*(n)$, $Z^{**}(n)$ and $Z^{***}(n)$ built using the modulo integers $Z_n$ to be SGs. All groupoids in these classes need not in general be SGs. We obtain conditions on n and on the pair t, u for the groupoid $Z_n$ (t, u) to be a SG. Many interesting results for the SGs built using $Z_n$ to be Bol, Moufang, alternative etc. are obtained.

## 5.1. Smarandache Groupoids in Z (n)

In this section, we completely devote our study on SGs from the class of groupoids Z (n) defined in chapter 3. We prove only those groupoids in Z (n) which have proper subsets that are semigroups will be SGs, other groupoids will not be SGs. We start with an example of a SG in Z (n).

*Example 5.1.1:* Let $Z_7$ (5, 3) ∈ Z (7) be the groupoid given by the following table:

| * | 0 | 1 | 2 | 3 | 4 | 5 | 6 |
|---|---|---|---|---|---|---|---|
| 0 | 0 | 3 | 6 | 2 | 5 | 1 | 4 |
| 1 | 5 | 1 | 4 | 0 | 3 | 6 | 2 |
| 2 | 3 | 6 | 2 | 5 | 1 | 4 | 0 |
| 3 | 1 | 4 | 0 | 3 | 6 | 2 | 5 |
| 4 | 6 | 2 | 5 | 1 | 4 | 0 | 3 |
| 5 | 4 | 0 | 3 | 6 | 2 | 5 | 1 |
| 6 | 2 | 5 | 1 | 4 | 0 | 3 | 6 |

This is a SG as {1}, {2}, {4}, {3}, {5} and {6} are subsets of $Z_7$ (5, 3) which are semigroups.

*Example 5.1.2:* Consider the groupoid $Z_4$ (3, 2) ∈ Z (4) given by the following table:

| * | 0 | 1 | 2 | 3 |
|---|---|---|---|---|
| 0 | 0 | 2 | 0 | 2 |
| 1 | 3 | 1 | 3 | 1 |
| 2 | 2 | 0 | 2 | 0 |
| 3 | 1 | 3 | 1 | 3 |

This groupoid is also a SG as {1}, {2} and {3} are semigroups under ∗.

These examples leads us to the following interesting theorem.

**THEOREM 5.1.1:** *Let $Z_n (t, u)$ be a groupoid in $Z (n)$, $n > 5$. $Z_n (t, u)$ is a SG if $(t, u) = 1$ and $t \neq u$ and $t + u \equiv 1 \pmod{n}$.*

*Proof:* Given $Z_n (t, u)$ is a groupoid from the class $Z (n)$ such that $t + u \equiv 1 \pmod{n}$. To show $Z_n (t, u)$ has semigroups consider for any $m \in Z_n \setminus \{0\}$, $m * m \equiv mt + mu \pmod{n} = m (t + u)$ as $(t + u) \equiv 1 \pmod{n}$ we have $m * m = m$ thus we have $\{m\}$ to be a semigroup.

This is true for all $m \in Z_n$. Hence $Z_n (t, u)$ is a Smarandache semigroup if $t + u \equiv 1 \pmod{n}$.

Now we are interested in finding when $Z_n (t, u) \in Z (n)$ is not a SG.

*Example 5.1.3:* Let $Z_5 (1, 3) \in Z (5)$ be given by the following table:

| * | 0 | 1 | 2 | 3 | 4 |
|---|---|---|---|---|---|
| 0 | 0 | 3 | 1 | 4 | 2 |
| 1 | 1 | 4 | 2 | 0 | 3 |
| 2 | 2 | 0 | 3 | 1 | 4 |
| 3 | 3 | 1 | 4 | 2 | 0 |
| 4 | 4 | 2 | 0 | 3 | 1 |

This is not a SG as it has no proper subset, which are semigroups under the operation ∗.

*Example 5.1.4:* Consider the groupoid. $Z_5 (2, 1) \in Z (5)$ given by the following table:

| * | 0 | 1 | 2 | 3 | 4 |
|---|---|---|---|---|---|
| 0 | 0 | 1 | 2 | 3 | 4 |
| 1 | 2 | 3 | 4 | 0 | 1 |
| 2 | 4 | 0 | 1 | 2 | 3 |
| 3 | 1 | 2 | 3 | 4 | 0 |
| 4 | 3 | 4 | 0 | 1 | 2 |

This is also not a SG as it has no proper subsets, which are semigroups.

***Example 5.1.5:*** Let $Z_9$ (5, 3) $\in$ Z (9) be a groupoid of order 9 defined by the following table:

| * | 0 | 1 | 2 | 3 | 4 | 5 | 6 | 7 | 8 |
|---|---|---|---|---|---|---|---|---|---|
| 0 | 0 | 3 | 6 | 0 | 3 | 6 | 0 | 3 | 6 |
| 1 | 5 | 8 | 2 | 5 | 8 | 2 | 5 | 8 | 2 |
| 2 | 1 | 4 | 7 | 1 | 4 | 7 | 1 | 4 | 7 |
| 3 | 6 | 0 | 3 | 6 | 0 | 3 | 6 | 0 | 3 |
| 4 | 2 | 5 | 8 | 2 | 5 | 8 | 2 | 5 | 8 |
| 5 | 7 | 1 | 4 | 7 | 1 | 4 | 1 | 4 | 7 |
| 6 | 3 | 6 | 0 | 3 | 6 | 0 | 3 | 6 | 0 |
| 7 | 8 | 2 | 5 | 8 | 2 | 5 | 8 | 2 | 5 |
| 8 | 4 | 7 | 1 | 4 | 7 | 1 | 4 | 7 | 1 |

This groupoid has only subgroupoids {0, 3, 6} and {1, 2, 4, 5, 7, 8} but $Z_9$ (5, 3) has no proper subsets which are semigroups. Thus, $Z_9$ (5, 3) is not a SG.

***Example 5.1.6:*** Consider the groupoid $Z_8$ (1, 6) $\in$ Z (8) given by the following table:

| * | 0 | 1 | 2 | 3 | 4 | 5 | 6 | 7 |
|---|---|---|---|---|---|---|---|---|
| 0 | 0 | 6 | 4 | 2 | 0 | 6 | 4 | 2 |
| 1 | 1 | 7 | 5 | 3 | 1 | 7 | 5 | 3 |
| 2 | 2 | 0 | 6 | 4 | 2 | 0 | 6 | 4 |
| 3 | 3 | 1 | 7 | 5 | 3 | 1 | 7 | 5 |
| 4 | 4 | 2 | 0 | 6 | 4 | 2 | 0 | 6 |
| 5 | 5 | 3 | 1 | 7 | 5 | 3 | 1 | 7 |
| 6 | 6 | 4 | 2 | 0 | 6 | 4 | 2 | 0 |
| 7 | 7 | 5 | 3 | 1 | 7 | 5 | 3 | 1 |

This is a SG as {4} $\subset Z_8$ (1, 6) is a semigroup.

***Example 5.1.7:*** Consider the groupoid $Z_{10}$ (1, 2) $\in$ Z (10) given by the following table:

| * | 0 | 1 | 2 | 3 | 4 | 5 | 6 | 7 | 8 | 9 |
|---|---|---|---|---|---|---|---|---|---|---|
| 0 | 0 | 2 | 4 | 6 | 8 | 0 | 2 | 4 | 6 | 8 |
| 1 | 1 | 3 | 5 | 7 | 9 | 1 | 3 | 5 | 7 | 9 |
| 2 | 2 | 4 | 6 | 8 | 0 | 2 | 4 | 6 | 8 | 0 |
| 3 | 3 | 5 | 7 | 9 | 1 | 3 | 5 | 7 | 9 | 1 |
| 4 | 4 | 6 | 8 | 0 | 2 | 4 | 6 | 8 | 0 | 2 |
| 5 | 5 | 7 | 9 | 1 | 3 | 5 | 7 | 9 | 1 | 3 |
| 6 | 6 | 8 | 0 | 2 | 4 | 6 | 8 | 0 | 2 | 4 |
| 7 | 7 | 9 | 1 | 3 | 5 | 7 | 9 | 1 | 3 | 5 |
| 8 | 8 | 0 | 2 | 4 | 6 | 8 | 0 | 2 | 4 | 6 |
| 9 | 9 | 1 | 3 | 5 | 7 | 9 | 1 | 3 | 5 | 7 |

This is a SG as $\{5\} \subset Z_{10}(1, 2)$ is a semigroup.

All these examples leads to the following theorem and some interesting open problems.

**THEOREM 5.1.2:** *Let $Z_{2p}(1, 2) \in Z(2p)$ be a groupoid (p is an odd prime). Then $Z_{2p}(1, 2)$ is a SG.*

*Proof:* Now to show $Z_{2p}(1, 2)$ is a SG we have to obtain a proper subset in $Z_{2p}$, which is a semigroup.

Now take $p \in Z_{2p}$ $p * p \equiv p + 2p \pmod{2p} \equiv p \pmod{2p}$. Thus, $\{p\}$ is a semigroup so $Z_{2p}(1, 2)$ is a SG.

**COROLLARY 5.1.3:** *Let $Z_{3p}(1, 3)$ be a groupoid in $Z(3p)$ ($p \neq 3$ and p is a prime). Then $Z_{3p}(1, 3)$ is a SG.*

*Proof:* Consider the subset $p \in Z_{3p}(1, 3)$; $p * p = p + 3p \equiv p \pmod{3p}$. Clearly, $\{p\}$ is a proper subset, which is a semigroup; hence, $Z_{3p}(1, 3)$ is a SG of order 3p.

**COROLLARY 5.1.4:** *Let $(1, p_2)$ and $Z_{p_1 p_2}(1, p_1)$ be groupoids in $Z(p_1 p_2)$ where $p_1$ and $p_2$ are two distinct prime. Then $Z_{p_1 p_2}(1, p_2)$ and $Z_{p_1 p_2}(1, p_1)$ are SGs.*

*Proof:* Take $p_1 \in Z_{p_1 p_2}(1, p_2)$. Now $p_1 * p_1 \equiv p_1 + p_1 p_2 \pmod{p_1 p_2} \equiv p_1 \pmod{p_1 p_2}$. So the set $\{p_1\}$ is a semigroup. Hence $Z_{p_1 p_2}(1, p_2)$ is a SG. Similarly $Z_{p_1 p_2}(1, p_1)$ is a SG as $p_2 \in Z_{p_1 p_2}$ is such that $p_2 * p_2 \equiv p_2 + p_1 p_2 \pmod{p_1 p_2} \equiv p_2 \pmod{p_1 p_2}$. So $\{p_2\}$ is a semigroup and $Z_{p_1 p_2}(1, p_1)$ is a SG.

**COROLLARY 5.1.5:** *Let $Z_n(1, p)$ be a groupoid in $Z(n)$ where p is a prime and p / n. Then $Z_n(1, p)$ is a SG.*

*Proof:* Let $n / p = m$ where $n = mp$. Now $m \in Z_n$ and $m * m \equiv m + pm \pmod{n = pm} \equiv m \pmod{n}$. Thus, the singleton set $\{m\}$ is a proper subset, which is a semigroup. Hence $Z_n(1, p)$ is a SG. This result leads to the following interesting theorem.

**THEOREM 5.1.6:** *The class of groupoids $Z(n)$ has SGs for all $n > 3$ and $n \neq 5$.*

*Proof:* Given $Z_n = \{0, 1, 2, \ldots, n-1\}$; $n > 3$ and $n \neq 5$.

$Z_n(t, u)$ with $t \neq u$, $(t, u) = 1$ $t, u \in Z_n \setminus \{0\}$. Now for $n = 4$. $Z_4(1, 2) \in Z(4)$, is a SG as $\{2\}$ is a semigroup, so $Z_4(1, 2)$ is a SG. Now for $n = 5$ we do not have any SG, this is left for the reader to verify. When $n \geq 6$ we see $Z(n)$ has SGs under two possibilities:

1. When we chose $t + u \equiv 1 \pmod{n}$ then $Z_n(t, u)$ is a SG.

2. Similarly when n is non prime that is a composite number and if p / n then $Z_n$ (1. p) is a SG.

Thus the class of groupoids Z (n) when n > 3 and n ≠ 5 has SGs.

*Note*: When n = 3 it is left for the reader to verify Z (3) has no SGs.

**THEOREM 5.1.7:** *The SG $Z_n$ (t, u) ∈ Z (n) with t + u ≡ 1 (mod n) is a Smarandache idempotent groupoid.*

*Proof:* Clearly x * x = x for all x ∈ $Z_n$ as x * x = tx + ux ≡ (t + u) x (mod n) as t + u ≡ 1 (mod n) we have x * x = x in $Z_n$. Hence the claim.

**THEOREM 5.1.8:** *The SG $Z_n$ (t, u) ∈ Z (n) with t + u ≡ 1(mod n) is a Smarandache P-groupoid if and only if $t^2 \equiv t$ (mod n) and $u^2 \equiv u$ (mod n).*

*Proof:* Clearly we have proved if t + u ≡ 1 (mod n), then $Z_n$ (t, u) is a SG in Z(n). To show $Z_n$ (t, u) is a Smarandache P-groupoid we have to prove (x * y) * x = x * (y * x) for all x, y ∈ $Z_n$(t,u). Now to show $Z_n$ (t, u) to be SG we have to prove (x * y) * x = x * (y * x).

Consider (x * y) * x = (tx + uy) * x = $t^2$x + tuy + ux. Now x * (y * x) = x * (ty + ux) = tx + tuy + $u^2$x. Now $t^2$x + tuy + ux = tx + tuy + $u^2$x if and only if $t^2 \equiv t$ (mod n) and $u^2 \equiv u$ (mod n).

*Example 5.1.8:* Let $Z_{12}$ (4, 9) ∈ Z (12). Clearly, $Z_{12}$ (4, 9) is a SG, which is also a Smarandache P-groupoid by theorem 5.1.8. $Z_{12}$ (4, 9) is defined by x * y = 4x + 9y (mod 12). It is left for the reader to prove $Z_{12}$ (4, 9) is a Smarandache P-groupoid.

**THEOREM 5.1.9:** *The groupoid $Z_n$ (t, u) ∈ Z (n) with t + u ≡ 1 (mod n) is a Smarandache alternative groupoid if and only if $t^2 \equiv t$ (mod n) and $u^2 \equiv u$ (mod n).*

*Proof:* $Z_n$ (t, u) when t + u ≡ 1 (mod n) is a SG. Clearly, (x * y) * y = x * (y * y) and (x * x) * y = x * (x * y) if and only if $t^2 \equiv t$ (mod n) and $u^2 \equiv u$ (mod n). Thus, we see from this theorem example 5.1.8 the groupoid $Z_{12}$ (4, 9) is also a Smarandache alternative groupoid.

**THEOREM 5.1.10:** *The groupoid $Z_n$ (t, u) ∈ Z (n) with t + u ≡1 (mod n) is a SG. This groupoid is a Smarandache strong Bol groupoid if and only if $t^3 \equiv t$ (mod n) and $u^2 \equiv u$ (mod n).*

*Proof:* We know $Z_n$ (t, u) ∈ Z (n) is a SG when t + u ≡ 1 (mod n). This groupoid is a Smarandache strong Bol groupoid if and only if $t^3 \equiv t$ (mod n) and $u^2 \equiv u$ (mod n). It essential to prove ((x * y) * z) * x = x * ((y * z) * x) for all x, y, z ∈ $Z_m$ (t, u).

((x * y) * z) * x     =     [(tx + uy) * z] * x
                       =     ($t^2$x + tuy + uz) * x
                       =     $t^3$x + $t^2$uy + tuz + ux.

Now, $x * ((y * z) * x) = x * [(ty + uz) * x]$
$= x * [t^2y + tuz + ux]$
$= tx + t^2uy + tu^2z + u^2x$.

$((x * y) * z) * x \equiv x * ((y * z) * x) \pmod{n}$ if and only if $t^3 \equiv t \pmod{n}$ and $u^2 \equiv u \pmod{n}$. Hence the claim.

*Example 5.1.9:* Let $Z_6 (3, 4) \in Z(6)$ be the SG given by the following table:

| * | 0 | 1 | 2 | 3 | 4 | 5 |
|---|---|---|---|---|---|---|
| 0 | 0 | 4 | 2 | 0 | 4 | 2 |
| 1 | 3 | 1 | 5 | 3 | 1 | 5 |
| 2 | 0 | 4 | 2 | 0 | 4 | 2 |
| 3 | 3 | 1 | 5 | 3 | 1 | 5 |
| 4 | 0 | 4 | 2 | 0 | 4 | 2 |
| 5 | 3 | 1 | 5 | 3 | 1 | 5 |

This groupoid is a Smarandache strong Bol groupoid as $4^3 \equiv 4 \pmod 6$ and $3^2 \equiv 3 \pmod 6$, which is the condition for a groupoid in $Z(n)$ to be a Smarandache strong Bol groupoid.

**THEOREM 5.1.11:** *Let $Z_n(t, u) \in Z(n)$ with $t + u \equiv 1 \pmod n$ be a SG. $Z_n(t, u)$ is a Smarandache strong Moufang groupoid if and only if $t^2 \equiv t \pmod n$ and $u^2 \equiv u \pmod n$.*

Proof: The Moufang identity is $(x * y) * (z * x) = (x * (y * z)) * x$ for all $x, y, z \in Z_n$.

Now $(x * y) * (z * x) = (tx + uy) * (tz + ux) = t^2x + tuy + tux + u^2x$.

Further $(x * (y * z)) * x = (x * (ty + uz)) * x$
$= (tx + tuy + u^2z) * x$
$= t^2x + t^2uy + tu^2z + ux$.

Now $t^2x + tuy + tuz + u^2x \equiv t^2x + t^2uy + tu^2z + ux \pmod n$ if and only if $u^2 \equiv u \pmod n$ and $t^2 \equiv t \pmod n$. Hence the claim.

Thus we see the class of groupoids $Z(n)$ contain SGs and also Smarandache strong Bol groupoids, Smarandache strong Moufang groupoids, Smarandache idempotent groupoids and Smarandache strong P-groupoids.

**PROBLEM 1:** Is $Z_8 (3, 5)$ a SG?

**PROBLEM 2:** Show $Z_{11} (5, 7) \in Z (11)$ is a SG. Does $Z_{11} (5, 7)$ have Smarandache subgroupoids?

**PROBLEM 3:** Does $Z (19)$ have SGs, which are Smarandache strong Bol groupoids?

**PROBLEM 4:** Give an example of a SG, which is not a Smarandache strong Moufang groupoid in Z (16).

**PROBLEM 5:** Find all the SGs in the class of groupoids Z (24).

**PROBLEM 6:** Find all SGs, which are Smarandache strong P-groupoids in the class of groupoids Z (18).

**PROBLEM 7:** Find all the Smarandache subgroupoids of $Z_{16}$ (9, 8).

**PROBLEM 8:** Does $Z_{27}$ (11, 17) have non trivial Smarandache subgroupoids?

**PROBLEM 9:** Find all SGs in Z (22).

**PROBLEM 10:** Is $Z_{22}$ (10, 13) a Smarandache Moufang groupoid? Justify!

**PROBLEM 11:** Prove $Z_{17}$ (11, 7) is a SG.

**PROBLEM 12:** Can Z (5) have SGs?

**PROBLEM 13:** Which groupoids in Z (14) are SGs?

**PROBLEM 14:** Prove $Z_{25}$ (11, 15) is a SG.

**PROBLEM 15:** Find all SGs in Z (15).

**PROBLEM 16:** Find all Smarandache P-groupoids in Z (9).

**PROBLEM 17:** Find all Smarandache strong Bol groupoids in Z (27).

**PROBLEM 18:** Can Z (8) have Smarandache strong Moufang groupoid?

**PROBLEM 19:** Find a Smarandache P- groupoid in Z (121).

## 5.2 Smarandache Groupoids in $Z^*(n)$

In this section, we study the conditions for the class of groupoids in $Z^*(n)$ to have SGs. Throughout this section we will be interested in studying those groupoids in $Z^*(n)$ which are not in Z (n) that is only groupoids in $Z^*(n) \setminus Z (n)$.

*Example 5.2.1*: Let $Z_5$ (2, 4) ∈ $Z^*$ (5). $Z_5$ (2, 4) is a SG given by the following table:

| * | 0 | 1 | 2 | 3 | 4 |
|---|---|---|---|---|---|
| 0 | 0 | 4 | 3 | 2 | 1 |
| 1 | 2 | 1 | 0 | 4 | 3 |
| 2 | 4 | 3 | 2 | 1 | 0 |

| 3 | 1 | 0 | 4 | 3 | 2 |
|---|---|---|---|---|---|
| 4 | 3 | 2 | 1 | 0 | 4 |

Clearly, this is a SG as {1}, {2}, {3} and {4} are subsets which are semigroups.

It is important and interesting to note that the class of groupoids $Z^*(3)$ and $Z^*(4)$ are such that $Z^*(3) = Z(3)$ and $Z^*(4) = Z(4)$. But $Z^*(5)$ has only two SGs given by $Z_5(2,4)$ and $Z_5(4, 2)$.

As in the case of groupoids in $Z(n)$, the groupoids in $Z^*(n)$ are SGs when $t + u \equiv 1$ (mod n).

***Example 5.2.2:*** Let $Z_6(2, 4) \in Z^*(6)$ be the groupoid given by the following table:

| * | 0 | 1 | 2 | 3 | 4 | 5 |
|---|---|---|---|---|---|---|
| 0 | 0 | 4 | 2 | 0 | 4 | 2 |
| 1 | 2 | 0 | 4 | 2 | 0 | 4 |
| 2 | 4 | 2 | 0 | 4 | 2 | 0 |
| 3 | 0 | 4 | 2 | 0 | 4 | 2 |
| 4 | 2 | 0 | 4 | 2 | 0 | 4 |
| 5 | 4 | 2 | 0 | 4 | 2 | 0 |

This is a SG as the set {0, 3} is a semigroup of $Z_6(2, 4)$.

Now yet another interesting example of a SG is given.

***Example 5.2.3:*** Consider the groupoid $Z_8(2, 6) \in Z^*(8)$ given by the following table:

| * | 0 | 1 | 2 | 3 | 4 | 5 | 6 | 7 |
|---|---|---|---|---|---|---|---|---|
| 0 | 0 | 6 | 4 | 2 | 0 | 6 | 4 | 2 |
| 1 | 2 | 0 | 6 | 4 | 2 | 0 | 6 | 4 |
| 2 | 4 | 2 | 0 | 6 | 4 | 2 | 0 | 6 |
| 3 | 6 | 4 | 2 | 0 | 6 | 4 | 2 | 0 |
| 4 | 0 | 6 | 4 | 2 | 0 | 6 | 4 | 2 |
| 5 | 2 | 0 | 6 | 4 | 2 | 0 | 6 | 4 |
| 6 | 4 | 2 | 0 | 6 | 4 | 2 | 0 | 6 |
| 7 | 6 | 4 | 2 | 0 | 6 | 4 | 2 | 0 |

This is a SG for the subset {0, 4} is a semigroup. Hence, $Z_8(2, 6)$ is a SG.

***Example 5.2.4:*** The groupoid $Z_8(3, 6)$ is also a SG from the class of groupoids $Z^*(8)$. Further the class $Z^*(8)$ has $Z_8(2,7)$, $Z_8(4, 5)$, $Z_8(3, 6)$, $Z_8(2, 6)$, $Z_8(4, 6)$, $Z_8(6, 2)$, $Z_8(6, 4)$ and $Z_8(6, 3)$ are SGs and the only SGs in $Z^*(8) \setminus Z(8)$ are $Z_8(2, 6)$, $Z_8(4, 6)$, $Z_8(6, 2)$, $Z_8(6, 4)$, $Z_8(3, 6)$ and $Z_8(6, 3)$. Thus there are 6 SGs in $Z^*(8) \setminus Z(8)$.

***Example 5.2.5:*** Let $Z^*(12)$ be the class of groupoids. We see $Z_{12}$ (3, 9) is a SG for it contains P = {0, 6} as a subset which is a semigroup. $Z_{12}$ (4, 6) is a SG as Q = {4} and P = {0, 6} are subsets which are semigroups.

$Z_{12}$ (4, 8) is a SG as the subset {0, 6} is a semigroup.
$Z_{12}$ (5, 10) is a SG as the subset {6} is a semigroup.
$Z_{12}$ (8, 10) is a SG as {0, 6} is a semigroup.
$Z_{12}$ (2, 4) is a SG as {0, 6} is a proper subset which is a semigroup.
$Z_{12}$ (2, 6) is a SG as {0, 6} is a proper subset which is a semigroup.
$Z_{12}$ (2, 8) is a SG as {0, 6} is a proper subset which is a semigroup.
$Z_{12}$ (2, 10) is a SG as {0, 6} is a proper subset which is a semigroup.
$Z_{12}$ (3, 6) is a SG as {0, 6} is a proper subset which is a semigroup.
$Z_{12}$ (3, 9) is a SG as {0, 6} is a proper subset which is a semigroup.
$Z_{12}$ (4, 10) is a SG as {0, 6} is a proper subset which is a semigroup.
$Z_{12}$ (6, 8) is a SG as {0, 6} is a proper subset which is a semigroup.
$Z_{12}$ (6, 10) is a SG as {0, 6} is a proper subset which is a semigroup.
$Z_{12}$ (6, 9) is a SG as {6} is a proper subset which is a semigroup.

Thus this class $Z^*(12) \setminus Z(12)$ has only 28 groupoids which are SGs.

This leads to open problems given in Chapter 7.

**THEOREM 5.2.1:** *Let $Z_n(t, u) \in Z^*(n) \setminus Z(n)$. If $t + u \equiv 1 \pmod{n}$. Then $Z_n(t, u)$ is a SG.*

*Proof*: If $t + u \equiv 1 \pmod{n}$ then for every $x \in Z_n$. $x * x = tx + ux \equiv (t + u) x \equiv x \pmod{n}$ thus every singleton is a semigroup. Hence $Z_n(t,u) \in Z^*(n) \setminus Z(n)$ is a SG if $t + u \equiv 1 \pmod{n}$.

**THEOREM 5.2.2:** *The SG $Z_n(t, u) \in Z^*(n) \setminus Z(n)$ with $t + u \equiv 1 \pmod{n}$ is a Smarandache strong Moufang groupoid, Smarandache strong P-groupoid, Smarandache idempotent groupoid, Smarandache strong Bol groupoid and Smarandache strong alternative groupoid if and only if $t^2 \equiv t \pmod{n}$ and $u^2 \equiv u \pmod{n}$.*

*Proof:* It can be easily verified that the SG $Z_n(t,u) \in Z^*(n) \setminus Z(n)$ with $t + u \equiv 1 \pmod{n}$ satisfies the following identities:

$$
\begin{aligned}
x * x &\equiv x \pmod{n} \\
(x * y) * y &\equiv x * (y * x) \pmod{n} \\
(x * x) * y &\equiv x * (x * y) \pmod{n} \\
(x * y) * x &\equiv x * (y * x) \pmod{n} \\
(x * y) * (z * x) &\equiv (x * (y * z)) * x \pmod{n} \\
((x * y) * z) * y &\equiv x * ((y * z) * y) \pmod{n}
\end{aligned}
$$

for all $x, y, z \in Z_n$ if and only if $t^2 \equiv t \pmod{n}$ and $u^2 \equiv u \pmod{n}$. Hence the claim.

Thus, we see not all SGs in $Z^*(n) \setminus Z(n)$ satisfy all the identities.

**PROBLEM 1:** Prove $Z^*(5)$ has no SGs.

**PROBLEM 2:** Does $Z^*(7) \setminus Z(7)$ have SGs?

**PROBLEM 3:** Is $Z_{11}(3, 9)$ a SG? Justify your answer.

**PROBLEM 4:** Is $Z_{13}(2, 8)$ a Smarandache strong Bol groupoid?

**PROBLEM 5:** Prove $Z_{13}(6, 8)$ is SG but it is not a Smarandache right alternative groupoid.

**PROBLEM 6:** Is $Z_{19}(6, 12)$ a SG? Justify.

**PROBLEM 7:** Prove $Z_{16}(8, 6)$ is a SG.

**PROBLEM 8:** Prove $Z_{16}(t, u) \in Z^*(16) \setminus Z(16)$ with $t + u \equiv 1 \pmod{16}$ can never occur.

**PROBLEM 9:** Find a SG in $Z_{36}(t, u) \in Z^*(36) \setminus Z(36)$ with $t + u \equiv 1 \pmod{36}$.

**PROBLEM 10:** Find all SGs in $Z^*(11)$.

**PROBLEM 11:** Find all SGs in $Z^*(18) \setminus Z(18)$.

**PROBLEM 12:** Prove $Z_{16}(2, 8)$ is a SG.

**PROBLEM 13:** Which has more number of SGs $Z^*(8)$ or $Z^*(7)$?

**PROBLEM 14:** Find a SG in $Z^*(49)$.

**PROBLEM 15:** Find all SGs in $Z^*(50)$.

**PROBLEM 16:** Find all Smarandache strong Bol groupoid in $Z^*(24)$.

**PROBLEM 17:** Find a Smarandache strong Moufang groupoid in $Z^*(17) \setminus Z(17)$.

**PROBLEM 18:** Find all SGs in $Z^*(20) \setminus Z(20)$.

**PROBLEM 19:** Find all the SGs in $Z^*(22) \setminus Z(22)$. Which class $Z^*(20)$ or $Z^*(22)$ have more number of SGs?

## 5.3 Smarandache Groupoids in $Z^{**}(n)$

In this section, we study yet another new class of groupoids $Z^{**}(n)$ to contain SGs. We are interested in only finding SGs in $Z^{**}(n) \setminus Z^*(n)$ as we have completely elaborated about SGs in $Z^*(n)$. Thus we will characterize and obtain some interesting properties about SGs in $Z^{**}(n) \setminus Z^*(n)$.

We just recall the definition of groupoids in $Z^{**}(n)$.

**DEFINITION:** *Let $Z_n = \{0, 1, 2, \ldots, n-1\}$, $n \geq 3$; $n - 1 < \infty$. Define operation $*$ on $Z_n$ by $a * b = ta + ub \pmod{n}$ where t can be equal to u also, $t, u \in Z_n \setminus \{0\}$. We denote this collection of groupoids for varying u and t by $Z^{**}(n)$.*

*Example 5.3.1:* $Z_6 (2, 2)$ is a groupoid in the class $Z^{**}(6)$ given by the following table:

| * | 0 | 1 | 2 | 3 | 4 | 5 |
|---|---|---|---|---|---|---|
| 0 | 0 | 2 | 4 | 0 | 2 | 4 |
| 1 | 2 | 4 | 0 | 2 | 4 | 0 |
| 2 | 4 | 0 | 2 | 4 | 0 | 2 |
| 3 | 0 | 2 | 4 | 0 | 2 | 4 |
| 4 | 2 | 4 | 0 | 2 | 4 | 0 |
| 5 | 4 | 0 | 2 | 4 | 0 | 2 |

This is also a SG as {4} is a subset, which is a semigroup under the operation $*$.

*Example 5.3.2:* Let $Z_9 (4, 4)$ be the groupoid given by the following table:

| * | 0 | 1 | 2 | 3 | 4 | 5 | 6 | 7 | 8 |
|---|---|---|---|---|---|---|---|---|---|
| 0 | 0 | 4 | 8 | 3 | 7 | 2 | 6 | 1 | 5 |
| 1 | 4 | 8 | 3 | 7 | 2 | 6 | 1 | 5 | 0 |
| 2 | 8 | 3 | 7 | 2 | 6 | 1 | 5 | 0 | 4 |
| 3 | 3 | 7 | 2 | 6 | 1 | 5 | 0 | 4 | 8 |
| 4 | 7 | 2 | 6 | 1 | 5 | 0 | 4 | 8 | 3 |
| 5 | 2 | 6 | 1 | 5 | 0 | 4 | 8 | 3 | 7 |
| 6 | 6 | 1 | 5 | 0 | 4 | 8 | 3 | 7 | 2 |
| 7 | 1 | 5 | 0 | 4 | 8 | 3 | 7 | 2 | 6 |
| 8 | 5 | 0 | 4 | 8 | 3 | 7 | 2 | 6 | 1 |

This is a commutative groupoid but this has no proper subsets, which are semigroups so $Z_9 (4, 4)$ is not a SG in $Z^{**}(9)$.

We see $Z^{**}(n) \setminus Z^*(n)$ contains SGs. The number of groupoids in $Z^{**}(n) \setminus Z^*(n)$ are (n – 1) for a given n and all of them are commutative. Now consider the following two groupoids in $Z^{**}(3)$.

*Example 5.3.3:* The 2 groupoids in $Z^{**}(3)$ are given by the following tables:
Table for $Z_3 (1, 1)$

| * | 0 | 1 | 2 |
|---|---|---|---|
| 0 | 0 | 1 | 2 |
| 1 | 1 | 2 | 0 |
| 2 | 2 | 0 | 1 |

Clearly, $Z_3 (1, 1)$ is a semigroup.

Table for $Z_3 (2, 2)$

| * | 0 | 1 | 2 |
|---|---|---|---|
| 0 | 0 | 1 | 2 |
| 1 | 1 | 2 | 0 |
| 2 | 2 | 0 | 1 |

The groupoid $Z_3$ (2, 2) is a SG as the proper subset {1} and {2} are semigroups. It is interesting to see the 2 groupoids in $Z^{**}(3)$ only one is a SG other being a semigroup. But in the class $Z^{**}(4)$ we see the groupoids $Z_4$ (1, 1) is a semigroup and $Z_4$ (2, 2) and $Z_4$ (3, 3) are SGs. It is left for the reader to verify the fact.

Since all semigroups are trivially SGs we are interested in studying only (n – 2) of the groupoids in $Z_n$ (m, m) given by $Z_n$ (2, 2), $Z_n$ (3, 3), … , $Z_n$ (n – 1, n – 1) as $Z_n$ (1, 1) is a semigroup. Here also it is proper to mention when n is a prime $Z_p$ (1, 1) has no proper subset which is a semigroup as $Z_p$ (1, 1) has no proper subsemigroups. So we can only say $Z_p$ (1, 1) is a trivial SG as $Z_p$ (1, 1) is itself a semigroup and $Z_p$(1, 1) may or may not have proper subsemigroups.

**THEOREM 5.3.1:** *Let p be a prime* $Z_p\left(\frac{p+1}{2}, \frac{p+1}{2}\right)$ *is a SG in $Z^{**}(p)$.*

*Proof*: Now $Z_p\left(\frac{p+1}{2}, \frac{p+1}{2}\right)$ is a SG as

$$x * x \quad = \quad \frac{(p+1)}{2}x + \frac{(p+1)}{2}x$$
$$= \quad (p+1)x \ [\bmod p]$$
$$= \quad x \ (\bmod p).$$

Hence all singleton sets are proper subsets and they are semigroups. So $Z_p\left(\frac{p+1}{2}, \frac{p+1}{2}\right)$ is a SG.

**COROLLARY 5.3.2:** $Z_p\left(\frac{p+1}{2}, \frac{p+1}{2}\right)$ *is a Smarandache idempotent groupoid of $Z^{**}(p)$.*

*Proof*: Follows by the very definition of $Z^{**}(p)$ and theorem 5.3.1.

**COROLLARY 5.3.3:** *Let n be odd. Then* $Z_n\left(\frac{n+1}{2}, \frac{n+1}{2}\right)$ *is a SG in $Z^{**}(n)$.*

*Proof*: It is left for the reader to prove.

**COROLLARY 5.3.4:** *Let $n$ be even and $\frac{n}{2}$ be such that $\left(\frac{n}{2}\right)^2 \equiv 0 \pmod{n}$. Then $Z_n\left(\frac{n+1}{2}, \frac{n+1}{2}\right)$ in $Z^{**}(n)$ is a SG.*

*Proof*: Now the subset $(0, n/2)$ is a semigroup. Hence the claim.

*Example 5.3.4:* Consider $Z_{12}$ (6, 6) $\in Z^{**}(12)$; $6^2 \equiv 0 \pmod{12}$. Clearly, the table for the subset $\{0, 6\}$ is a semigroup.

| * | 0 | 6 |
|---|---|---|
| 0 | 0 | 0 |
| 6 | 0 | 0 |

*Example 5.3.5:* Consider $Z_{14}$ (7, 7) now $7^2 \equiv 7 \pmod{14}$. The subset $\{0, 7\}$ is given by the following table:

| * | 0 | 7 |
|---|---|---|
| 0 | 0 | 7 |
| 7 | 7 | 0 |

is a semigroup. Hence, $Z_{14}$ (7, 7) is a SG in $Z^{**}(14)$.

**THEOREM 5.3.5:** *Let $Z_n$ (m, m) be a groupoid such that $n$ is even $m^2 \equiv m \pmod{n}$ and $m + m \equiv 0 \pmod{n}$. Then $Z_n$ (m, m) is a SG of $Z^{**}(n)$.*

*Proof*: Using the fact $n$ is even $m^2 \equiv m \pmod{n}$ and $m + m \equiv 0 \pmod{n}$ we see the subset $\{0, m\}$ given by the following table is a semigroup.

| * | 0 | m |
|---|---|---|
| 0 | 0 | m |
| m | m | 0 |

Thus $Z_n$ (m, n) in $Z^{**}(n)$ is a SG.

*Example 5.3.6:* Let $Z_{18}$ (9, 9) in $Z^{**}(18)$ is a SG as $\{0, 9\}$ is given by the following table:

| * | 0 | 9 |
|---|---|---|
| 0 | 0 | 9 |
| 9 | 9 | 0 |

is a semigroup.

**THEOREM 5.3.6:** *Let $Z_n$ (m, m) be a groupoid in $Z^{**}(n)$. $Z_n$ (m, m) is a SG only if $m + m \equiv 1 \pmod{n}$.*

*Proof*: Now in $Z_n$ (m, m) $\in Z^{**}$(n) we have m + m $\equiv$ 1 (mod n). To show $Z_n$ (m, m) is a SG we have to obtain a proper subset in $Z_n$ (m, m) which is a semigroup. Now take any element r $\in Z_n$.

$$r * r \equiv rm + rm \pmod{n}$$
$$\equiv r(m + m) \pmod{n}$$
$$\equiv 1 \pmod{n}$$

Thus $r * r \equiv r \pmod{n}$ so {r} is a semigroup in $Z_n$ (m, m), so $Z_n$ (m, m) is a SG in $Z^{**}$(n).

**THEOREM 5.3.7:** *The SG $Z_n$ (m, m) in $Z^{**}$(n) where m + m $\equiv$ 1 (mod n) and $m^2 \equiv m$ (mod n) is*

1. *Smarandache idempotent groupoid.*
2. *Smarandache strong P-groupoid.*
3. *Smarandache strong Bol groupoid.*
4. *Smarandache strong Moufang groupoid.*
5. *Smarandache strong alternative groupoid.*

*Proof*: The proof is left as an exercise to the reader as it involves only simple number theoretic techniques.

**PROBLEM 1:** Find all SGs in $Z^{**}$(18) \ $Z^*$(18).

**PROBLEM 2:** Find all SGs in $Z^{**}$(19) \ $Z^*$(19).

**PROBLEM 3:** Find Smarandache strong Bol groupoids in $Z^{**}$(22).

**PROBLEM 4:** Find Smarandache strong Moufang groupoids in $Z^{**}$(23).

**PROBLEM 5:** Does $Z^{**}$(22) or $Z^{**}$(23) have more number of SGs? Justify.

**PROBLEM 6:** Find the SGs in $Z^{**}$(42) \ $Z^*$(42) and in $Z^{**}$(42) \ Z (42).

**PROBLEM 7:** Is $Z_{19}$ (10, 10) a SG? Justify.

**PROBLEM 8:** Is $Z_{19}$ (7, 7) a SG? Justify.

**PROBLEM 9:** Find whether the groupoid $Z_{22}$ (13, 13) is a SG.

**PROBLEM 10:** Find whether all the groupoids in $Z^{**}$(25) \ $Z^*$(25) are SGs.

**PROBLEM 11:** Find a SG in $Z^{**}$(26). Find a groupoid in $Z^{**}$(26) \ $Z^*$(26), which is not a SG.

**PROBLEM 12:** Does there exist atleast 3 groupoids in $Z^{**}$(81) \ $Z^*$(81) which are SGs?

**PROBLEM 13:** Find atleast 2 groupoids in $Z^{**}(16) \setminus Z^{*}(16)$, which are Smarandache Bol groupoids.

**PROBLEM 14:** Find a Smarandache strong Moufang groupoid in $Z^{**}(32) \setminus Z^{*}(32)$.

**PROBLEM 15:** Does there exist a Smarandache strong P-groupoid in $Z^{**}(27)$?

**PROBLEM 16:** Find all the SGs in $Z^{**}(35) \setminus Z^{*}(35)$.

**PROBLEM 17:** Find a SG in $Z^{**}(29)$.

**PROBLEM 18:** Does there exist a SG which is not Bol or Moufang or P-groupoid or idempotent groupoid in $Z^{**}(37) \setminus Z^{*}(37)$?

**PROBLEM 19:** Find SGs in $Z^{**}(10)$.

**PROBLEM 20:** Find SGs in $Z^{**}(11)$.

**PROBLEM 21:** Which class $Z^{**}(10)$ or $Z^{**}(11)$ have more number of SGs?

**PROBLEM 22:** Find all the SGs in $Z^{**}(30) \setminus Z^{*}(30)$.

**PROBLEM 23:** Will $Z^{**}(29)$ have more than one SG?

**PROBLEM 24:** Can $Z^{**}(28)$ have more than 2 SGs?

**PROBLEM 25:** Which has more SGs $Z^{**}(28)$ or $Z^{**}(29)$?

## 5.4 Smarandache Groupoids in $Z^{***}(n)$

In this section, we find which of the groupoids in $Z^{***}(n)$ are SGs. In our study, we concentrate SGs only in the class $Z^{***}(n) \setminus Z^{**}(n)$. We obtain some interesting results about them. We just recall the definition of the class of groupoids in $Z^{***}(n)$.

**DEFINITION:** *Let $Z_n = \{0, 1, 2, \ldots, n-1\}$, $n \geq 3$, $n < \infty$. Define $*$ for $a, b \in Z_n$ by $a * b = ta + ub \pmod{n}$ where $t, u \in Z_n$. Here we permit the possibility of $t = 0$ or $u = 0$ also. Thus the class of groupoids in $Z^{***}(n)$ is such that $Z^{**}(n) \subset Z^{***}(n)$.*

Further, we have $Z(n) \subset Z^{*}(n) \subset Z^{**}(n) \subset Z^{***}(n)$. Here in this section we study only the SGs in $Z^{***}(n) \setminus Z^{**}(n)$.

*Example 5.4.1:* Let $Z_5(3,0)$ be the groupoid given by the following table:

| * | 0 | 1 | 2 | 3 | 4 |
|---|---|---|---|---|---|
| 0 | 0 | 0 | 0 | 0 | 0 |
| 1 | 3 | 3 | 3 | 3 | 3 |
| 2 | 1 | 1 | 1 | 1 | 1 |
| 3 | 4 | 4 | 4 | 4 | 4 |
| 4 | 2 | 2 | 2 | 2 | 2 |

This is not a SG. This groupoid is both non-associative and non-commutative.

We have in the class of groupoids $Z^{***}(n)$ the groupoids $Z_n$ (1, 0) and $Z_n$ (0, 1) to be semigroups so we are not interested in this section to study these groupoids.

***Example 5.4.2:*** Consider $Z_6$ (2, 0) $\in Z^{***}(6)$ given by the following table:

| * | 0 | 1 | 2 | 3 | 4 | 5 |
|---|---|---|---|---|---|---|
| 0 | 0 | 0 | 0 | 0 | 0 | 0 |
| 1 | 2 | 2 | 2 | 2 | 2 | 2 |
| 2 | 4 | 4 | 4 | 4 | 4 | 4 |
| 3 | 0 | 0 | 0 | 0 | 0 | 0 |
| 4 | 2 | 2 | 2 | 2 | 2 | 2 |
| 5 | 4 | 4 | 4 | 4 | 4 | 4 |

The subset {0, 3} is a semigroup. Hence, $Z_6$ (2, 0) is a SG.

**THEOREM 5.4.1:** *The groupoid $Z_n$ (2, 0) $\in Z^{***}(n)$ where n = 2m is a SG.*

*Proof*: The set {0, m) is a semigroup given by the table:

| * | 0 | m |
|---|---|---|
| 0 | 0 | 0 |
| m | 0 | 0 |

Hence, $Z_n$ (2, 0) is a SG.

**COROLLARY 5.4.2:** *The groupoid $Z_n$ (0, 2) in $Z^{***}(n)$ where n = 2m is a SG.*

*Proof:* Left for the reader to verify.

**THEOREM 5.4.3:** *The groupoid $Z_n$ (m, 0) $\in Z^{***}(n)$ where n = 2m is a SG.*

*Proof*: The subset {0, 2} is a semigroup; given by the following table:

| * | 0 | 2 |
|---|---|---|
| 0 | 0 | 0 |
| 2 | 0 | 0 |

**COROLLARY 5.4.4:** *The groupoid $Z_n$ (0, m) in $Z^{***}(n)$ where n = 2m is a SG.*

*Proof:* For the subset {0, 2} is given by the table:

| * | 0 | 2 |
|---|---|---|
| 0 | 0 | 0 |
| 2 | 0 | 0 |

is a semigroup. Hence $Z_n (m, 0)$ is a SG.

**THEOREM 5.4.5:** *Let $Z^{***}(n)$ be the class of groupoids. Let p be a prime such that, p / n. then $Z_n (p, 0)$ is a SG.*

*Proof*: Now to show $Z_n (p, 0)$ is a SG we have to find a subset which is a semigroup. Consider $\{0, n/p\}$ the subset given by the following table:

| * | 0 | n/p |
|---|---|-----|
| 0 | 0 | 0 |
| n/p | 0 | 0 |

This is a semigroup. Hence $Z_n (p, 0)$ is a SG.

**COROLLARY 5.4.6:** *The groupoid $Z_n (0, p) \in Z^{***}(n)$ is a SG if p / m and p is a prime.*

*Proof:* It is easily verified on similar lines.

**COROLLARY 5.4.7:** $Z_n (n/p, 0)$ *is a SG.*

*Proof*: Consider the set $\{0, p\}$ in $Z_n (n/p, 0)$. Clearly, this set is a semigroup as evident from the table:

| * | 0 | p |
|---|---|---|
| 0 | 0 | 0 |
| p | 0 | 0 |

Thus $Z_n (p, 0)$ is a SG.

From this, we get following theorem.

**THEOREM 5.4.8:** *Let $Z_n = \{0, 1, 2, ... , n – 1\}$ , $n > 3$, $n < \infty$. Suppose $n = p_1, ... , p_m$ where each $p_i$ is a distinct prime. Then the class $Z^{***}(n)$ has atleast $m + mC_2 + ... + mC_{m-1}$ SGs.*

*Proof*: Clearly we see for each prime $p_i$ we get m distinct SGs as the sets $\{0, n/p_i\}$ is a semigroup. Now we get atleast $mC_2$ distinct SGs as the sets $\{0, n/p_ip_j\}$, $i \neq j$ and $i,j \in \{1, 2, ..., m\}$ is a semigroup. (We have $mC_2$ such distinct sets collected from $p_1, p_2, ... , p_m$). Similarly we get $mC_3$ distinct SGs for the set $\{0, n / p_ip_jp_k\}$, $i \neq j$, $j \neq k$ and $i \neq k$ is a semigroup. $i, j, k \in \{1, 2, ... , m\}$. Continuing this process we see atleast $mC_{m-1}$ distinct SGs as the set $\left\{0, n/p_1p_2...\hat{p}_ip_{i+1}...p_m\right\}$ is a semigroup. $\hat{p}_i$ means it does not occur in the product $p_1 p_2 p_3 ... \hat{p}_i p_{i+1} ... p_m$. Thus if $n = p_1, p_2 ... p_m$ where m primes are distinct, we have atleast $m + mC_2 + ...+ mC_{m-1}$ SGs in $Z^{***}(n)$.

**COROLLARY 5.4.9:** *Let $Z_n = \{0, 1, 2, \ldots, n-1\}$ if $n = p_1 p_2 \ldots p_m$ where $p_i$ are distinct primes then $Z^{***}(n)$ has atleast $2(m + mC_1 + mC_2 + mC_{m-1})$ groupoids.*

*Proof:* We know if $Z_n$ (0, t) is a SG then $Z_n$ (t, 0) is also a SG. Using this principle and theorem 5.4.8 we see $Z^{***}(n)$ has atleast $2(m + mC_1 + mC_2 + \ldots + mC_{m-1})$ SGs.

**COROLLARY 5.4.10:** *Let $Z^{***}(n)$ be the class of groupoids with $n = p_1^{\alpha_1}, \ldots, p_m^{\alpha_m}$ where $p_i$ are primes, then $Z^{***}(n)$ has SGs.*

*Proof:* $Z_n(p_i^{\alpha_i}, 0)$ is a SG. For the set $(0, n/p_i^{\alpha_i})$ is a semigroup given by the following table:

| * | 0 | $n/p_i^{\alpha_i}$ |
|---|---|---|
| 0 | 0 | 0 |
| $n/p_i^{\alpha_i}$ | 0 | 0 |

Hence the claim.

*Example 5.4.3:* Let $Z^{***}(30)$ be a groupoid. Clearly $Z_{30}$ (3, 0), $Z_{30}$ (2, 0), $Z_{30}$ (5, 0), $Z_{30}$ (6, 0), $Z_{30}$ (10, 0) and $Z_{30}$ (15, 0) are SGs in $Z^{***}(30)$. Thus, there are 12 SGs in $Z^{***}(30)$. It is important to note that these are not only the SGs.

*Example 5.4.4:* Let $Z_6$ (3, 0) belong to $Z^{***}(6)$ given by the following table:

| * | 0 | 1 | 2 | 3 | 4 | 5 |
|---|---|---|---|---|---|---|
| 0 | 0 | 0 | 0 | 0 | 0 | 0 |
| 1 | 3 | 3 | 3 | 3 | 3 | 3 |
| 2 | 0 | 0 | 0 | 0 | 0 | 0 |
| 3 | 3 | 3 | 3 | 3 | 3 | 3 |
| 4 | 0 | 0 | 0 | 0 | 0 | 0 |
| 5 | 3 | 3 | 3 | 3 | 3 | 3 |

{3} is a semigroup so $Z_6$ (3, 0) is a SG.

**THEOREM 5.4.11:** *Let $Z_n$ (m, 0) belong to $Z^{***}(n)$ with $m * m = m \pmod{n}$. Then $Z_n$ (m, 0) is a SG, which is Smarandache strong Bol groupoid, Smarandache strong Moufang groupoid, Smarandache strong P-groupoid and Smarandache strong alternative groupoid.*

*Proof:* Since $m \in Z_n \setminus \{0\}$ is such that $m * m \equiv m \pmod{n}$ so $Z_n$ (m, 0) and $Z_n$ (0, m) are SGs as {m} is a proper subset which is a semigroup. It is left for the reader to verify $Z_n$ (m, 0) and $Z_n$ (0, m) are Smarandache strong Bol groupoid, Smarandache strong Moufang groupoid, Smarandache strong P-groupoid and Smarandache strong alternative groupoid.

**COROLLARY 5.4.12:** *Let $Z_n = \{0, 1, 2, \ldots, n-1\}$, $n > 3$, $n < \infty$. $Z_n(m, 0)$ is a SG for every idempotent $m \in Z_n.(m^2 \equiv m \pmod{n})$.*

*Proof*: If $m^2 \equiv m \pmod{n}$ we have just proved in Theorem 5.4.11 $Z_n(m, 0)$ and $Z_n(0, m)$ are SGs.

***Example 5.4.5:*** The SGs in $Z^{***}(12) \setminus Z^{**}(12)$ are $Z_{12}(4,0)$, $Z_{12}(0, 4)$, $Z_{12}(9, 0)$ and $Z_{12}(0, 9)$.

***Example 5.4.6:*** Consider the groupoids $Z_{15}(6, 0)$, $Z_{15}(0, 6)$, $Z_{15}(10, 0)$ and $Z_{15}(0, 10)$ in $Z^{***}(15) \setminus Z^{**}(15)$ they are SGs which are Smarandache Bol groupoids, Smarandache Moufang groupoids, Smarandache alternative groupoids and Smarandache P-groupoid.

This is left for the reader to verify the above facts.

**PROBLEM 1:** Find the number of SGs in $Z^{***}(15)$.

**PROBLEM 2:** Find the SG in $Z^{***}(8)$.

**PROBLEM 3:** Prove $Z_{33}(11,0)$ is a SG.

**PROBLEM 4:** Find all the SGs in $Z^{***}(32)$.

**PROBLEM 5:** Does $Z^{***}(31)$ have SGs?

(Note: A SG must be only in $Z^{***}(31) \setminus Z^{**}(31)$ and it is not $Z_{31}(1, 0)$ or $Z_{31}(0, 1)$).

**PROBLEM 6:** Find all SGs in $Z^{***}(18) \setminus Z^{**}(18)$.

**PROBLEM 7:** Is the groupoid $Z_7(3, 0)$ a SG? Justify.

**PROBLEM 8:** Find all the SGs in $Z^{***}(20) \setminus Z^{**}(20)$.

**PROBLEM 9:** Find the Smarandache Bol groupoid in $Z^{***}(20) \setminus Z^{**}(20)$.

**PROBLEM 10:** Find the Smarandache Moufang groupoid in $Z^{***}(26) \setminus Z^{**}(26)$.

## 5.5 Smarandache Direct Product Using the New Class of Smarandache Groupoids

In this section, we study the Smarandache direct product of groupoids using new classes of SGs. The main use of direct product is that in all the new classes of SGs built using $Z_n$ we saw the proper subsets which were semigroups just had order 2 or 1 elements, but we are able to overcome this problem once we make use of Smarandache direct products. We can get semigroups of desired order by using desired number of groupoids. To get both heterogeneous and mixed structures we can use groupoids from the different classes $Z(n)$, $Z^*(n)$, $Z^{**}(n)$ and $Z^{***}(n)$. Further, we can have proper Smarandache subgroupoids in them.

***Example 5.5.1:*** Let $Z_6(3, 4) \in Z(6)$ and $Z_6(5, 2) \in Z(6)$. The direct product $Z_6(3, 4) \times Z_6(5, 2) = \{(x, y) / x \in Z_6(3, 4) \text{ and } y \in Z_6(5, 2)\}$ has many subsets, which are semigroups.

***Example 5.5.2:*** $G = Z_7(5, 3) \times Z_6(3, 0)$ is a SG of order 42. This has atleast 2 proper Smarandache subgroupoids of order 6.

Thus we can study about Smarandache normal subgroupoids and so on. Thus, we have an interesting result.

**THEOREM 5.5.1:** *Let $G_1, G_2, \ldots, G_n$ be n SGs from the class of groupoids $Z^{***}(m)$ for varying m. Let $G = G_1 \times G_2 \times \ldots \times G_n$. Then G has atleast $n + nC_2 + \ldots + nC_{n-1}$ Smarandache subgroupoids.*

*Proof*: Given $G = G_1 \times G_2 \times \ldots \times G_n$ and each $G_i$ is a SG. Now $\overline{G}_i = \{0\} \times \{0\} \times \ldots \times G_i \times \{0\} \times \ldots \times \{0\}$. So for $i = 1, 2, \ldots, n$ we get n Smarandache subgroupoids of G. Take $\overline{G}_i \, \overline{G}_j = \{0\} \times \{0\} \times \ldots \times G_i \times \ldots \times G_j \times \ldots \times \{0\}$ with $i \neq j$ we have $nC_2$ distinct Smarandache subgroupoids of G. Similarly $\overline{G}_i \, \overline{G}_j \, \overline{G}_k = \{0\} \times \{0\} \times \ldots \times G_i \times \ldots \times G_j \times \ldots \times G_k \times \ldots \times \{0\}$, i, j and k are distinct. We get from this $nC_3$ Smarandache subgroupoids of G and so on. Thus $\overline{G}_1 \times \ldots \times \{0\} \times \ldots \times \overline{G}_n$ gives $nC_{n-1}$ that is n Smarandache subgroupoids of G. Hence we have in G at least $n + nC_1 + \ldots + nC_{n-1}$ number of Smarandache subgroupoids of G. Thus even if none of the groupoids $G_1, G_2, \ldots, G_n$ have proper Smarandache subgroupoids still we may at least get $n + nC_2 + \ldots + nC_{n-1}$ number of Smarandache subgroupoids in $G = G_1 \times \ldots \times G_n$.

Hence this theorem helps us to obtain just Smarandache Bol groupoids, Smarandache Moufang groupoids, Smarandache alternative groupoids and Smarandache P-groupoids which is got even if one of the groupoids in this set $G_1, G_2, \ldots, G_n$ happens to be even Smarandache strong Bol groupoid or Smarandache strong Moufang groupoid and so on. Thus, we see from the technique of Smarandache direct product of groupoids we are able to get varieties of SGs. Just as in the case of direct product of groups, we see even if one of the SGs in this class is a Smarandache normal groupoid, G has Smarandache normal subgroupoids. Thus, Smarandache direct products helps us in many ways to get a desired form of SGs. This concept of Smarandache direct product of groupoids will find its application in the Smarandache semi-automaton and Smarandache automaton respectively discussed in chapter 6.

***Example 5.5.3:*** Let $Z_3(2, 1)$ and $Z_4(1, 3)$ be two groupoids. $G = Z_3(2, 1) \times Z_4(1, 3)$. The groupoid G is even by the following table:

| * | (00) | (01) | (02) | (03) | (10) | (11) | (12) | (13) | (20) | (21) | (22) | (23) |
|---|---|---|---|---|---|---|---|---|---|---|---|---|
| (00) | (00) | (03) | (02) | (01) | (10) | (13) | (12) | (11) | (20) | (23) | (22) | (21) |
| (01) | (01) | (00) | (03) | (02) | (11) | (10) | (13) | (12) | (21) | (20) | (23) | (22) |
| (02) | (02) | (01) | (00) | (03) | (12) | (11) | (10) | (13) | (22) | (21) | (20) | (23) |
| (03) | (03) | (02) | (01) | (00) | (13) | (12) | (11) | (10) | (23) | (22) | (21) | (20) |
| (10) | (20) | (23) | (22) | (21) | (00) | (03) | (02) | (01) | (10) | (13) | (12) | (11) |
| (11) | (21) | (20) | (23) | (22) | (01) | (00) | (03) | (02) | (11) | (10) | (13) | (12) |
| (12) | (22) | (21) | (20) | (23) | (02) | (01) | (00) | (03) | (12) | (11) | (10) | (13) |

| (13) | (23) | (22) | (21) | (20) | (03) | (02) | (01) | (00) | (13) | (12) | (11) | (10) |
|------|------|------|------|------|------|------|------|------|------|------|------|------|
| (20) | (10) | (13) | (12) | (11) | (20) | (23) | (22) | (21) | (00) | (03) | (02) | (01) |
| (21) | (11) | (10) | (13) | (12) | (21) | (20) | (23) | (22) | (01) | (00) | (03) | (02) |
| (22) | (12) | (11) | (10) | (13) | (22) | (21) | (20) | (23) | (02) | (01) | (00) | (03) |
| (23) | (13) | (12) | (11) | (10) | (23) | (22) | (21) | (20) | (03) | (02) | (01) | (00) |

We denote the pairs (x, y) by (xy) in the above table for notational convenience.

Let $G = Z_3 (2, 1) \times Z_4 (1, 3)$. This groupoid has 12 elements given by the table above. This has proper subgroupoids $\{(0, 0), (0, 1), (0, 2), (0, 3)\} = \{(00), (01), (02), (03)\}$, which is a Smarandache subgroupoid as $\{(0, 0), (0, 2)\} = \{(00), (02)\}$ is a semigroup. This example is specially given because we see the groupoid $Z_3 (2, 1)$ is not a SG still we can using it with a SG get a SG having subgroupoids which are also Smarandache subgroupoids. Thus, direct product helps us to find many interesting results.

**PROBLEM 1:** Find all Smarandache subgroupoids in the direct product $Z_{12} (3, 8) \times Z_8 (3, 6)$.

**PROBLEM 2:** Find all the normal groupoids of $Z^{**}(3) \times Z^{***}(3)$.

**PROBLEM 3:** Can the groupoid $Z_5 (2, 2) \times Z_{12} (2, 2)$ be Smarandache alternative? Justify your answer.

**PROBLEM 4:** How many Smarandache Bol groupoids are in the direct product $G = Z^{***}(6) \times Z^{**}(5) \times Z_8 (12)$?

**PROBLEM 5:** Find the Smarandache Moufang groupoids in the direct product $G = Z^{*}(6) \times Z^{**}(12) \times Z^{**}(5)$.

**PROBLEM 6:** Find all Smarandache P-groupoids in the direct product $G = Z^{***}(7) \times Z^{***}(8)$.

**PROBLEM 7:** Find all SGs in the direct product $Z (7) \times \{Z^{***}(7) \setminus Z (7)\}$.

**PROBLEM 8:** Find all Smarandache idempotent groupoids in the direct product $Z^{***}(6) \times Z^{***}(7) \times Z^{***}(8)$.

**PROBLEM 9:** How many SGs are in the direct product $Z (8) \times Z (7)$?

**PROBLEM 10:** Find all SGs in $Z (11) \times Z (13)$.

## 5.6 Smarandache Groupoids with Identity Using $Z_n$

Now we studied four new classes of groupoids $Z (n)$, $Z^{*}(n)$, $Z^{**}(n)$ and $Z^{***}(n)$ using $Z_n$ but none of them which are not semigroups have identity in them, so we have to define once again a new class of groupoids using $Z_n \cup \{e\} = G$ where $e \notin Z_n$.

We define operation on G by

1. $a * a = e$.
2. $a * e = e * a = a$ for all $a \in G$.
3. $a \neq b$, $a * b = ta + ub \pmod{n}$ where $t, u \in Z_n \setminus \{0\}$, $(t, u) = 1$ and $t \neq u$.

We denote this class of groupoids by $G(n)$. If we make $t, u \in Z_n \setminus \{0\}$, $(t, u) = d$, with $t \neq u$ then we denote the class by $G^*(n)$. Similarly if we permit $t, u \in Z_n \setminus \{0\}$ such that $t = u$ can also occur then we denote this class by $G^{**}(n)$. If we have one of $t$ or $u$ to be zero we denote it by $G^{***}(n)$.

Thus, we have groupoids with identity as $G(n) \subset G^*(n) \subset G^{**}(n) \subset G^{***}(n)$. Now our aim is to study which of these groupoids are SGs, Smarandache P-groupoids and so on.

***Example 5.6.1:*** Let $G = (Z_4 \cup \{e\})$ define $*$ on $G$ by $a * b = 2a + 3b \pmod{4}$ for $a, b \in Z_4$ is given by the following table:

| * | e | 0 | 1 | 2 | 3 |
|---|---|---|---|---|---|
| e | e | 0 | 1 | 2 | 3 |
| 0 | 0 | e | 3 | 2 | 1 |
| 1 | 1 | 2 | e | 0 | 3 |
| 2 | 2 | 0 | 3 | e | 1 |
| 3 | 3 | 2 | 1 | 0 | e |

This is a SG with identity of order 5 as $\{e, m\}$ is a semigroup for $m = 1, 2, 3$.

**THEOREM 5.6.1:** *All groupoids in $G^{***}(n)$ are SGs.*

*Proof:* Clearly by the definition of $G^{***}(n)$ we have $\{e, m\}$ for $m = 1, 2, 3, \ldots, n - 1$ are semigroups. So every groupoid in $G^{***}(n)$ are SGs.

*Remark:* Since $G(n) \subset G^*(n) \subset G^{**}(n) \subset G^{***}(n)$ for every n we see all groupoids in $G(n)$, $G^*(n)$, $G^{**}(n)$ are SGs as we have proved in Theorem 5.6.1 that all groupoids in $G^{***}(n)$ is a SG.

**THEOREM 5.6.2:** *No groupoid in $G^{***}(n)$ is a Smarandache idempotent groupoid.*

*Proof:* Since by the very definition of G we see $x * x = e$, so no element in G is an idempotent element, so no groupoid in $G^{***}(n)$ is a Smarandache idempotent groupoid.

***Example 5.6.2:*** Let $G = Z_6 \cup \{e\}$ for $a, b \in Z_6$ define $a * b = 5a + 3b$. The groupoid is given by the following table:

| * | e | 0 | 1 | 2 | 3 | 4 | 5 |
|---|---|---|---|---|---|---|---|
| e | e | 0 | 1 | 2 | 3 | 4 | 5 |
| 0 | 0 | e | 3 | 0 | 3 | 0 | 3 |
| 1 | 1 | 5 | e | 5 | 2 | 5 | 2 |
| 2 | 2 | 4 | 1 | e | 1 | 4 | 1 |

| 3 | 3 | 3 | 0 | 3 | e | 3 | 0 |
|---|---|---|---|---|---|---|---|
| 4 | 4 | 2 | 5 | 2 | 5 | e | 5 |
| 5 | 5 | 1 | 4 | 1 | 4 | 1 | e |

This groupoid is Smarandache strong right alternative, but it is not even Smarandache left alternative. Now we give an example of a left alternative groupoid in G (6), which is not right alternative.

*Example 5.6.3:* The groupoid in G (6) is given by the following table:

| * | e | 0 | 1 | 2 | 3 | 4 | 5 |
|---|---|---|---|---|---|---|---|
| e | e | 0 | 1 | 2 | 3 | 4 | 5 |
| 0 | 0 | e | 5 | 4 | 3 | 2 | 1 |
| 1 | 1 | 4 | e | 2 | 1 | 0 | 5 |
| 2 | 2 | 2 | 1 | e | 5 | 4 | 3 |
| 3 | 3 | 0 | 5 | 3 | e | 2 | 1 |
| 4 | 4 | 4 | 3 | 2 | 1 | e | 5 |
| 5 | 5 | 2 | 1 | 0 | 5 | 4 | e |

This is a Smarandache strong left alternative groupoid for $(x * x) * y = x * (x * y)$. Consider $(x * x) * y = y$. Now $x * (x * y) = x * [4x + 5y] = 4x + 20x + 25y = y$. Hence the claim.

This is not Smarandache right alternative, it is left for the reader to verify. In view of this, we have the following theorem.

**THEOREM 5.6.3:** *Groupoids in $G^{***}(n)$ are Smarandache strong right alternative groupoid if and only if $t^2 \equiv 1 \pmod{n}$ and $tu + u \equiv 0 \pmod{n}$.*

*Proof:* A SG G in $G^{***}(n)$ is Smarandache strong right alternative if for x, y ∈ G we must have $(x * y) * y = x * (y * y)$. $(x * y) * y = (tx + uy) * y = t^2x + tuy + uy \pmod{n}$. Now $x * (y * y) = x$. $(x * y) * y \equiv x * (y * y) \pmod{n}$ if and only if $t^2 \equiv 1 \pmod{n}$ and $tu + u \equiv 0 \pmod{n}$. Similarly we get the following characterization theorem for a groupoid in $G^{***}(n)$ to be Smarandache strong left alternative.

**THEOREM 5.6.4:** *Let G be a groupoid in $G^{***}(n)$. G is a Smarandache strong left alternative if and only if $u^2 \equiv 1 \pmod{n}$ and $(t + tu) \equiv 0 \pmod{n}$.*

*Proof:* On similar lines, the proof can be given as in Theorem 5.6.3.

**THEOREM 5.6.5:** *Let G (n) be a groupoid in $G^{***}(n)$, n not a prime. G (n) is a Smarandache strong Bol groupoid, Smarandache strong Moufang groupoid, Smarandache strong P-groupoid only when $t^2 \equiv t \pmod{n}$ and $u^2 \equiv u \pmod{n}$.*

*Proof:* The verification is left for the reader.

**PROBLEM 1:** Find all SGs with identity, which are Smarandache strong Moufang groupoids in $G^{***}(7) \setminus G^{**}(7)$.

**PROBLEM 2:** Find Smarandache P-groupoids in $G^{**}(8) \setminus G(8)$.

**PROBLEM 3:** Find Smarandache Bol groupoids in $G^{***}(9) \setminus G^{*}(9)$.

**PROBLEM 4:** Find the number of SGs in $G^{*}(12)$.

**PROBLEM 5:** Find the number of Smarandache Bol groupoids in $G^{*}(11) \setminus G(11)$.

**PROBLEM 6:** Find the number of Smarandache alternative groupoids in $G^{*}(5) \setminus G(5)$.

**PROBLEM 7:** How many Smarandache strong right alternative groupoids exist in $G^{***}(12) \setminus G^{**}(12)$?

## Supplementary Reading

1. G. Birkhoff and S. S. Maclane, *A Brief Survey of Modern Algebra*, New York, The Macmillan and Co. (1965).

2. R. H. Bruck, *A Survey of Binary Systems*, Springer Verlag, (1958).

3. W.B.Vasantha Kandasamy, *New classes of finite groupoids using $Z_n$*, Varahmihir Journal of Mathematical Sciences, Vol 1, 135-143, (2001).

4. W. B. Vasantha Kandasamy, *Smarandache Groupoids*, http://www.gallup.unm.edu/~smarandache/Groupoids.pdf.

CHAPTER SIX
# SMARANDACHE SEMI AUTOMATON AND SMARANDACHE AUTOMATON

The semi automaton and automaton are built using the fundamental algebraic structure semigroups. In this chapter we use the generalized concept of semigroups viz. Smarandache groupoids which always contains a semigroup in them are used to construct Smarandache semi automaton and Smarandache automaton. Thus, the Smarandache groupoids find its application in the construction of finite machines. Here we introduce the concept of Smarandache semi automaton and Smarandache automaton using Smarandache free groupoids. This chapter starts with the definition of Smarandache free groupoids.

## 6.1 Basic Results

**DEFINITION:** *Let S be non empty set. Generate a free groupoid using S and denote it by <S>. Clearly the free semigroup generated by the set S is properly contained in <S>; as in <S> we may or may not have the associative law to be true.*

*Remark:* Even $(ab) c \neq a (bc)$ in general for all a, b, c $\in$ <S>. Thus unlike a free semigroup where the operation is associative, in case of free groupoid we do not assume the associativity while placing them in justra position.

**THEOREM 6.1.1:** *Every free groupoid is a Smarandache free groupoid.*

*Proof:* Clearly if S is the set which generates the free groupoid then it will certainly contain the free semigroup generated by S, so every free groupoid is a Smarandache free groupoid.

We just recall the definition of semi automaton and automaton from the book of R.Lidl and G.Pilz.

**DEFINITION:** *A Semi Automaton is a triple Y = (Z, A, $\delta$) consisting of two non - empty sets Z and A and a function $\delta$: Z $\times$ A $\to$ Z, Z is called the set of states, A the input alphabet and $\delta$ the "next state function" of Y.*

Let $A = \{a_1, ..., a_n\}$ and $Z = \{z_1, ..., z_k\}$. The semi automaton $Y = (Z, A, \delta)$. The semi automaton can also be described by the transition table.

Description by Table:

| $\delta$ | $a_1$ | $\cdots$ | $a_n$ |
|---|---|---|---|
| $z_1$ | $\delta(z_1, a_1)$ | $\cdots$ | $\delta(z_1, a_n)$ |
| $\vdots$ | $\vdots$ | | $\vdots$ |
| $z_k$ | $\delta(z_k, a_1)$ | $\cdots$ | $\delta(z_k, a_n)$ |

as $\delta : Z \times A \to Z$, $\delta(z_i, a_j) \in Z$. The semi automaton can also be described by graphs.

Description by graphs:

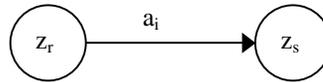

We depict $z_1, ..., z_k$ as 'discs' in the plane and draw an arrow labeled $a_i$ from $z_r$ to $z_s$ if $\delta(z_r, a_i) = z_r$. This graph is called the state graph.

***Example 6.1.1:*** $Z$ – set of states and $A$ – input alphabet. Let $Z = \{0, 1, 2\}$ and $A = \{0, 1\}$. The function $\delta: Z \times A \to Z$ defined by $\delta(0, 1) = 1 = \delta(2, 1) = \delta(1, 1)$, $\delta(0, 0) = 0$, $\delta(2, 0) = 1$, $\delta(1, 0) = 0$.

This is a semi automaton, having the following graph.

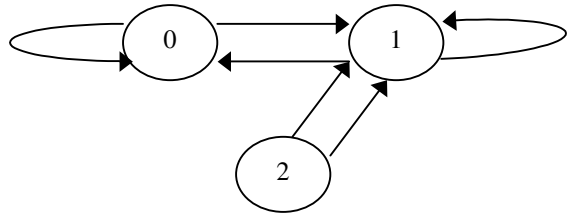

The description by table is:

| $\delta$ | 0 | 1 |
|---|---|---|
| 0 | 0 | 1 |
| 1 | 0 | 1 |
| 2 | 1 | 1 |

**DEFINITION:** *An Automaton is a quintuple $K = (Z, A, B, \delta, \Lambda)$ where $(Z, A, \delta)$ is a semi automaton, B is a non-empty set called the output alphabet and $\lambda: Z \times A \to B$ is the output function.*

If $z \in Z$ and $a \in A$, then we interpret $\delta(z, a) \in Z$ as the next state into which z is transformed by the input a, $\lambda(z, a) \in B$ is the output resulting from the input a.

Thus if the automaton is in state z and receives input a, then it changes to state $\delta(z, a)$ with output $\lambda(z, a)$.

If $A = \{a_1, ..., a_n\}$, $B = \{b_1, b_2, ..., b_m\}$ and $Z = \{z_1, ..., z_k\}$, $\delta : Z \times A \to Z$ and $\lambda : Z \times A \to B$ given by the description by tables where $\delta(z_k, a_i) \in Z$ and $\lambda(z_k, a_j) \in B$.

Description by Tables:

| $\delta$ | $a_1$ | ... | $a_n$ |
|---|---|---|---|
| $z_1$ | $\delta(z_1, a_1)$ | ... | $\delta(z_1, a_n)$ |
| $\vdots$ | $\vdots$ | | $\vdots$ |
| $z_k$ | $\delta(z_k, a_1)$ | ... | $\delta(z_k, a_n)$ |

In case of automatan, we also need an output table.

| $\lambda$ | $a_1$ | ... | $a_n$ |
|---|---|---|---|
| $z_1$ | $\lambda(z_1, a_1)$ | ... | $\lambda(z_1, a_n)$ |
| $\vdots$ | $\vdots$ | | $\vdots$ |
| $z_k$ | $\lambda(z_k, a_1)$ | ... | $\lambda(z_k, a_n)$ |

Description by Graph:

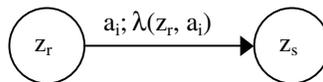

*Example 6.1.2:* Let $Z = \{z_1, z_2\}$ $A = B = \{0, 1\}$ we have the following table and graph:

| $\delta$ | 0 | 1 |
|---|---|---|
| $z_0$ | $z_0$ | $z_1$ |
| $z_1$ | $z_1$ | $z_0$ |

| $\lambda$ | 0 | 1 |
|---|---|---|
| $z_0$ | 0 | 1 |
| $z_1$ | 0 | 1 |

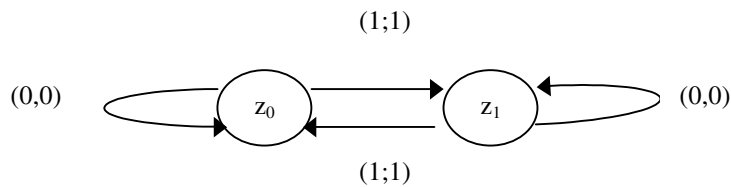

This automaton is know as the parity check automaton.

Now it is important and interesting to note that Z, A and B are only non-empty sets. They have no algebraic operation defined on them. The automatons and semi automatons defined in this manner do not help one to perform sequential operations. Thus, it is reasonable to consider the set of all finite sequences of elements of the set A including the empty sequence $\Lambda$.

In other words, in our study of automaton we extend the input set A to the free monoid $\overline{A}$ and $\overline{\lambda}: Z \times \overline{A} \to \overline{B}$ where $\overline{B}$ is the free monoid generated by B. We also extend functions $\delta$ and $\lambda$ from $Z \times A$ to $Z \times \overline{A}$ by defining $z \in Z$ and $a_1, ..., a_n \in A$ by $\overline{\delta}: Z \times \overline{A} \to Z$.

$$\overline{\delta}(z, \Lambda) = z$$
$$\overline{\delta}(z, a_1) = \delta(z, a_1)$$
$$\overline{\delta}(z, a_1 a_2) = \delta(\overline{\delta}(z, a_1), a_2)$$
$$\vdots$$
$$\overline{\delta}(z, a_1 a_2 ... a_n) = \delta(\overline{\delta}(z, a_1 a_2 ... a_{n-1}), a_n)$$

and
$\lambda: Z \times A \to B$ by $\overline{\lambda}: Z \times \overline{A} \to \overline{B}$
by
$$\overline{\lambda}(z, \Lambda) = \Lambda$$
$$\overline{\lambda}(z, a_1) = \lambda(z, a_1)$$
$$\overline{\lambda}(z, a_1 a_2) = \lambda(z, a_1) \overline{\lambda}(\delta(z, a_1), a_2)$$
$$\vdots$$
$$\overline{\lambda}(z, a_1 a_2 ... a_n) = \lambda(z, a_1) \overline{\lambda}(\delta(z, a_1), a_2 ... a_r)$$

The semi automaton $Y = (Z, A, \delta)$ and automaton $K = (Z, A, B, \delta, \lambda)$ is thus generalized to the new semi automaton $Y = (Z, \overline{A}, \overline{\delta})$ and automaton $K = (Z, \overline{A}, \overline{B}, \overline{\delta}, \overline{\lambda})$.

We throughout in this book mean by new semi automaton $Y = (Z, \overline{A}, \overline{\delta})$ and new automaton $K = (Z, \overline{A}, \overline{B}, \overline{\delta}, \overline{\lambda})$.

## 6.2 Smarandache Semi Automaton and Smarandache Automaton

The concept of Smarandache semi automaton and Smarandache automaton was first defined by the author in the year 2002. The study of Smarandache semi automaton and automata use the concept of groupoids and free groupoids.

Now we define Smarandache semi automaton and Smarandache automaton as follows:

**DEFINITION:** $Y_s = (Z, \overline{A}_s, \overline{\delta}_s)$ is said to be a Smarandache semi automaton if $\overline{A}_s = \langle A \rangle$ is the free groupoid generated by A with $\Lambda$ the empty element adjoined with it and $\overline{\delta}_s$ is the function from $Z \times \overline{A}_s \to Z$. Thus the Smarandache semi automaton contains $Y = (Z, \overline{A}, \overline{\delta})$ as a new semi automaton which is a proper sub-structure of $Y_s$.

*Or equivalently, we define a Smarandache semi automaton as one, which has a new semi automaton as a sub-structure.*

The advantages of the Smarandache semi automaton are if for the triple $Y = (Z, A, \delta)$ is a semi automaton with Z, the set of states, A the input alphabet and $\delta : Z \times A \to Z$ is the next state function.

When we generate the Smarandache free groupoid by A and adjoin with it the empty alphabet $\Lambda$ then we are sure that $\overline{A}$ has all free semigroups. Thus, each free semigroup will give a new semi automaton. Thus by choosing a suitable A we can get several new semi automaton using a single Smarandache semi automaton.

We now give some examples of Smarandache semi automaton using finite groupoids. When examples of semi automaton are given usually the books use either the set of modulo integers $Z_n$ under addition or multiplication we use groupoid built using $Z_n$.

**DEFINITION:** $\overline{Y}_s^{'} = (Z_1, \overline{A}_s, \overline{\delta}_s^{'})$ is called the Smarandache sub semi automaton of $\overline{Y}_s = (Z_2, \overline{A}_s, \overline{\delta}_s^{'})$ denoted by $\overline{Y}_s^{'} \leq \overline{Y}_s$ if $Z_1 \subset Z_2$ and $\overline{\delta}_s^{'}$ is the restriction of $\overline{\delta}_s$ on $Z_1 \times \overline{A}_s$ and $\overline{Y}_s^{'}$ has a proper subset $\overline{H} \subset \overline{Y}_s^{'}$ such that $\overline{H}$ is a new semi automaton.

*Example 6.2.1:* Let $Z = Z_4(2, 1)$ and $A = Z_6(2, 1)$. The Smarandache semi automaton (Z, A, $\delta$) where $\delta : Z \times A \to Z$ is given by $\delta(z, a) = z \bullet a \pmod 4$. We get the following table:

| δ | 0 | 1 | 2 | 3 |
|---|---|---|---|---|
| 0 | 0 | 1 | 2 | 3 |
| 1 | 2 | 3 | 0 | 1 |
| 2 | 0 | 1 | 2 | 3 |
| 3 | 2 | 3 | 0 | 1 |

We get the following graph for this Smarandache semi automaton.

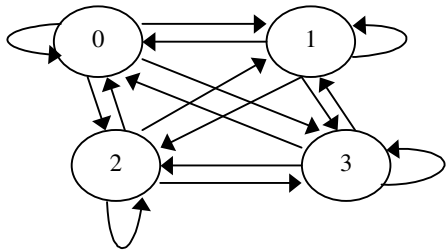

This has a nice Smarandache sub semi automaton given by the table.

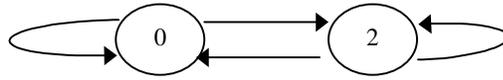

This Smarandache semi automaton is nothing but given by the states {0, 2}.

***Example 6.2.2:*** Let $Z = Z_3(1, 2)$ and $A = Z_4(2, 2)$ the triple $(Z, A, \delta)$ is a Smarandache semi automaton with $\delta(z, a) = (z * a) \pmod 3$ '$*$' the operation in A given by the following table:

| $\delta$ | 0 | 1 | 2 | 3 |
|---|---|---|---|---|
| 0 | 0 | 2 | 1 | 0 |
| 1 | 2 | 1 | 0 | 2 |
| 2 | 1 | 0 | 2 | 1 |

The graphical representation of the Smarandache semi automaton is

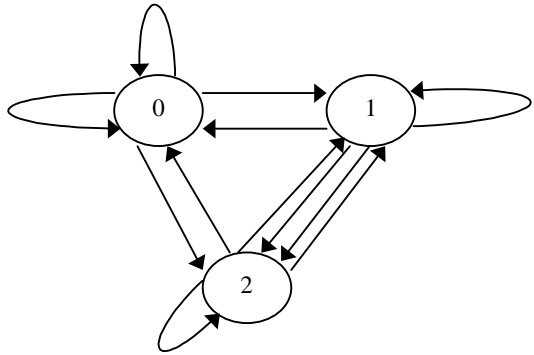

***Example 6.2.3:*** Now let $Z = Z_4(2, 2)$ and $A = Z_3(1, 2)$ we define $\delta(z, a) = (z * a) \pmod 4$ '$*$' as in Z. The table for the semi automaton is given by

| $\delta$ | 0 | 1 | 2 |
|---|---|---|---|
| 0 | 0 | 2 | 0 |
| 1 | 2 | 0 | 2 |
| 2 | 0 | 2 | 0 |
| 3 | 2 | 0 | 2 |

The graph for it is

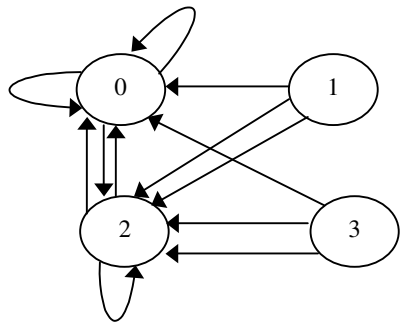

Thus this has a Smarandache sub semi automaton $Z_1$ given by $Z_1 = \{0, 2\}$ states.

**DEFINITION:** $\overline{K}_s = (Z, \overline{A}_s, \overline{B}_s, \overline{\delta}_s, \overline{\lambda}_s)$ *is defined to be a Smarandache automaton if* $\overline{K} = (Z, \overline{A}_s, \overline{B}_s, \overline{\delta}_s, \overline{\lambda}_s)$ *is the new automaton and* $\overline{A}_s$ *and* $\overline{B}_s$, *the Smarandache free groupoids so that* $\overline{K} = (Z, \overline{A}_s, \overline{B}_s, \overline{\delta}_s, \overline{\lambda}_s)$ *is the new automaton got from K and* $\overline{K}$ *is strictly contained in* $\overline{K}_s$.

Thus Smarandache automaton enables us to adjoin some more elements which is present in A and freely generated by A, as a free groupoid; that will be the case when the compositions may not be associative.

Secondly, by using Smarandache automaton we can couple several automaton as

$$Z = Z_1 \cup Z_2 \cup ... \cup Z_n$$
$$A = A_1 \cup A_2 \cup ... \cup A_n$$
$$B = B_1 \cup B_2 \cup ... \cup B_n$$
$$\lambda = \lambda_1 \cup \lambda_2 \cup ... \cup \lambda_n$$
$$\delta = \delta_1 \cup \delta_2 \cup ... \cup \delta_n.$$

where the union of $\lambda_i \cup \lambda_j$ and $\delta_i \cup \delta_j$ denote only extension maps as '$\cup$' has no meaning in the composition of maps, where $K_i = (Z_i, A_i, B_i, \delta_i, \lambda_i)$ for $i = 1, 2, 3, ..., n$ and $\overline{K} = \overline{K}_1 \cup \overline{K}_2 \cup ... \cup \overline{K}_n$. Now $\overline{K}_s = (\overline{Z}_s, \overline{A}_s, \overline{B}_s, \overline{\lambda}_s, \overline{\delta}_s)$ is the Smarandache Automaton.

A machine equipped with this Smarandache Automaton can use any new automaton as per need.

We give some examples of Smarandache automaton using Smarandache groupoids.

***Example 6.2.4:*** Let $Z = Z_4 (3, 2)$, $A = B = Z_5 (2, 3)$. $K = (Z, A, B, \delta, \lambda)$ is a Smarandache automaton defined by the following tables where $\delta(z, a) = z * a \pmod 4$ and $\lambda(z, a) = z * a \pmod 5$.

| λ | 0 | 1 | 2 | 3 | 4 |
|---|---|---|---|---|---|
| 0 | 0 | 3 | 1 | 4 | 2 |
| 1 | 2 | 0 | 3 | 1 | 4 |
| 2 | 4 | 2 | 0 | 3 | 1 |
| 3 | 1 | 4 | 2 | 0 | 3 |

| δ | 0 | 1 | 2 | 3 | 4 |
|---|---|---|---|---|---|
| 0 | 0 | 2 | 0 | 2 | 0 |
| 1 | 3 | 1 | 3 | 1 | 3 |
| 2 | 2 | 0 | 2 | 0 | 2 |
| 3 | 1 | 3 | 1 | 3 | 1 |

We obtain the following graph:

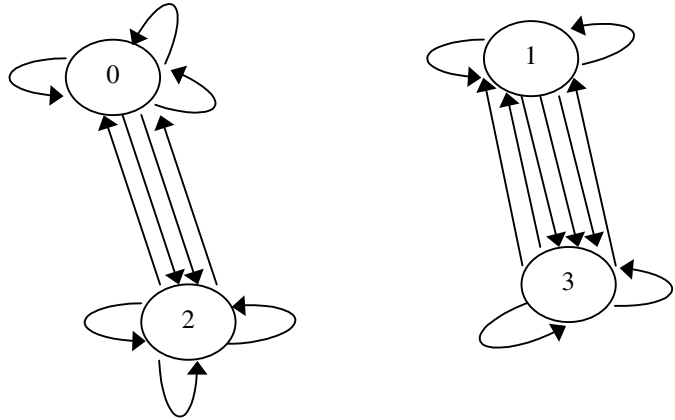

Thus we see this automaton has 2 Smarandache sub automatons given by the sates {0, 2} and {1, 3}.

***Example 6.2.5:*** Consider $Z_5 = Z_5(3, 2)$, $A = Z_3(0, 2)$, $B = Z_4(2, 3)$, $\delta(z, a) = z * a \pmod 5$, $\lambda(z, a) = z * a \pmod 4$. The Smarandache automaton $(Z, A, B, \delta, \lambda)$ is given by the following tables:

| δ | 0 | 1 | 2 |
|---|---|---|---|
| 0 | 0 | 3 | 1 |
| 1 | 3 | 1 | 4 |
| 2 | 1 | 4 | 2 |
| 3 | 4 | 2 | 0 |
| 4 | 2 | 0 | 3 |

| λ | 0 | 1 | 2 | 3 | 4 |
|---|---|---|---|---|---|
| 0 | 0 | 3 | 2 | 1 | 0 |
| 1 | 2 | 1 | 0 | 3 | 2 |
| 2 | 0 | 3 | 2 | 1 | 0 |
| 3 | 2 | 1 | 0 | 3 | 2 |
| 4 | 0 | 3 | 2 | 1 | 0 |

This has no proper Smarandache sub automaton.

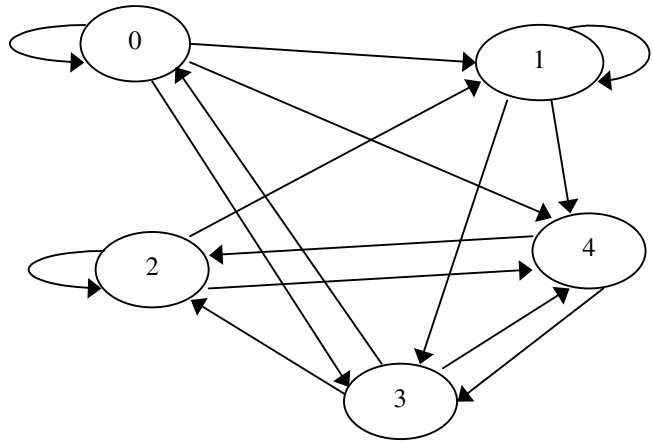

**DEFINITION:** $\overline{K}_s' = (Z_1, \overline{A}_s, \overline{B}_s, \overline{\delta}_s, \overline{\lambda}_s)$ is called Smarandache sub-automaton of $\overline{K}_s = (Z_2, \overline{A}_s, \overline{B}_s, \overline{\delta}_s, \overline{\lambda}_s)$ denoted by $\overline{K}_s' \leq \overline{K}_s$ if $Z_1 \subseteq Z_2$ and $\overline{\delta}_s'$ and $\overline{\lambda}_s'$ are the restriction of $\overline{\delta}_s$ and $\overline{\lambda}_s$ respectively on $Z_1 \times \overline{A}_s$ and has a proper subset $\overline{H} \subset \overline{K}_s'$ such that $\overline{H}$ is a new automaton.

**DEFINITION:** Let $\overline{K}_1$ and $\overline{K}_2$ be any two Smarandache automaton where $\overline{K}_1 = (Z_1, \overline{A}_s, \overline{B}_s, \overline{\delta}_s, \overline{\lambda}_s)$ and $\overline{K}_2 = (Z_2, \overline{A}_s, \overline{B}_s, \overline{\delta}_s, \overline{\lambda}_s)$. A map $\phi : \overline{K}_1$ to $\overline{K}_2$ is a Smarandache automaton homomorphism if $\phi$ restricted from $K_1 = (Z_1, A_1, B_1, \delta_1, \lambda_1)$ and $K_2 = (Z_2, A_2, B_2, \delta_2, \lambda_2)$ denoted by $\phi_r$ is a automaton homomorphism from $K_1$ to $K_2$. $\phi$ is called a monomorphism (epimorphism or isomorphism) if there is an isomorphism $\phi_r$ from $K_1$ to $K_2$.

**DEFINITION:** Let $\overline{K}_1$ and $\overline{K}_2$ be two Smarandache automatons, where $\overline{K}_1 = (Z_1, \overline{A}_s, \overline{B}_s, \overline{\delta}_s, \overline{\lambda}_s)$ and $\overline{K}_2 = (Z_2, \overline{A}_s, \overline{B}_s, \overline{\delta}_s, \overline{\lambda}_s)$.

The Smarandache automaton direct product of $\overline{K}_1$ and $\overline{K}_2$ denoted by $\overline{K}_1 \times \overline{K}_2$ is defined as the direct product of the automaton $K_1 = (Z_1, A_1, B_1, \delta_1, \lambda_1)$ and $K_2 = (Z_2, A_2, B_2, \delta_2, \lambda_2)$ where $K_1 \times K_2 = (Z_1 \times Z_2, A_1 \times A_2, B_1 \times B_2, \delta, \lambda)$ with $\delta((z_1, z_2), (a_1, a_2)) = (\delta_1(z_1, a_1), \delta_2(z_2, a_2))$, $\lambda((z_1, z_2), (a_1, a_2)) = (\lambda_1(z_1, a_2), \lambda_2(z_2, a_2))$ for all $(z_1, z_2) \in Z_1 \times Z_2$ and $(a_1, a_2) \in A_1 \times A_2$.

*Remark:* Here in $\overline{K}_1 \times \overline{K}_2$ we do not take the free groupoid to be generated by $A_1 \times A_2$ but only free groupoid generated by $\overline{A}_1 \times \overline{A}_2$

Thus the Smarandache automaton direct product exists wherever a automaton direct product exists.

We have made this in order to make the Smarandache parallel composition and Smarandache series composition of automaton extendable in a simple way.

**DEFINITION:** *A Smarandache groupoid $G_1$ divides a Smarandache groupoid $G_2$ if the corresponding semigroups $S_1$ and $S_2$ of $G_1$ and $G_2$ respectively divides, that is, if $S_1$ is a homomorphic image of a sub-semigroup of $S_2$.*

In symbols $G_1 | G_2$. The relation divides is denoted by '|'.

**DEFINITION:** *Let $\overline{K}_1 = (Z_1, \overline{A}_s, \overline{B}_s, \overline{\delta}_s, \overline{\lambda}_s)$ and $\overline{K}_2 = (Z_2, \overline{A}_s, \overline{B}_s, \overline{\delta}_s, \overline{\lambda}_s)$ be two Smarandache automaton. We say the Smarandache automaton $\overline{K}_1$ divides the Smarandache automaton $\overline{K}_2$ if in the automatons $K_1 = (Z_1, A, B, \delta_1, \lambda_1)$ and $K_2 = (Z_2, A, B, \delta_2, \lambda_2)$, if $K_1$ is the homomorphic image of a sub-automaton of $K_2$. Notationally $K_1 | K_2$.*

**DEFINITION:** *Two Smarandache Automaton $\overline{K}_1$ and $\overline{K}_2$ are said to be equivalent if they divide each other. In symbols $\overline{K}_1 \sim \overline{K}_2$.*

## 6.3 Direct Product of Smarandache Automaton

We proceed on to define direct product Smarandache automaton and study about them. We can extend the direct product of semi automaton to more than two Smarandache automatons.

Using the definition of direct product of two automaton $K_1$ and $K_2$ with an additional assumption we define Smarandache series composition of automaton.

**DEFINITION:** *Let $K_1$ and $K_2$ be any two Smarandache automatons where $\overline{K}_1 = (Z_1, \overline{A}_s, \overline{B}_s, \overline{\delta}_s, \overline{\lambda}_s)$ and $\overline{K}_2 = (Z_2, \overline{A}_s, \overline{B}_s, \overline{\delta}_s, \overline{\lambda}_s)$ with an additional assumption $A_2 = B_1$.*

*The Smarandache automaton composition series denoted by $\overline{K}_1 \#\!\!\!+ \overline{K}_2$ of $\overline{K}_1$ and $\overline{K}_2$ is defined as the series composition of the automaton $K_1 = (Z_1, A_1, B_1, \delta_1, \lambda_1)$ and $K_2 = (Z_2, A_2, B_2, \delta_2, \lambda_2)$ with $\overline{K}_1 \#\!\!\!+ \overline{K}_2 = (Z_1 \times Z_2, A_1, B_2, \delta, \lambda)$ where $\delta((z_1, z_2), a_1) = (\delta_1(z_1, a_1), \delta_2(z_2, \lambda_1(z_1, a_1))$ and $\lambda((z_1, z_2), a_1) = (\lambda_2(z_2, \lambda_1(z_1, a_1))$ $((z_1, z_2) \in Z_1 \times Z_2, a_1 \in A_1)$.*

*This automaton operates as follows: An input $a_1 \in A_1$ operates on $z_1$ and gives a state transition into $z_1' = \delta_1(z_1, a_1)$ and an output $b_1 = \lambda_1(z_1, a_2) \in B_1 = A_2$. This output $b_1$ operates on $Z_2$ transforms a $z_2 \in Z_2$ into $z_2' = \delta_2(a_2, b_1)$ and produces the output $\lambda_2(z_2, b_1)$.*

*Then $\overline{K}_1 \#\!\!\!+ \overline{K}_2$ is in the next state $(z_1', z_2')$ which is clear from the following circuit:*

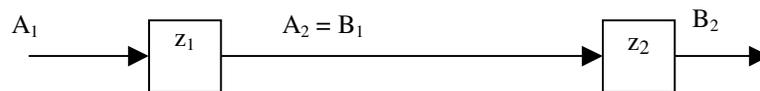

Now a natural question would be do we have a direct product, which corresponds to parallel composition, of the 2 Smarandache automatons $\overline{K}_1$ and $\overline{K}_2$.

Clearly the Smarandache direct product of automatons $\overline{K}_1 \times \overline{K}_2$ since $Z_1$, and $Z_2$ can be interpreted as two parallel blocks. $A_1$ operates on $Z_i$ with output $B_i$ ($i \in \{1,2\}$), $A_1 \times A_2$ operates on $Z_1 \times Z_2$, the outputs are in $B_1 \times B_2$.

The circuit is given by the following diagram:

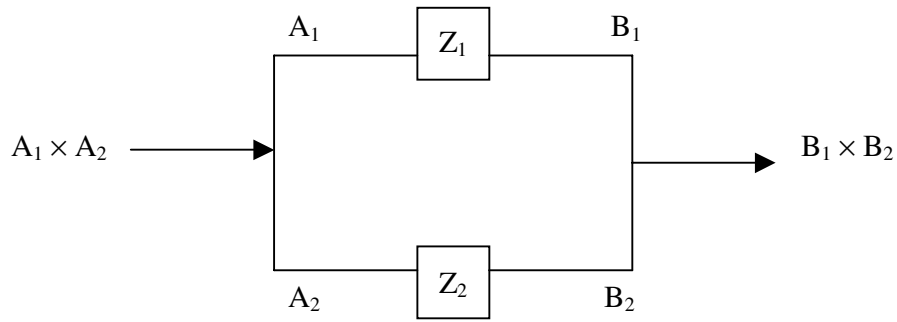

**PROBLEM 1:** Define Smarandache minimal automaton and give an example.

**PROBLEM 2:** Let $K_2 = \{(z_1, z_2, z_3)\}, \{a_1, a_2\}, \{0, 1\}, \delta, \lambda\}$

| $\delta$ | $a_1$ | $a_2$ |
|---|---|---|
| $z_1$ | $z_1$ | $z_3$ |
| $z_2$ | $z_2$ | $z_3$ |
| $z_3$ | $z_1$ | $z_2$ |
| $z_4$ | $z_1$ | $z_2$ |

| $\lambda$ | $a_1$ | $a_2$ |
|---|---|---|
| $z_1$ | 0 | 0 |
| $z_2$ | 0 | 0 |
| $z_3$ | 0 | 0 |
| $z_4$ | 1 | 0 |

Define the equivalence classes on Z as $\sim_1$ and $\sim_2$. Find $Z/\sim_1$ and $Z/\sim_2$.

**PROBLEM 3:** Give two examples of a Smarandache automaton using $A = \{a_1, a_2\}$ and $B = \{b_1, b_2, b_3\}$ with $Z = \{z_1, z_2, z_3, z_4, z_5\}$ with varying $\lambda$ and $\delta$.

**PROBLEM 4:** $A = B = Z = Z_4(2, 3)$. $\delta(z, a) = z * a \pmod 4$. $\lambda(z, a) = a * z \pmod 4$ Draw the graph of the Smarandache automaton.

**PROBLEM 5:** $A = B = Z = Z_4(3, 2)$ with $\delta$ and $\lambda$ defined as in problem 4. Draw the graph and compare the Smarandache automatons in Problem 4 and Problem 5.

**PROBLEM 6:** $A = B = Z_7(3, 2)$ is a Smarandache groupoid of order 7 and $Z = Z_4(2, 2)$, Smarandache groupoid of order 4. $\delta(z, a) = z * a \pmod 4$. Define $\lambda(z, a) = z * a \pmod 7$. Draw its graph.

**PROBLEM 7:** $A = Z_6(3, 2)$, $B = Z_5(3, 2)$ and $Z = Z_7(3, 2)$, $\delta(z, a) = z * a \pmod 7$. $\lambda(z, a) = z * a \pmod 5$. Draw its graph.

## Supplementary Reading

1. R. H. Bruck, *A Survey of Binary Systems*, Springer Verlag, (1958).

# CHAPTER SEVEN
# RESEARCH PROBLEMS

The Smarandache algebraic structure is a very new subject and there are a number of unsolved problems in this new subject. Here we concentrate only on problems connected with SGs. A list of problems is given to attract both algebraist students and researchers in the SGs. Some of the problems are difficult to solve.

**PROBLEM 1:** For a fixed integer n, n a prime; obtain the number of SGs in the class Z (n).

**PROBLEM 2:** How many SGs are there in the class Z (n) for a fixed n?

*Hint:* We know that the number of groupoids in the class Z (n) is bounded by $(n-1)(n-2)$, that is $|Z(n)| < (n-1)(n-2)$.

**PROBLEM 3:** If p is a prime. Prove $Z_p$ (m, m) (m < p) has no proper subgroupoid.

**PROBLEM 4:** Find conditions on $Z_n$ (m, m), m and n so that $Z_n$ (m, m) has normal subgroupoids.

*Hint:* If m is a prime, n any composite number. Can we say $Z_n$ (m, m) have no subgroupoids. For we see by this example.

*Example*: $Z_8$ (3, 3) given by the following table:

| * | 0 | 1 | 2 | 3 | 4 | 5 | 6 | 7 |
|---|---|---|---|---|---|---|---|---|
| 0 | 0 | 3 | 6 | 1 | 4 | 7 | 2 | 5 |
| 1 | 3 | 6 | 1 | 4 | 7 | 2 | 5 | 0 |
| 2 | 6 | 1 | 4 | 7 | 2 | 5 | 0 | 3 |
| 3 | 1 | 4 | 7 | 2 | 5 | 6 | 0 | 3 |
| 4 | 4 | 7 | 2 | 5 | 6 | 0 | 3 | 1 |
| 5 | 7 | 2 | 5 | 6 | 0 | 3 | 1 | 4 |
| 6 | 2 | 5 | 6 | 0 | 3 | 1 | 4 | 7 |
| 7 | 5 | 6 | 0 | 3 | 1 | 4 | 7 | 2 |

This is a commutative and non associative and has no subgroupoid.

**PROBLEM 5:** Suppose $Z_p(t, u)$ be a groupoid in $Z(p)$ for $0 < t + u < p$ prove $Z_p(t, u)$ is a SG; if p a prime.

**PROBLEM 6:** Obtain some interesting results on n, t, u in the groupoid $Z_n(t, u) \in Z(n)$ so that

1. $Z_n(t, u)$ is not a SG
2. $Z_n(t, u)$ is a SG.

*Note:* Condition for $Z_n(t, u) \in Z(n)$ to be a SG is given in this text but we want some other interesting conditions for $Z_n(t, u) \in Z(n)$ to be a SG.

**PROBLEM 7:** Let $Z_{p^m}(t, u)$ be a groupoid in $Z(p^m)$ (p a prime $m \geq 1$) with $(p, t) = 1$ and $(p, u) = 1$ and $t + u \equiv p^m - 1 \pmod{p^m}$. Can $Z_{p^m}(t, u)$ be a SG?

**PROBLEM 8:** Find SGs in $Z(n)$ other than the groupoids mentioned in this book.

**PROBLEM 9:** Can $Z_n(1, u) \in Z(n)$ be a SG when $(u, n) = 1$?

**PROBLEM 10:** Find the number of SGs in $Z(n)$.

**PROBLEM 11:** Find all Smarandache strong Bol groupoids in $Z(n)$.

**PROBLEM 12:** Can a Smarandache strong P-groupoid in $Z(n)$ be not a Smarandache strong alternative groupoid?

**PROBLEM 13:** Find the number of Smarandache strong Moufang groupoids in $Z(n)$.

**PROBLEM 14:** Find the SGs in

1. $Z^*(n)$.
2. $Z^*(n) \setminus Z(n)$.
3. $Z^{***}(n) \setminus Z^{**}(n)$.

**PROBLEM 15:** Does $Z^*(p^m)$ where p is a prime $m > 1$ have atleast $p^m$ SGs in them.

*Hint:* For instance in $Z^*(2^3)$ that is in $Z^*(8)$, we have atleast 8 SGs.

**PROBLEM 16:** Let n be a composite number. Is all groupoids in $Z^*(n) \setminus Z(n)$ SGs?

**PROBLEM 17:** Which class will have more number of SGs in $Z^*(n) \setminus Z(n)$

1. When n + 1 is prime.
2. When n is prime.

3. When n is a non prime and n + 1 is also a non prime.

**PROBLEM 18:** How many groupoids in $Z^{**}(n) \setminus Z^*(n)$ are SGs barring $Z_n(1, 1)$ which is a semigroup?
Discuss the case

1. When n is prime.
2. When n is odd.
3. When n is even.

**PROBLEM 19:** Find all SGs in $Z^{**}(n)$.

1. How many are Smarandache strong Bol groupoids?
2. How many are Smarandache strong Moufang groupoids?
3. How many are Smarandache strong P-groupoids?
4. How many are Smarandache strong alternative groupoids?
5. How many are Smarandache strong idempotent groupoids?

**PROBLEM 20:** When will the number of SGs be greater when n is prime or when n is even?

**PROBLEM 21:** Let $Z_n = \{0, 1, 2, \ldots, n-1\}$, $n > 3$, $n < \infty$. If $n = p_1^{\alpha_1} \ldots p_m^{\alpha_m}$ where each $p_i$ is a prime. How many SGs does $Z^{***}(n) \setminus Z^{**}(n)$ contain?

**PROBLEM 22:** For any n find all SGs in $Z^{***}(n)$ which are got by $Z_n(m, 0)$ such that $m^2 \equiv m \pmod{n}$ or equivalently we term the problem as find for any general n the number of idempotents in $Z_n$.

**PROBLEM 23:** Find the number of Smarandache normal groupoids in $Z^{***}(n) \times Z^{***}(m) \times Z^{***}(p)$ where

1. n is composite.
2. m is even.
3. p is a prime.

**PROBLEM 24:** Find the number of Smarandache Bol groupoids in $Z^{***}(n) \times Z^{***}(n+1)$.

**PROBLEM 25:** Find the number of Smarandache Moufang groupoids in $G = Z^{***}(p-1) \times Z^{***}(p) \times Z^{***}(p+1)$ where p is any prime.

**PROBLEM 26:** Find the number of SGs in $G = Z^{***}(m) \times Z^{***}(m^2)$.

1. When m is a prime.
2. When m is even.
3. When m is odd.

**PROBLEM 27:** Find the number of SGs in the class G (n) defined in chapter five by adjoining e with $Z_n$ where $G(n) = \{(G(t, u) = Z_n \cup \{e\})\ (t, u)$ where $t \neq u$, $(t, u) = 1$; $t, u \in Z_n \setminus \{0\}\}$.

*Hint:* The operation $*$ on $G = Z_n \cup \{e\}$ is $a * a = e$ for all $a \in Z_n \cup \{e\}$ $a * e = e * a = a$ for all $a \in Z_n$ and $a * b = at + bu \pmod{n}$ where $a \neq b$; $a, b \in Z_n$.

$G(t, u)$ is a groupoid, when we take the collection of all such groupoids using $Z_n$ for varying t and u we denote it by $G(n)$.
Find

1. Number of loops in $G(n)$.
2. Number of Smarandache strong Bol groupoids and number of Smarandache Bol groupoids in $G(n)$.
3. Number of Smarandache strong Moufang groupoids and number of Smarandache Moufang groupoids in $G(n)$.
4. Number of Smarandache P-groupoids and number of Smarandache strong P-groupoids in $G(n)$.
5. Number of Smarandache strong alternative groupoid and number of Smarandache alternative groupoids in $G(n)$.

**PROBLEM 28:** Study the problem 27 for $G^*(n)$, $G^{**}(n)$ and $G^{***}(n)$.

**PROBLEM 29:** Can $G(n)$ have any loop of

1. Prime order.
2. Odd order.

*Hint:* We know a class of loops exists in $G^{**}(n)$ when n is odd so that even order loops exist.

**PROBLEM 30:** How many loops are in $G(n)$?

**PROBLEM 31:** Define Smarandache cascades using Smarandache automatons and describe them.

**PROBLEM 32:** Obtain any interesting results about Smarandache automaton.

# INDEX

*About the Author*

Dr. W. B. Vasantha is an Associate Professor in the Department of Mathematics, Indian Institute of Technology Madras, Chennai, where she lives with her husband Dr. K. Kandasamy and daughters Meena and Kama. Her current interests include Smarandache algebraic structures, fuzzy theory, coding/ communication theory. In the past decade she has guided seven Ph.D. scholars in the different fields of non-associative algebras, algebraic coding theory, transportation theory, fuzzy groups, and applications of fuzzy theory of the problems faced in chemical industries and cement industries. Currently, five Ph.D. scholars are working under her guidance. She has to her credit 250 research papers of which 218 are individually authored.  Apart from this she and her students have presented around 230 papers in national and international conferences. She teaches both undergraduate and post-graduate students and has guided over 35 M.Sc. and M.Tech projects. She has worked in collaboration projects with Indian Space Research Organization and with the Tamil Nadu State AIDS Control Society.

She can be contacted at wbvasantha@indiainfo.com and vasantak@md3.vsnl.net.in

113